\documentclass[a4paper,11pt,twoside]{article}
\usepackage{pst-fill,pst-grad}
\usepackage{textcomp}
\usepackage[english]{babel}
\usepackage[utf8x]{inputenc}
\usepackage{graphicx}
\usepackage{amsmath}
\usepackage{multicol}
\usepackage{float}
\usepackage{fancyhdr}
\usepackage[matrix,arrow,curve]{xy}
\usepackage{pstricks} 
\usepackage{amsmath,amsfonts,verbatim,afterpage,theorem,euscript,mathrsfs,amssymb}
\usepackage{amsfonts}
\usepackage{amssymb}
\usepackage{array}
\usepackage{dsfont}
\usepackage{hyperref}

\numberwithin{equation}{section}
\def \vu{\vec{u}}
\def \vv{\vec{v}}

\def \vU{\vec{U}}
\def \vB{\vec{B}}
 
\def \vUU{\vec{\mathcal{U}}}

\def \vn{\vec{\nabla}}
\def \vb{\vec{b}}
\def \vf{\vec{f}}
\def \vg{\vec{g}}
\def \vF{\vec{F}}
\def \vG{\vec{G}}

\def \vphi{\vec{\varphi}}
\def \vPhi{\vec{\Phi}}

\def \vpphi{\vec{\phi}}

\def \R{\mathbb{R}^{3}}
 
\def \q{q_{0}}
\def \vVV{\vec{\mathcal{V}}}
\def \NN{\vec{\mathcal{N}}}

\def \pj{\partial_j}
\def \pai{\partial_i}

\def \Q1{\mathcal{Q}_{1}}
\def \B1{\mathcal{B}_{1}}

\def \uuo{(|\vu|^2)}
\def \bbo{(|\vb|^2)}
\newcommand{\norm}[1]{\lVert#1\rVert}
\def \Ar{\mathbf{A}_{r}}
\def \BBr{\mathbf{B}_{r}}

\def \Dr{\mathbf{D}_{r}}
\def \Prq{\mathbf{P}_{r}}
\def \Ao{\mathbf{A}_{\rho}}
\def \BBo{\mathbf{B}_{\rho}}

\def \Do{\mathbf{D}_{\rho}}
\def \Poq{\mathbf{P}_{\rho}}
\def \Tr{\mathbf{\Theta}_{r}}
\def \TTo{\mathbf{\Theta}_{\rho}}

\def \Qrq{\mathbf{Q}_{r}}
\def \ftau{\frac{5}{\tau_{0}}}

\newtheorem{Definition}{Definition}[section]
\newtheorem{Proposition}{Proposition}[section]
\newtheorem{Lemme}{Lemma}[section]
\newtheorem{Theoreme}{Theorem}
\newtheorem{Corollaire}{Corollary}[section]
\newtheorem{Remarque}{Remark}[section]

\title{\bf  On the partial regularity theory for the MHD equations}
\usepackage{authblk}
\author{Diego Chamorro\footnote{\emph{diego.chamorro@univ-evry.fr}} }
\author{Jiao He\footnote{\emph{jiao.he@univ-evry.fr}}}
\affil{\footnotesize LaMME, Univ. Evry, CNRS, Universit\'e Paris-Saclay, 91025, Evry, France.}

\begin{document}
\maketitle
\begin{scriptsize}
\abstract{We generalize here the celebrated Partial Regularity Theory of Caffarelli, Kohn and Nirenberg to the MHD equations in the framework of parabolic Morrey spaces. This type of parabolic generalization using Morrey spaces appears to be crucial when studying the role of the pressure in the regularity theory for the classical Navier-Stokes equations as well as for the MHD equations.}\\

\end{scriptsize}

\section{Introduction}
In this article we study \emph{regularity} results for the incompressible 3D magnetohydrodynamic (MHD) equations which are given by the following system:
\begin{equation}\label{EquationMHDoriginal}
\begin{cases}
\partial_{t}\vU=\Delta \vU -(\vU\cdot\vn)\vU+(\vB\cdot\vn)\vB-\vn P+\vF,\quad div(\vU) = div(\vF)=0,\\[3mm]
\partial_{t}\vB=\Delta \vB -(\vU\cdot\vn)\vB+(\vB\cdot\vn)\vU+\vG,\quad div(\vB)=div(\vG)=0,\\[3mm]
\vU(0,x)=\vU_{0}(x), div(\vU_0)=0 \mbox{ and } \vB(0,x)=\vB_{0}(x), div(\vB_0)=0,\qquad x\in \mathbb{R}^3,
\end{cases}
\end{equation}
where $\vU, \vB:[0,T]\times \R\longrightarrow \R$ are two divergence-free vector fields which represent the velocity and the magnetic field, respectively, and the scalar function $P:[0,T]\times \R\longrightarrow \mathbb{R}$ stands for the pressure. The initial data $\vU_0, \vB_0: \R\longrightarrow \R$ and the external forces $\vF, \vG:[0,T]\times \R\longrightarrow \R$ are given.\\

Of course, when the magnetic field $\vB$ becomes the zero vector, the MHD equations \eqref{EquationMHDoriginal} are reduced to the 3D classical Navier-Stokes equations 
\begin{equation}\label{NSequations}
\partial_{t}\vU=\Delta \vU -(\vU\cdot\vn)\vU-\vn P+\vF, \qquad div(\vU)=div(\vF)=0.
\end{equation}
It is worth noting here that for the Navier-Stokes equations there are two different regularity theories. The first one, known as the Serrin \emph{local} theory \cite{Serrin1}, is essentially based on a control of the velocity vector field of the type $\vU\in (L^p_tL^q_x)_{loc}$ with $\tfrac{2}{p}+\tfrac{3}{q}\leq1$ (the case $\tfrac{2}{p}+\tfrac{3}{q}=1$ was proved by Struwe \cite{Struwe} and Takahashi \cite{Takahashi}) and with this assumption it is possible to obtain a local gain of regularity of the solutions of (\ref{NSequations}). One very important feature of this theory is the fact that no particular restrictions are asked to the pressure $P$ which can be a very general object (for example we can ask $P\in \mathcal{D}'$). However, this generality implies paradoxically some constraints and the gain of regularity is only obtained in the spatial variable as the temporal regularity is linked to some information on the pressure (see Section 13.1 of \cite{PGLR1} for this particular point). \\

The second regularity theory,  known as the \emph{partial regularity theory},  is due to Caffarelli, Kohn and Nirenberg and it was developed in \cite{CKN}. In this case the local boundedness assumption is replaced by local energy estimates and with some additional hypothesis on the pressure $P$ (usually $P\in (L^{q_{0}}_{t,x})_{loc}$ for some $1<q_0<+\infty$) we can deduce a gain of regularity in \emph{both} variables, space and time. \\

These two points of view are of course quite different since they rely on different techniques and require different hypotheses. However  it is important to point out that a common treatment of these two theories can be performed by using the framework of parabolic Morrey spaces $\mathcal{M}_{t,x}^{p,q}$ (see formula (\ref{Morreyparabolic}) below for a precise definition of these functional spaces). Indeed, O'Leary \cite{OLeary} generalized Serrin's theory by replacing the local Lebesgue hypothesis with a local information expressed in terms of parabolic Morrey spaces, while Kukavica \cite{Kukavica} proposed a generalization of Caffarelli-Kohn-Nirenberg's theory using this parabolic framework. An interesting point of this common framework appears clearly when studying the role of the pressure in the Caffarelli-Kohn-Nirenberg theory for the classical Navier-Stokes equations, indeed, as it is shown in \cite{CML}, the language of parabolic Morrey spaces is a powerful tool which allows to mix, in a very specific sense, these two regularity theories.\\

In a recent article \cite{ChCHJ}, we have generalized to the MHD equations (\ref{EquationMHDoriginal}) the local regularity theory using parabolic Morrey spaces. The aim of this article is now to generalize these techniques in order to study the partial regularity theory for the MHD equations and we will deduce here parabolic H\"older regularity for the solutions of theses equations in small neighborhoods (see Theorem \ref{Teorem1} for a precise statement).\\

Note that the Caffarelli-Kohn-Nirenberg theory has been investigated for the MHD equations (see \cite{He1}, \cite{He2}), but to the best of our knowledge the generalization using parabolic Morrey spaces is new and we find this approach interesting since this framework admits some important applications.\\

The plan of the article is as follows: in Section \ref{Secc_NotationandResults} we introduce some notation and we state the main theorem. In Section \ref{Secc_ProofMainTheorem} we present, under some particular assumptions, the general strategy for proving Theorem \ref{Teorem1} while in Sections \ref{Secc_AveragedBounds},  \ref{Secc_InductiveArgument} and \ref{Secc_LastHypotheses} we deduce these assumptions from the general hypothesis of the partial regularity theory. Finally, in Appendix \ref{Secc_AppendixA}  we recall some useful properties of Morrey spaces and in Appendix \ref{Secc_AppendixB} we give the proof of a technical lemma. 
\section{Notation and presentation of the results}\label{Secc_NotationandResults}
The starting point of this work relies on the use of the Elsasser formulation  for the MHD equations (see \cite{elsasser1950}) which enables us to obtain a more symmetric expression of the problem considered here \footnote{This approach is interesting since we will assume the \emph{same} hypotheses on the variables $\vU$ and $\vB$.}. More precisely, if we define $\vu= \vU + \vB$, $\vb= \vU-\vB$, $\vf= \vF + \vG$ and $\vg = \vF - \vG$, then the original system  \eqref{EquationMHDoriginal} becomes  
\begin{equation}\label{EquationMHD}
\begin{cases}
\partial_{t}\vu=\Delta \vu -(\vb\cdot\vn)\vu-\vn P+\vf,\quad div (\vu) = div (\vf)=0,\\[3mm]
\partial_{t}\vb=\Delta \vb -(\vu\cdot\vn)\vb- \vn P +\vg,\quad div (\vb)=div (\vg)=0,\\[3mm]
\vu(0,x)=\vu_{0}(x), div(\vu_0)=0,\quad \vb(0,x)=\vb_{0}(x), div(\vb_0)=0,
\end{cases}
\end{equation}
where, since $div(\vu)=div(\vb)=0$, we have that $P$ satisfies the equation 
\begin{equation}\label{EquationPression}
\Delta P= - \sum^3_{i,j= 1} \partial_i \partial_j (u_i b_j),
\end{equation}
and from this equation, we remark that the pressure $P$ is only determined by the couple  $(\vu, \vb)$. We will see that in Section \ref{Secc_AveragedBounds} this equation is crucial when doing pressure estimates. Remark in particular that the solution $(\vu, \vb)$ to \eqref{EquationMHD} has the same regularity as the solution $(\vU, \vB)$ to \eqref{EquationMHDoriginal}, so all the rest of the article will be devoted to study the regularity of the couple $(\vu, \vb)$. \\

Now, let $\Omega$ be a bounded domain of $]0, +\infty[\times \mathbb{R}^3$, we assume the following (local) hypotheses:
\begin{equation}\label{HypothesesLocal1}
\begin{split}
&\vu, \vb \in L^{\infty}_t L^{2}_x \cap L^{2}_t \dot{H}^{1}_x(\Omega),\\[3mm]
&P \in L^{q_{0}}_{t,x}(\Omega) \mbox{ with } 1<q_0\leq \tfrac{3}{2},\\[3mm]
&\vf, \vg \in L^{\frac{10}{7}}_{t,x}(\Omega).
\end{split}
\end{equation}
\begin{Remarque}\label{RemarquePression1}
Without loss of generality, we can assume that $P \in L^{\frac{3}{2}}_{t,x}(\Omega)$, as this sole condition implies the hypotheses  $P \in L^{q_{0}}_{t,x}(\Omega)$ in the whole range $1\leq q_0\leq \tfrac{3}{2}$.
\end{Remarque}
We will say that the couple $(\vu, \vb)\in L^\infty_tL^2_x\cap L^2_t\dot{H}^1_x(\Omega)$ satisfies the MHD equations (\ref{EquationMHD}) in the weak sense if for all $\vphi, \vpphi \in \mathcal{D}(\Omega)$ such that $div(\vphi)=div(\vpphi)=0$, we have\\
$$
\begin{cases}
\langle\partial_{t}\vu-\Delta \vu +(\vb\cdot\vn)\vu-\vf|\vphi \rangle_{\mathcal{D}'\times \mathcal{D}}=0,\\[3mm]
\langle\partial_{t}\vb-\Delta \vb +(\vu\cdot\vn)\vb-\vg|\vpphi \rangle_{\mathcal{D}'\times \mathcal{D}}=0,
\end{cases}
$$
note that if $(\vu, \vb)$ are solutions of the previous system, then due to the expression (\ref{EquationPression}) there exists a pressure $P$ such that (\ref{EquationMHD}) is fulfilled in $\mathcal{D}'$.\\

The class of weak solutions is too wide for our purposes and we need to reduce the set of admissible solutions and actually we will only work with a very specific subset given by the following definition. 

\begin{Definition}[Suitable solution]\label{Def_SuitableSolutions} Let $(\vu, P, \vb)$ be a weak solution over $\Omega$ of equations (\ref{EquationMHD}). We will say that the $(\vu, P, \vb)$ is a \emph{suitable} solution if the distribution $\mu$ given by the expression
\begin{eqnarray*}
\mu&=&-\partial_t(|\vu|^2 + |\vb|^2 )+ \Delta (|\vu|^2 + |\vb|^2 )  - 2 (|\vn \otimes \vu|^2 + |\vn \otimes \vb|^2 )\\  
& &- div \left( (|\vu|^2 + 2P )\vb + (|\vb|^2 + 2P )\vu \right) + 2 (\vf \cdot \vu + \vg \cdot\vb),
\end{eqnarray*}
is a non-negative locally finite measure on $\Omega$.
\end{Definition}
It is worth noting here that from the set of hypotheses (\ref{HypothesesLocal1}) we can deduce that $\mu$ is well defined as a distribution but we will need to assume its positivity, which is the whole point of \emph{suitable} solutions.\\

We still need to introduce one more ingredient which is related to the parabolic structure of the functional spaces we are going to work with. We  consider the homogeneous space $(\mathbb{R}\times \R, d, \lambda)$ where $d$ is the parabolic quasi-distance given by 
\begin{equation}\label{Def_QuasiDistance}
d\big((t,x), (s,y)\big)=|t-s|^{\frac{1}{2}}+|x-y|,
\end{equation}
and where $\lambda$ is the usual Lebesgue measure $d\lambda=dtdx$. Remark that the homogeneous dimension is now $N=5$. See \cite{Folland} for more details concerning the general theory of homogeneous spaces. Associated to this distance, we can define homogeneous (parabolic) H\"older spaces $\dot{\mathcal{C}}^\alpha(\mathbb{R}\times \R, \R)$ where $\alpha\in ]0,1[$ by the following condition:
\begin{equation}\label{Holderparabolic}
\|\vphi\|_{\dot{\mathcal{C}}^\alpha}=\underset{(t,x)\neq (s,y)}{\sup}\frac{|\vphi(t,x)-\vphi(s,y)|}{\left(|t-s|^{\frac{1}{2}}+|x-y|\right)^\alpha}<+\infty,
\end{equation}
as we can see this quantity captures H\"older regularity in both time and space variables. 

Now, for $1< p\leq q<+\infty$, we define the parabolic Morrey spaces $\mathcal{M}_{t,x}^{p,q}$ as the set of measurable functions $\vphi:\mathbb{R}\times\R\longrightarrow \R$ that belong to the space $(L^p_{t,x})_{loc}$ such that $\|\vphi\|_{\mathcal{M}_{t,x}^{p,q}}<+\infty$ where
\begin{equation}\label{Morreyparabolic}
\|\vphi\|_{\mathcal{M}_{t,x}^{p,q}}=\underset{(t_{0}, x_{0}) \in \mathbb{R}\times \R  , r>0}{\sup}\left(\frac{1}{r^{5(1-\frac{p}{q})}}\int_{|t-t_{0}|<r^{2}}\int_{B(x_{0},r)}|\vphi(t,x)|^{p}dxdt\right)^{\frac{1}{p}}.
\end{equation}
Remark in particular that we have $\mathcal{M}_{t,x}^{p,p}=L^p_{t,x}$ and we will list  some useful properties of these spaces in Appendix \ref{Secc_AppendixA}.\\ 

These parabolic spaces are very useful in the analysis of the properties of the solutions of the Navier-Stokes equations, hence of the MHD equations and their properties appears to be more and more useful in the study of some PDEs. See for example \cite{CML}, \cite{ChCHJ}, \cite{Chen},  \cite{Kukavica},  \cite{OLeary}, \cite{Robinson} and the book \cite{PGLR1}.\\

We can now state our main theorem which studies the H\"older regularity of suitable solutions of the MHD equations (\ref{EquationMHD}). 
\begin{Theoreme}\label{Teorem1}
Let $\Omega$ be a bounded domain of $]0, +\infty[ \times \mathbb{R}^{3}$. Let $(\vu, P, \vb)$ be a weak solution on $\Omega$ of the MHD equations \eqref{EquationMHD}.
Assume that
\begin{itemize}
\item[1)] $(\vu, \vb, P, \vf, \vg)$ satisfies the conditions (\ref{HypothesesLocal1}),
\item[2)] $(\vu, P, \vb)$ is suitable in the sense of Definition \ref{Def_SuitableSolutions}, 
\item[3)] we have the following local information on $\vf$ and $\vg$: $\mathds{1}_{\Omega}\vf \in \mathcal{M}_{t,x}^{\frac{10}{7}, \tau_{a}}$ and $ \mathds{1}_{\Omega} \vg \in \mathcal{M}_{t,x}^{\frac{10}{7}, \tau_{b}}$ for some $\tau_{a}, \tau_b>\frac{5}{2-\alpha}$ with $0<\alpha<\frac{1}{3}$.
\end{itemize}
There exists a positive constant $\epsilon^{*}$ which depends only on $\tau_{a}$ and $\tau_b$ such that, if for some $\left(t_{0}, x_{0}\right) \in \Omega$, we have
\begin{equation}\label{HypothesePetitesseGrad}
\limsup _{r \rightarrow 0} \frac{1}{r} \iint_{]t_{0}-r^{2}, t_{0}+r^{2}[ \times B\left(x_{0}, r\right)}|\vn \otimes \vu|^{2} + |\vn \otimes \vb|^{2} dxds <\epsilon^{*},
\end{equation}
then $(\vu, \vb)$ is H\"olderian of exponent $\alpha$ in a neighborhood of $\left(t_{0}, x_{0}\right)$ (in the sense of (\ref{Holderparabolic})) for some small $\alpha$ in the interval $0<\alpha<\frac{1}{3}$.
\end{Theoreme}
Some remarks are in order here. First, since we are assuming some control in the pressure (see Remark \ref{RemarquePression1} above) we can obtain regularity results in time and space variables which is expressed here in terms of parabolic H\"older spaces. Remark also that the parameters $\tau_a$ and $\tau_b$ that define the Morrey spaces for the forces $\vf$ and $\vg$ are linked to the exponent $\alpha$ of the expected H\"olderian regularity and this is somehow natural as the information given by the external forces is not involved in the nonlinear terms and must be taken into account. Note now that since $\frac{10}{7}\leq \frac{5}{2-\alpha}<\tau_a,\tau_b$, then by Lemma \ref{Lemme_locindi} the local conditions $\mathds{1}_{\Omega}\vf \in \mathcal{M}_{t,x}^{\frac{10}{7}, \tau_{a}}$ and $ \mathds{1}_{\Omega} \vg \in \mathcal{M}_{t,x}^{\frac{10}{7}, \tau_{b}}$ are stronger than the conditions $\vf, \vg \in L^{\frac{10}{7}}_{t,x} (\Omega)$ given in  (\ref{HypothesesLocal1}) and this fact explains the third hypothesis above. Finally, we note that the regularity obtained is only valid on small neighborhoods of points for which we have (\ref{HypothesePetitesseGrad}), 
and this justifies the designation of \emph{partial regularity} associated to this theory. \\

To end this section, let us remark that the set $\Sigma_0$ of points for which we have 
$$\limsup _{r \rightarrow 0} \frac{1}{r} \iint_{]t_{0}-r^{2}, t_{0}+r^{2}[ \times B\left(x_{0}, r\right)}|\vn \otimes \vu|^{2} + |\vn \otimes \vb|^{2}   dx ds>0,$$
is called the set of \emph{large gradients} and from Theorem \ref{Teorem1} and a standard Vitali covering Lemma it can be deduced (see Section 13.10 of the book \cite{PGLR1}) that the parabolic Hausdorff measure of the set $\Sigma_0$ is null, which means that this set is actually very small. 

\section{Proof of the main Theorem}\label{Secc_ProofMainTheorem}

The strategy of the proof of Theorem \ref{Teorem1} is based on regularity results of solutions of parabolic equations. Indeed, following a classical result given in the book \cite{Ladyzhenskaya1} (see also the article \cite{Ladyzhenskaya2}) we have the following lemma (stated using parabolic Morrey spaces and borrowed from Proposition 13.4 of  the book \cite{PGLR1}):
\begin{Lemme}[H\"older regularity]\label{Lemma_holderian}
For $\vv, \vPhi:[0, +\infty[\times \mathbb{R}^3\longrightarrow \mathbb{R}^3$ two vector fields, we consider the following equation
\begin{equation}\label{ParabolicEquation}
\begin{cases}
\partial_t\vv(t,x)=\Delta \vv(t,x)+\vPhi(t,x),\\[3mm]
\vv(0, x)=0.
\end{cases}
\end{equation}
Assume moreover that we have the information $\vPhi \in \mathcal{M}_{t,x}^{\mathfrak{p}_0, \mathfrak{q}_{0}}$ with $1 \leq \mathfrak{p}_0 \leq \mathfrak{q}_0$, $\frac{5}{2}<\mathfrak{q}_0<3$ and $\frac{1}{\mathfrak{q}_0} = \frac{2-\alpha}{5}$, $0< \alpha < \frac{1}{3}$. Then the function $\vv$ equal to $0$ for $t\leq 0$ and to 
$$\vv (t,x) = \int_0^t e^{(t-s) \Delta}\vPhi(s, \cdot) \,ds,$$
for $t > 0$, is a solution of equation (\ref{ParabolicEquation}) that is Hölderian of exponent $\alpha$ with respect to the parabolic distance (\ref{Def_QuasiDistance}). 
\end{Lemme}
Of course, as we only have local hypotheses, we can not apply directly this lemma to the MHD equations (\ref{EquationMHD}) and the first step is to localize our framework: we fix the point $(t_0, x_0)$ considered in the hypotheses of Theorem \ref{Teorem1} and we construct two auxiliary non-negative functions $\varphi, \psi:\mathbb{R}\times \mathbb{R}^3\longrightarrow \mathbb{R}$ by the conditions $\varphi,\psi\in \mathcal{C}^{\infty}_0(\mathbb{R}\times \mathbb{R}^3, \mathbb{R})$
\begin{equation*}
\begin{split}
& supp(\varphi)\subset ]-\tfrac{1}{16},\tfrac{1}{16}[\times B(0,\tfrac14),\\[3mm]
& supp(\psi)\subset ]-\tfrac{1}{4},\tfrac{1}{4}[\times B(0,\tfrac12),
\end{split}
\end{equation*}
and such that 
\begin{equation}\label{Definition_FonctionLocalisantes}
\begin{split}
\varphi\equiv1 \mbox{ on } ]-\tfrac{1}{64}, \tfrac{1}{64}[\times B(0, \tfrac18),\\[3mm]
\psi\equiv1 \mbox{ on } ]-\tfrac{1}{16}, \tfrac{1}{16}[\times B(0, \tfrac14).
\end{split}
\end{equation}
Remark in particular that $\psi\equiv 1$ on the support of $\varphi$ and thus we have the pointwise identity $\psi \varphi=\varphi$ in $\mathbb{R}\times \mathbb{R}^3$. Now, for a small fixed $R_0$ such that $0<4R_0<t_0$ we define 
\begin{equation}\label{Def_VariableU}
\phi(t,x)=\varphi\left(\frac{t-t_0}{R_0^2}, \frac{x-x_0}{R_0}\right)\qquad \mbox{and}\quad \vUU= \phi (\vu+ \vb),
\end{equation}
as we can observe, we have the identity $\vUU=\vu+\vb$ on a small neighborhood of the point $(t_0, x_0)$ and the support of the variable $\vUU$ is contained in the parabolic ball of the form
\begin{equation}\label{Def_ParabolicBalldeBase}
Q_{R_0}(t_0,x_0)=]t_0-R_0^2,t_0+R_0^2[\times B(x_0,R_0).
\end{equation}
When the context is clear, we will write $Q_{R_0}$ instead of $Q_{R_0}(t_0,x_0)$ and for usual (euclidean) balls we will write $B_{R_0}$ instead of $B(x_0,R_0)$. We will assume moreover that $R_0$ is small enough to grant that 
\begin{equation}\label{Defi_BouleParaboliquedeBase}
Q_{4R_0}(t_0, x_0)\subset \Omega. 
\end{equation}
Note now that, since $0<4R_0<t_0$  and $supp(\phi)\subset ]t_0-\tfrac{R_0^2}{16}, t_0 +\tfrac{R_0^2}{16}[\times B(x_0,\tfrac{R_0}{4})$, we have $\vUU(0, x)=0$ and we obtain the following equation
\begin{equation}\label{Equation_ChaleurLocale}
\begin{cases}
\partial_t\vUU(t,x)=\Delta \vUU(t,x)+\vPhi(t,x),\\[3mm]
\vUU(0, x)=0,
\end{cases}
\end{equation}
where 
\begin{equation}\label{TermeGeneral_a_Estimer}
\vPhi=\underbrace{(\partial_t \phi-\Delta \phi) (\vu+\vb)}_{(a)} -2\sum_{i=1}^{3}\underbrace{(\partial_{i}\phi) (\partial_i(\vu+\vb))}_{(b)} - \underbrace{\phi \left( (\vb\cdot\vn)\vu+(\vu\cdot\vn)\vb\right)}_{(c)}- 2 \underbrace{\phi(\vn P)}_{(d)} +\underbrace{\phi (\vf+\vg)}_{(e)}. 
\end{equation}
Thus, in order to apply Lemma \ref{Lemma_holderian}, we only need to proof that the function $\vPhi$ belongs to the Morrey space $\mathcal{M}_{t,x}^{\mathfrak{p}_0, \mathfrak{q}_{0}}$ with $1 \leq \mathfrak{p}_0 \leq \mathfrak{q}_0$, $\frac{5}{2}< \mathfrak{q}_0<3$ and $\frac{1}{\mathfrak{q}_0} = \frac{2-\alpha}{5}$, $0< \alpha <  \frac{1}{3}$: then due to the formulas (\ref{Definition_FonctionLocalisantes}) and (\ref{Def_VariableU}),  it is straightforward to deduce that the function $\vu+\vb$ is H\"older regular of order $\alpha$ on a small neighborhood of $(t_0, x_0)$ contained in the parabolic ball $Q_{R_0}$.
\begin{Remarque}
Since the hypotheses on $\vu$ and $\vb$ are completely symmetric, it is possible to perform a separated study of $\vu$ and $\vb$ in order to obtain the H\"older regularity for each one of these variables. As the computations are exactly the same, for the sake of simplicity, we prefer to study the function $\vu+\vb$. 

\end{Remarque}

The fact that $\vPhi\in\mathcal{M}_{t,x}^{\mathfrak{p}_0, \mathfrak{q}_{0}}$ will made possible as long as we have some interesting estimates of the constitutive terms of (\ref{TermeGeneral_a_Estimer}). In this sense we have the following proposition:
\begin{Proposition}\label{Proposition_Principale1}
Let $R_1, R_2, R_3$ be three real numbers such that 
$$0<R_0<R_1<R_2<R_3<2R_0<t_0,$$
and consider $(\vu, P,\vb)$ a suitable solution of MHD equations (\ref{EquationMHD}) over $\Omega$ in the sense of Definition \ref{Def_SuitableSolutions}. In the framework of the general assumptions of Theorem \ref{Teorem1}, assume moreover that on some parabolic neighborhood $Q_{R_3}(t_0, x_0)$, $Q_{R_2}(t_0, x_0)$ and $Q_{R_1}(t_0, x_0)$ of type (\ref{Def_ParabolicBalldeBase}) we have the following information:
\begin{itemize}
\item[1)] $\mathds{1}_{Q_{R_3}}\vu,\; \mathds{1}_{Q_{R_3}} \vb\in \mathcal{M}_{t,x}^{3, \tau_{0}}$ for some $\tau_{0} > \frac{5}{1-\alpha}$,
\item[2)] $\mathds{1}_{Q_{R_3}} \vn \otimes \vu, \;\mathds{1}_{Q_{R_3}}\vn \otimes \vb \in \mathcal{M}_{t,x}^{2, \tau_{1}}$ with $\frac{1}{\tau_{1}} = \frac{1}{\tau_{0}} + \frac{1}{5}$,
\item[3)] $\mathds{1}_{Q_{R_2}}\vu, \;\mathds{1}_{Q_{R_2}} \vb  \in \mathcal{M}_{t,x}^{3, \delta}$ with $\frac{1}{\delta} + \frac{1}{\tau_{0}} \leq \frac{1-\alpha}{5}$,
\item[4)] for all $1\leq i,j\leq 3$ we have $\mathds{1}_{Q_{R_1}} \frac{\vn\partial_i \partial_j}{(-\Delta)}(u_ib_j)\in\mathcal{M}_{t,x}^{\mathfrak{p}, \mathfrak{q}}$ with $\mathfrak{p}_0\leq \mathfrak{p}<+\infty$ and $\mathfrak{q}_0\leq \mathfrak{q}<+\infty$,
\item[5)] $\mathds{1}_{Q_{R_3}} \vf\in \mathcal{M}_{t,x}^{\frac{10}{7},\tau_{a}}$ and $\mathds{1}_{Q_{R_3}}\vg \in \mathcal{M}_{t,x}^{\frac{10}{7},\tau_{b}}$ for some $\tau_a,\tau_b > \frac{5}{2-\alpha}$,
\end{itemize}
then we have that all the terms of (\ref{TermeGeneral_a_Estimer}), and therefore the function $\vPhi$ itself, belong to the Morrey space $\mathcal{M}_{t,x}^{\mathfrak{p}_0, \mathfrak{q}_0}$ with $1 \leq \mathfrak{p}_0\leq \tfrac{6}{5}$ and $\tfrac52 <\mathfrak{q}_0<3$ where $\frac{1}{\mathfrak{q}_0} = \frac{2-\alpha}{5}$ with $0< \alpha < \frac{1}{3}$.
\end{Proposition}
\begin{Remarque}\label{Remark_Point2}
Note that Theorem \ref{Teorem1} follows at once if we have the conclusion of this proposition: we only need to apply Lemma \ref{Lemma_holderian} to obtain that the function $\vUU$ defined in (\ref{Def_VariableU}) is H\"olderian of exponent $\alpha$ and since the information over $\vu$ and $\vb$ is symmetric, it is easy to obtain that the couple $(\vu, \vb)$ is itself H\"olderian of exponent $\alpha$.
\end{Remarque}
\begin{Remarque}
The upper bound $1\leq \mathfrak{p}_0\leq \frac{6}{5}$ given in Proposition \ref{Proposition_Principale1} is technical and ensures the condition $\mathfrak{p}_0\leq \mathfrak{q}_0$. Note in particular that in Lemma \ref{Lemma_holderian} the H\"older regularity exponent $0<\alpha< \frac{1}{3}$ is only related to the parameter $\mathfrak{q}_0$ and not to $\mathfrak{p}_0$.
\end{Remarque}
\begin{Remarque}
It is important to mention here that the term $ \frac{\vn\partial_i \partial_j}{(-\Delta)}(u_ib_j)$ in the hypothesis 4) is related to the pressure term $\vn P$ in \eqref{TermeGeneral_a_Estimer}. Note however that in Proposition \ref{Proposition_Principale1} we do not state any particular assumption on the pressure $P$ but, as we will see later on, in order to obtain the hypotheses 1)-5) of this proposition we will need the information $P \in L^{q_{0}}_{t,x}(\Omega)$ with $1<q_0\leq\frac{3}{2}$ as stated in the general framework (\ref{HypothesesLocal1}).
\end{Remarque}
\begin{Remarque}\label{Remark_Point5} The points 3) and 4) of the hypotheses will be deduced later on from the points 1), 2) and 5) and this explains the fact that we need to reduce the support of the information as some extra localization properties are needed here. See Section \ref{Secc_LastHypotheses} for more details.
\end{Remarque}
{\bf Proof of Proposition \ref{Proposition_Principale1}.} Assuming for the moment the information stated in the points \emph{1)-5)} we will study each term (a)-(e) of (\ref{TermeGeneral_a_Estimer}) separately. 
\begin{itemize}
\item[(a)] Since we have by the point \emph{1)} the information $\mathds{1}_{Q_{R_3}}\vu,\; \mathds{1}_{Q_{R_3}} \vb\in \mathcal{M}_{t,x}^{3, \tau_{0}}$ for some $\tau_{0} > 5$, then it is easy to obtain that $(\partial_t \phi - \Delta \phi)(\vu + \vb)\in \mathcal{M}_{t,x}^{\mathfrak{p}_0, \mathfrak{q}_{0}}$. Indeed, since $\phi$ is a smooth function, then due to its support properties (see (\ref{Def_VariableU})), from the first point of Lemma \ref{Lemme_Product}, from Lemma \ref{Lemme_locindi} and since $1<\mathfrak{p}_{0}\leq \frac{6}{5}$, $\frac{5}{2}<\mathfrak{q}_0<3$, we have
$$\left\|(\partial_t \phi - \Delta \phi)  (\vu + \vb)\right\|_{\mathcal{M}_{t,x}^{\mathfrak{p}_0, \mathfrak{q}_{0}}}\leq C\|\mathds{1}_{Q_{R_0}}(\vu + \vb)\|_{\mathcal{M}_{t,x}^{\mathfrak{p}_0, \mathfrak{q}_{0}}}\leq C\|\mathds{1}_{Q_{R_3}}(\vu + \vb)\|_{\mathcal{M}_{t,x}^{3, \tau_{0}}}<+\infty.$$
\item [(b)] For the second term of (\ref{TermeGeneral_a_Estimer}) we use the information given by the point \emph{3)} of the hypotheses of Proposition \ref{Proposition_Principale1}. Thus, by the H\"older inequalities in Morrey spaces (see Lemma \ref{Lemme_Product}) we obtain
$$\left\|(\partial_{i}\phi) (\partial_i(\vu+\vb))\right\|_{\mathcal{M}_{t,x}^{\mathfrak{p}_0, \mathfrak{q}_0}}\leq \|\mathds{1}_{Q_{R_0}}\partial_{i}\phi\|_{\mathcal{M}_{t,x}^{p_1, q_1}} \left(\|\mathds{1}_{Q_{R_0}}\partial_i\vu\|_{\mathcal{M}_{t,x}^{2, q_2}}+\|\mathds{1}_{Q_{R_0}} \partial_i\vb\|_{\mathcal{M}_{t,x}^{2, q_2}}\right),$$
where $\tfrac{1}{p_1}+\tfrac{1}{2}\leq \tfrac{1}{\mathfrak{p}_0}$ and $\tfrac{1}{q_1}+\tfrac{1}{q_2}= \tfrac{1}{\mathfrak{q}_0}$, moreover, by Lemma \ref{Lemme_locindi} we have for $\tau_{1}\geq q_2$:
$$\left\|(\partial_{i}\phi) (\partial_i(\vu+\vb))\right\|_{\mathcal{M}_{t,x}^{\mathfrak{p}_0, \mathfrak{q}_0}}\leq C\left(\|\mathds{1}_{Q_{R_3}}\vn\otimes\vu\|_{\mathcal{M}_{t,x}^{2, \tau_{1}}}+\|\mathds{1}_{Q_{R_3}} \vn\otimes\vb\|_{\mathcal{M}_{t,x}^{2, \tau_{1}}}\right)<+\infty.$$
Note that since $1< \mathfrak{p}_0\leq \tfrac{6}{5}$, the condition $p_1\geq 3$ is enough to satisfy $\tfrac{1}{p_1}+\tfrac{1}{2}\leq \tfrac{1}{\mathfrak{p}_0}$. On the other hand, since $\frac{1}{\mathfrak{q}_0}=\tfrac{2-\alpha}{5}$ and $\frac{1}{\tau_{1}} = \frac{1}{\tau_{0}} + \frac{1}{5}$ we should have $\tfrac{1}{q_1}= \tfrac{2-\alpha}{5} - \tfrac{1}{q_2} \leq \tfrac{2-\alpha}{5} - \frac{1}{\tau_{1}} = \tfrac{1-\alpha}{5} - \frac{1}{\tau_{0}}$, thus since $\tau_{0} > \tfrac{5}{1-\alpha}$, this is possible as long as $q_1$ is big enough.
\item [(c)] We study the term $\left\|\phi \left( (\vb\cdot\vn)\vu+(\vu\cdot\vn)\vb\right)\right\|_{\mathcal{M}_{t,x}^{\mathfrak{p}_0, \mathfrak{q}_0}}$. Since $1< \mathfrak{p}_0\leq \tfrac{6}{5}$ and $\tfrac{5}{2}<\mathfrak{q}_0<3$, by Lemma \ref{Lemme_locindi}, by the H\"older inequalities in Morrey spaces and using the information of points \emph{2)-3)}, we have:
\begin{eqnarray*}
\left\|\phi \left( (\vb\cdot\vn)\vu+(\vu\cdot\vn)\vb\right)\right\|_{\mathcal{M}_{t,x}^{\mathfrak{p}_0, \mathfrak{q}_0}}& \leq & C\left\|\mathds{1}_{Q_{R_0}} \left( (\vb\cdot\vn)\vu+(\vu\cdot\vn)\vb\right)\right\|_{\mathcal{M}_{t,x}^{\frac{6}{5}, \mathfrak{q}_0}}\notag\\
& \leq &C\left( \|\mathds{1}_{Q_{R_2}}\vb\|_{\mathcal{M}_{t,x}^{3, \delta}}\|\mathds{1}_{Q_{R_3}}\vn\otimes\vu\|_{\mathcal{M}_{t,x}^{2, \tau_{1}}}\right.\\
& &\left. +\|\mathds{1}_{Q_{R_2}}\vu\|_{\mathcal{M}_{t,x}^{3, \delta}}\|\mathds{1}_{Q_{R_3}}\vn\otimes\vb\|_{\mathcal{M}_{t,x}^{2, \tau_{1}}}\right)<+\infty,\notag
\end{eqnarray*}
where we have $\frac{1}{\delta}+\frac{1}{\tau_{1}}\leq \frac{1}{\mathfrak{q}_0}$, but since $\frac{1}{\mathfrak{q}_0}=\frac{2-\alpha}{5}$ and $\frac{1}{\tau_{1}}=\frac{1}{\tau_{0}}+\frac{1}{5}$, the previous conditions is equivalent to $\frac{1}{\delta}+\frac{1}{\tau_{0}}\leq \frac{1-\alpha}{5}$, which is exactly the condition stated in the point \emph{2)}.
\item [(d)] The term that contains the pressure can be treated as follows: by the formula (\ref{EquationPression}) we have\\ $\displaystyle{P=\frac{1}{(-\Delta)}\sum_{i,j=1}^3\partial_i\partial_j(u_ib_j)}$, so we need to study the quantity
$$\|\phi \vn P\|_{\mathcal{M}_{t,x}^{\mathfrak{p}_0, \mathfrak{q}_0}}\leq\sum_{i,j=1}^3\left\|\phi\left (\frac{\vn}{(-\Delta)} \partial_i\partial_j(u_ib_j)\right)\right\|_{\mathcal{M}_{t,x}^{\mathfrak{p}_0, \mathfrak{q}_0}},$$
but since we assumed in \emph{4)} that $\mathds{1}_{Q_{R_1}}\frac{\vn}{(-\Delta)} \partial_i\partial_j(u_ib_j)\in\mathcal{M}_{t,x}^{\mathfrak{p}, \mathfrak{q}}$ with $\mathfrak{p}_0\leq \mathfrak{p}<+\infty$ and $\mathfrak{q}_0\leq \mathfrak{q}<+\infty$, then by Lemma \ref{Lemme_locindi} we obtain for all $1\leq i,j\leq 3$:
$$\left\|\phi\left (\frac{\vn\partial_i\partial_j}{(-\Delta)} (u_ib_j)\right)\right\|_{\mathcal{M}_{t,x}^{\mathfrak{p}_0, \mathfrak{q}_0}}\leq \left\|\mathds{1}_{Q_{R_1}}\frac{\vn\partial_i\partial_j}{(-\Delta)} (u_ib_j)\right\|_{\mathcal{M}_{t,x}^{\mathfrak{p}, \mathfrak{q}}}<+\infty.$$
\item [(e)] For the last term of (\ref{TermeGeneral_a_Estimer}), we need to study $\|\phi (\vf+\vg)\|_{\mathcal{M}_{t,x}^{\mathfrak{p}_0, \mathfrak{q}_{0}}}$, but since $1<\mathfrak{p}_0\leq \tfrac{10}{7}$ and since $\mathfrak{q}_0=\frac{5}{2-\alpha}<\tau_a, \tau_b$, then from the first point of Lemma \ref{Lemme_Product} and from Lemma \ref{Lemme_locindi} we have
$$\|\phi (\vf+\vg)\|_{\mathcal{M}_{t,x}^{\mathfrak{p}_0, \mathfrak{q}_{0}}}\leq C\|\mathds{1}_{Q_{R_1}} (\vf+\vg)\|_{\mathcal{M}_{t,x}^{\frac{10}{7}, \min\{\tau_{a}, \tau_{b}\}}}\leq C\left(\|\mathds{1}_{Q_{R_3}} \vf\|_{\mathcal{M}_{t,x}^{\frac{10}{7}, \tau_{a}}}+\|\mathds{1}_{Q_{R_3}} \vg\|_{\mathcal{M}_{t,x}^{\frac{10}{7}, \tau_{b}}}\right)<+\infty.$$
\end{itemize}
This completes the proof of Proposition \ref{Proposition_Principale1}.\hfill$\blacksquare$\\

As mentioned in the Remark \ref{Remark_Point2}, Theorem \ref{Teorem1} follows from the Proposition \ref{Proposition_Principale1}, so it suffices to show the five hypotheses stated in Proposition \ref{Proposition_Principale1},  \emph{i.e.} we will prove points \emph{1)-5)} from the general hypotheses of Theorem \ref{Teorem1}. To be more precise, the points \emph{1)} and \emph{2)} will be shown in Section \ref{Secc_InductiveArgument} by using the estimates given in Section \ref{Secc_AveragedBounds}, while the points \emph{3)} and \emph{4)} will be verified in Section \ref{Secc_LastHypotheses}. Note also that since $R_3<2R_0$ and $Q_{4R_0}\subset \Omega$ by (\ref{Defi_BouleParaboliquedeBase}), then hypothesis \emph{5)} of Proposition \ref{Proposition_Principale1} follows from the third hypothesis of Theorem \ref{Teorem1}.

\section{Local bounds}\label{Secc_AveragedBounds}
Remark that all the information assumed in the hypotheses of Proposition \ref{Proposition_Principale1} is presented in the framework of Morrey spaces, thus to carry on our study it will be useful to fix some averaged quantities: for a point $(t,x)\in Q_{R_0}(t_0, x_0)$ and for a general radius $0<r<R_3$, following the notation (\ref{Def_ParabolicBalldeBase}) we consider the parabolic ball
$$Q_r(t,x)=]t-r^2,t+r^2[\times B(x,r),$$
and when the context is clear, we will write $Q_r$ and $B_r$ instead of $Q_r(t,x)$ and $B(x,r)$.\\ 

We define now the following dimensionless quantities (in the sense that they are scale invariant):
\begin{align}
& \mathcal{A}_r(t,x)=\sup_{t-r^2 <s< t+r^2} \frac{1}{r} \int_{B(x,r)} |\vu(s,y)|^2dy,  & \alpha_r(t,x)&=\sup_{t-r^2 < s< t+r^2} \frac{1}{r} \int_{B(x,r)} |\vb(s,y)|^2dy, \nonumber\\
& \mathcal{B}_r(t,x)=\frac{1}{r}  \iint_{Q_r(t,x)} |\vn \otimes \vu(s,y)|^2dyds,  &\beta_r(t,x) &= \frac{1}{r}  \iint_{Q_r(t,x)} |\vn \otimes \vb(s,y)|^2dyds, \nonumber\\
& \mathcal{C}_r(t,x) = \frac{1}{r^2}  \iint_{Q_r(t,x)} |\vu(s,y)|^3dyds,  \quad &\gamma_r(t,x) &= \frac{1}{r^2}  \iint_{Q_r(t,x)} |\vb(s,y)|^3dyds,\label{Definition_Averaged_Quantities} \\
& \mathcal{D}_r(t,x) = \frac{1}{r^{\frac57}}  \iint_{Q_r(t,x)} |\vf(s,y)|^{\frac{10}{7}}dyds, \quad\; &\delta_r(t,x) &= \frac{1}{r^{\frac 57}} \iint_{Q_r(t,x)} |\vg(s,y)|^{\frac{10}{7}}dyds,\nonumber\\
& \mathcal{P}_{r}(t,x) = \frac{1}{r^{5-2q_0}}  \iint_{Q_r(t,x)} |P(s,y)|^{q_0}dyds\mbox{ with } \tfrac{10}{7}<q_0\leq \tfrac32.\nonumber
\end{align}
The aim of this section is to obtain two inequalities (given in Proposition \ref{Propo_FirstEstimate} and in Proposition \ref{Propo_SecondEstimate} below) that involves all the previous quantities. These inequalities are necessary to apply an inductive procedure that will lead us to some of the controls assumed in Proposition \ref{Proposition_Principale1}. This inductive argument will be displayed in the next section.\\

In the following lemma we exhibit a first relationship between some of the terms in (\ref{Definition_Averaged_Quantities}) that will be used in Proposition \ref{Propo_FirstEstimate}.
\begin{Lemme}\label{Lemme_L3norm} Under the general hypotheses of Theorem \ref{Teorem1}, for any $0<r<R_3$, there exists an absolutely constant $C$, which does not depend on $r$, such that we have
$$\mathcal{C}_r^{\frac 13} \leq  C  (\mathcal{A}_r + \mathcal{B}_r)^{\frac12}, \quad \mbox{and}\quad \gamma_r^{\frac13} \leq C(\alpha_r + \beta_r)^{\frac12}.$$
\end{Lemme}
{\bf Proof.} We only detail the proof of the first estimate as the second follows the same computations. Thus, by the definition of $\mathcal{C}_r$ given in (\ref{Definition_Averaged_Quantities}) and Hölder's inequality, we have 
$$\mathcal{C}_r^{\frac 13} = \frac{1}{r^{\frac 23}} \|\vu\|_{L_{t,x}^{3} (Q_r)} 
\leq C \frac{1}{r^{\frac 12}} \|\vu\|_{L_{t,x}^{\frac{10}{3}} (Q_r)}.$$
Now we remark that we have the interpolation inequality $\|\vu(t,\cdot)\|_{L^{\frac{10}{3}}(B_r)}\leq \|\vu(t,\cdot)\|_{L^{2}(B_r)}^{\frac{2}{5}}\|\vu(t,\cdot)\|_{L^{6}(B_r)}^{\frac{3}{5}}$ and applying the Hölder inequality with respect to the time variable, we obtain 
$$\|\vu\|_{L_{t,x}^{\frac{10}{3}} (Q_r)} \leq \|\vu\|_{L_t^{\infty}L_x^{2} (Q_r)}^{\frac25}\|\vu\|_{L_t^{2}L_x^{6} (Q_r)}^{\frac{3}{5}}.$$
For the $L_t^2L_x^6$ norm of $\vu$, we use the classical Gagliardo-Nirenberg inequality (see \cite{Brezis}) to obtain 
$$\|\vu\|_{L_t^{2}L_x^{6} (Q_r)} \leq C\left(\|\vn \otimes\vu\|_{L_t^{2}L_x^{2} (Q_r)} +\|\vu\|_{L_t^{\infty}L_x^{2} (Q_r)}\right),$$
and using Young's inequalities we have
\begin{eqnarray}
\|\vu\|_{L_{t,x}^{\frac{10}{3}} (Q_r)} &\leq& C \|\vu\|_{L_t^{\infty}L_x^{2} (Q_r)}^{\frac 25}\left(\|\vn \otimes \vu\|_{L_t^{2}L_x^{2} (Q_r)}^{\frac 35}+\|\vu\|_{L_t^{\infty}L_x^{2} (Q_r)}^{\frac35} \right)\notag\\
& \leq&C\left(\|\vu\|_{L_t^{\infty}L_x^{2} (Q_r)}+\|\vn\otimes \vu\|_{L_t^{2}L_x^{2} (Q_r)}\right).\label{EstimationUtileUL103}
\end{eqnarray}
Noting that $\|\vu\|_{L_t^{\infty}L_x^{2} (Q_r)}=r^{\frac12}\mathcal{A}_r^{\frac12}$ and $\|\vn\otimes \vu\|_{L^2_tL_x^{2} (Q_r)}=r^{\frac12}\mathcal{B}_r^{\frac12}$, we finally obtain $\mathcal{C}_r^{\frac 13} \leq  C  (\mathcal{A}_r + \mathcal{B}_r)^{\frac12}$ and Lemma \ref{Lemme_L3norm} is proven. \hfill$\blacksquare$\\

We give now the first general inequality that bounds all the term defined in formula (\ref{Definition_Averaged_Quantities}).
\begin{Proposition}[First Estimate]\label{Propo_FirstEstimate}
Under the hypotheses of Theorem \ref{Teorem1}, for $0<r<\frac{\rho}{2} \leq  \frac{R_3}{2}$, we have
\begin{equation*}
\begin{split}
\mathcal{A}_{r} +  \mathcal{B}_{r}+ \alpha_{r} + \beta_{r} & \leq  C \frac{r^2}{\rho^2}  (\mathcal{A}_\rho + \alpha_\rho)  +C \frac{\rho^2}{r^2} \left( (\mathcal{A}_\rho + \alpha_\rho + \beta_\rho)  \mathcal{B}_\rho^{\frac 12} + (\alpha_\rho + \mathcal{A}_\rho + \mathcal{B}_\rho)  \beta_\rho^{\frac12} \right) \\
& +  C \frac{\rho^2}{r^2} \mathcal{P}_{\rho}^{\frac{1}{q_0}} \left( (\mathcal{A}_\rho+ \mathcal{B}_\rho)^{\frac 12} + (\alpha_\rho + \beta_\rho)^{\frac12}\right) \\
& +  C\frac{\rho}{r} \left(\mathcal{D}_\rho^{\frac{7}{10}}  (\mathcal{A}_\rho + \mathcal{B}_\rho)^{\frac 12} + \delta_\rho^{\frac{7}{10}}  (\alpha_\rho + \beta_\rho)^{\frac12}  \right).
\end{split}
\end{equation*}
\end{Proposition}
{\bf Proof of Proposition \ref{Propo_FirstEstimate}.} To obtain this estimate we will use the local energy estimate satisfied by solutions of equation (\ref{EquationMHD}).  It is crucial to choose here a good test function and following Scheffer \cite{Scheffer} we will consider the non-negative function $\omega\in \mathcal{C}^\infty_0(\mathbb{R}\times\mathbb{R}^3)$ defined by the formula
\begin{equation}\label{Definition_FuncTestInegaliteEnergie}
\omega(s,y)= r^2\phi\left(\frac{s-t}{\rho^2}, \frac{y-x}{\rho}\right)\theta\left(\frac{s-t}{r^2}\right)\mathfrak{g}_{(4r^2+t-s)}(x-y), \qquad 0<r<\frac{\rho}{2}\leq \frac{R_3}{2}, 
\end{equation}
where $\phi\in  \mathcal{C}^\infty_0(\mathbb{R}\times\mathbb{R}^3)$ is a non-negative function supported on $]-1, 1[\times B(0, 1)$ and is equal to $1$ on $]-\tfrac{1}{4}, \tfrac{1}{4}[\times B(0, \tfrac12)$ and $\theta:\mathbb{R}\longrightarrow \mathbb{R}$ is a non-negative smooth function such that $\theta\equiv 1$ on $]-\infty, 1[$ and $\theta\equiv 0$ on $]2, +\infty[$ and $\mathfrak{g}_t(x)$ is the usual heat kernel.\\

We gather in the following lemma some properties of this test function:
\begin{Lemme}\label{Lemme_PropertiesTestFunctionEnergyInequality}
Recalling that $0<r<\frac{\rho}{2}$ (and thus $Q_r(t,x)\subset Q_{\rho}(t,x)$), we have
\begin{itemize}
\item[1)] the function $\omega$ is a bounded non-negative smooth function and its support is contained in the parabolic ball $Q_{\rho}(t,x)$ and for all $(s,y)\in Q_{r}(t,x)$ we have the lower bound $\omega(s,y)\geq \frac{C}{r}$,
\item[2)] for all $(s,y)\in Q_{\rho}(t,x)$ with $0<s< t + r^2$ we have $\omega(s,y)\leq \frac{C}{r}$,
\item[3)] for all $(s,y)\in Q_{\rho}(t,x)$ with $0<s< t + r^2$ we have $|\vn \omega(s,y)|\leq \frac{C}{r^2}$,
\item[4)] moreover, for all $(s,y)\in Q_{\rho}(t,x)$ with $0<s< t + r^2$ we have 
$|(\partial_s + \Delta) \omega(s,y)|\leq C \frac{r^2}{\rho^5}$.
\end{itemize}
\end{Lemme}
See the Appendix \ref{Secc_AppendixB} for a proof of this lemma. Now, with this particular test function $\omega$, we can construct the following local energy inequality
\begin{eqnarray}
\int_{\mathbb{R}^3} \left(|\vu(\tau,y)|^2 + |\vb(\tau,y)|^2 \right)\omega(\tau,y)dy+2 \int_{s<\tau}\int_{\mathbb{R}^3} \left(|\vn \otimes \vu(s,y)|^2 + |\vn \otimes \vb(s,y)|^2 \right) \omega(s,y)dyds\notag\\
\leq \int_{s<\tau}\int_{\mathbb{R}^3} \left(|\vu(s,y)|^2 + |\vb(s,y)|^2 \right) (\partial_t+ \Delta) \omega(s,y)dyds\qquad\qquad \notag\\
+\int_{s<\tau}\int_{\mathbb{R}^3} \left(|\vu(s,y)|^2 + 2P(s,y) \right) (\vb \cdot \vn) \omega(s,y)dyds\qquad\qquad\,\label{EstimationDEnergieLocale1}\\
+\int_{s<\tau}\int_{\mathbb{R}^3} \left(|\vb(s,y)|^2 + 2P(s,y) \right) (\vu \cdot \vn) \omega(s,y)dyds\qquad\qquad\notag\\
+2\int_{s<\tau}\int_{\mathbb{R}^3}  \left(\vf(s,y) \cdot \vu(s,y) + \vg(s,y) \cdot\vb(s,y)\right)\omega(s,y)dyds.\;\;\notag
\end{eqnarray}
Now, we define the quantities $\uuo_\rho$ and $\bbo_\rho$ as the following averages:
\begin{equation}\label{Def_Averages}
\uuo_\rho(t,x)=\frac{1}{|B(x,\rho)|}\int_{B(x,\rho)} |\vu(t,y)|^2 \, dy,\qquad  \bbo_\rho(t,x)=\frac{1}{|B(x,\rho)|}\int_{B(x,\rho)} |\vb(t,y)|^2 \, dy,
\end{equation}
and since $\vu$ and $\vb$ are divergence free, for any test function $\phi$ compactly supported within $B(x,\rho)$, we have 
$$\int_{B(x,\rho)} \uuo_\rho \,(\vb \cdot \vn) \phi(t,y) \,dy = 0 \quad \text{and} \quad\int_{B(x,\rho)}\bbo_\rho \, (\vu \cdot \vn) \phi(t,y) \,dy = 0,$$
these facts will allow us to introduce the averages $\uuo_\rho$ and $\bbo_\rho$ in inequality (\ref{EstimationDEnergieLocale1}) in order to use Poincar\'e's inequality. Indeed, we can rewrite the previous local energy inequality in the following manner
\begin{eqnarray*}
\int_{\mathbb{R}^3} \left(|\vu(\tau,y)|^2 + |\vb(\tau,y)|^2 \right)\omega(\tau,y)dy+2 \int_{s<\tau}\int_{\mathbb{R}^3} \left(|\vn \otimes \vu(s,y)|^2 + |\vn \otimes \vb(s,y)|^2 \right) \omega(s,y)dyds\notag\\\leq \int_{s<\tau}\int_{\mathbb{R}^3} \left(|\vu(s,y)|^2 + |\vb(s,y)|^2 \right) (\partial_t+ \Delta) \omega(s,y)dyds\qquad\qquad \notag\\+\int_{s<\tau}\int_{\mathbb{R}^3} \left(|\vu(s,y)|^2 -\uuo_\rho\right) (\vb \cdot \vn) \omega(s,y)dyds\qquad\qquad\\+\int_{s<\tau}\int_{\mathbb{R}^3} \left(|\vb(s,y)|^2-\bbo_\rho \right) (\vu \cdot \vn) \omega(s,y)dyds\qquad\qquad\notag\\+C\int_{s<\tau}\int_{\mathbb{R}^3} P(s,y)\left((\vb\cdot\vn)\omega(s,y)+(\vu\cdot\vn)\omega(s,y)\right)dyds\quad\notag\\+C\int_{s<\tau}\int_{\mathbb{R}^3}  \left(\vf(s,y) \cdot \vu(s,y) + \vg(s,y) \cdot\vb(s,y)\right)\omega(s,y)dyds.\notag
\end{eqnarray*}
Using the properties of the test function $\omega$ stated in Lemma \ref{Lemme_PropertiesTestFunctionEnergyInequality} we have: 
\begin{equation}\label{EstimationDEnergieLocale2}
\begin{split}
\frac{1}{r} 
& \int_{B_r} |\vu(\tau,y)|^2 
+ |\vb(\tau,y)|^2dy 
+\frac{1}{r} \iint_{Q_r} |\vn \otimes \vu(s,y)|^2+ |\vn \otimes \vb(s,y)|^2 dyds\\
& \leq C\underbrace{\frac{r^2}{\rho^5}\iint_{Q_\rho} |\vu(s,y)|^2 + |\vb(s,y)|^2dyds}_{(I)}\qquad\qquad\\
&+\underbrace{\frac{C}{r^2}\iint_{Q_\rho} \left||\vu(s,y)|^2 -\uuo_\rho \right| |\vb(s,y)|dyds}_{(II)}
+\underbrace{\frac{C}{r^2}\iint_{Q_\rho} \left||\vb(s,y)|^2 -\bbo_\rho \right| |\vu(s,y)|dyds}_{(III)}\\
& +\underbrace{\frac{C}{r^2}\iint_{Q_\rho} |P(s,y)|\left(|\vu(s,y)|+|\vb(s,y)|\right)dyds}_{(IV)}
+\underbrace{\frac{C}{r}\iint_{Q_\rho}  |\vf(s,y)||\vu(s,y)| + |\vg(s,y)||\vb(s,y)|dyds}_{(V)}.
\end{split}
\end{equation}
We will study each one of the previous terms separately. The first term on the right-hand side above is easy to bound: indeed, by definition of the quantities $\mathcal{A}_\rho$ and $\alpha_\rho$ given in (\ref{Definition_Averaged_Quantities}), we get directly
\begin{equation}
(I) \leq  C \frac{r^2}{\rho^2}  \left(\sup_{t-\rho^2<s< t+\rho^2} \frac{1}{\rho} \int_{B_\rho} |\vu(s,y)|^2dy + \sup_{t-\rho^2 <s< t+\rho^2} \frac{1}{\rho} \int_{B_\rho} |\vb(s,y)|^2dy \right) \leq C \frac{r^2}{\rho^2}  (\mathcal{A}_\rho + \alpha_\rho). 
\label{I1}
\end{equation}
The terms $(II)$ and $(III)$ can be treated in the same fashion since we have symmetric information on the functions $\vu$ and $\vb$, so we only study one of them: indeed, for $(II)$ we have
$$\frac{1}{r^2}\iint_{Q_\rho} \left||\vu(s,y)|^2 -\uuo_\rho \right| |\vb(s,y)|dyds
\leq \frac{1}{r^2}  \int_{t-\rho^2}^{t+\rho^2} \||\vu(s,\cdot)|^2 - \uuo_\rho\|_{L^{\frac{3}{2}} (B_\rho)}\|\vb(s, \cdot)\|_{L^{3} (B_\rho)}ds,$$
thus, by the Poincaré inequality we obtain 
\begin{eqnarray*}
(II)&\leq&\frac{C}{r^2}  \int_{t-\rho^2}^{t+\rho^2} \|\vn(|\vu(s,\cdot)|^2 )\|_{L^{1} (B_\rho)} \|\vb(s,\cdot)\|_{L^{3} (B_\rho)}ds \\
&\leq&\frac{C}{r^2}  \int_{t-\rho^2}^{t+\rho^2} \|\vu(s,\cdot)\|_{L^{2} (B_\rho)} \|\vn \otimes\vu\|_{L^{2} (B_\rho)} \|\vb(s, \cdot)\|_{L^{3} (B_\rho)} ds\\
&\leq&\frac{C}{r^2}  \|\vu\|_{L_t^{6}L_x^{2} (Q_\rho)} \|\vn \otimes \vu\|_{L_{t,x}^{2} (Q_\rho)} \|\vb\|_{L_{t,x}^{3} (Q_\rho)},
\end{eqnarray*}
where we used the H\"older inequality in the time variable in the last estimate. Now we remark that we have the following bounds for $ \|\vu\|_{L_t^{6}L_x^{2} (Q_\rho)}$, $\|\vn \otimes \vu\|_{L^{2}_{t,x} (Q_\rho)}$ and $\|\vb\|_{L^{3}_{t,x} (Q_\rho)}$ (recall the expressions given in (\ref{Definition_Averaged_Quantities})):
\begin{eqnarray*}
\|\vu\|_{L_t^{6}L_x^{2} (Q_\rho)}\leq C \rho^{\frac13} \|\vu\|_{L_t^{\infty}L_x^{2} (Q_\rho)}\leq C \rho^{\frac56}  \left(\sup_{t-\rho^2 <s<t+\rho^2} \frac{1}{\rho} \int_{B_\rho} |\vu(s,y)|^2dy \right)^{\frac12} = C \rho^{\frac56} \mathcal{A}_\rho^{\frac12},\\
\|\vn \otimes \vu\|_{L_{t,x}^{2} (Q_\rho)}=\rho^{\frac12}\mathcal{B}_\rho^{\frac12}\quad \mbox{and}\quad \|\vb\|_{L_{t,x}^{3}(Q_\rho)}=\rho^{\frac23}  \gamma_\rho^{\frac13},\hspace{3cm}
\end{eqnarray*}
we obtain then 
$$(II)\leq C\frac{\rho^2}{r^2}\mathcal{A}_\rho^{\frac12}\mathcal{B}_\rho^{\frac12} \gamma_\rho^{\frac13} \leq C \frac{\rho^2}{r^2}\mathcal{A}_\rho^{\frac12}\mathcal{B}_\rho^{\frac12}(\alpha_\rho + \beta_\rho)^{\frac12}$$
where we used Lemma \ref{Lemme_L3norm} to estimate the term $\gamma_\rho^{\frac13}$. Since the same computations can be performed for $(III)$ we have
\begin{eqnarray}
(II)+(III)&\leq &C\frac{\rho^2}{r^2} \left(\mathcal{A}_\rho^{\frac12}\mathcal{B}_\rho^{\frac12} (\alpha_\rho + \beta_\rho)^{\frac12} +   \alpha_\rho^{\frac12} \beta_\rho^{\frac12} (\mathcal{A}_\rho + \mathcal{B}_\rho)^{\frac12} \right)\notag\\
& \leq  &C \frac{\rho^2}{r^2} \left( (\mathcal{A}_\rho + \alpha_\rho + \beta_\rho)  \mathcal{B}_\rho^{\frac12} + (\alpha_\rho + \mathcal{A}_\rho + \mathcal{B}_\rho)  \beta_\rho^{\frac12} \right).\label{I2}
\end{eqnarray}
We study now the term $(IV)$ of (\ref{EstimationDEnergieLocale2}). Using Hölder's inequality, we have with $\frac{1}{q_0}+\frac{1}{q_0'}=1$:
$$(IV)=\frac{C}{r^2}\iint_{Q_\rho} |P(s,y)|\left(|\vu(s,y)|+ |\vb(s,y)|\right)dyds \leq\frac{C}{r^2} \|P\|_{L_{t,x}^{q_0}(Q_\rho)}\left(\|\vu\|_{L_{t,x}^{q_0'}  (Q_\rho)} +  \|\vb\|_{L_{t,x}^{q_0'} (Q_\rho)}\right).$$
Since we have $\tfrac{10}{7}<q_0\leq \tfrac{3}{2}$ and $3\leq q_0'<\tfrac{10}{3}$ we can write\\
$$(IV)\leq \frac{C}{r^2} \|P\|_{L_{t,x}^{q_0}(Q_\rho)}\rho^{5(\frac{7}{10}-\frac{1}{q_0})}\left(\|\vu\|_{L_{t,x}^{\frac{10}{3}}  (Q_\rho)} +  \|\vb\|_{L_{t,x}^{\frac{10}{3}} (Q_\rho)}\right)$$
Since by definition (see expression (\ref{Definition_Averaged_Quantities})) we have 
$\rho^{(\frac{5}{q_0}-2)}\mathcal{P}_\rho^{\frac{1}{q_0}}=\|P\|_{L_{t,x}^{q_0}(Q_\rho)}$ and since by (\ref{EstimationUtileUL103}) we have the estimates $\|\vu\|_{L_{t,x}^{\frac{10}{3}} (Q_r)}\leq C\rho^{\frac12}(\mathcal{A}_\rho + \mathcal{B}_\rho)^{\frac12}$ and $\|\vb\|_{L_{t,x}^{\frac{10}{3}} (Q_r)}\leq C\rho^{\frac12}(\alpha_\rho + \beta_\rho)^{\frac12}$, then we obtain
\begin{equation} \label{I3}
(IV)\leq C\frac{\rho^2}{r^2}\mathcal{P}_\rho^{\frac{1}{q_0}}\left( (\mathcal{A}_\rho + \mathcal{B}_\rho)^{\frac12} + (\alpha_\rho + \beta_\rho)^{\frac12}\right). 
\end{equation}
Finally for the last term $(V)$ of (\ref{EstimationDEnergieLocale2}) we have by the H\"older inequality
$$(V)=\frac{C}{r}\iint_{Q_\rho}  |\vf||\vu| + |\vg||\vb|dyds\leq C \frac{1}{r}  \left( \|\vf\|_{L_{t,x}^{\frac{10}{7}}(Q_\rho)}\|\vu\|_{L_{t,x}^{\frac{10}{3}}(Q_\rho)} +  \|\vg\|_{L_{t,x}^{\frac{10}{7}}(Q_\rho)}\|\vb\|_{L_{t,x}^{\frac{10}{3}}(Q_\rho)} \right).$$
Recalling the control $\|\vu\|_{L_{t,x}^{\frac{10}{3}} (Q_\rho)}\leq C\left(\|\vu\|_{L_t^{\infty}L_x^{2} (Q_\rho)}+\|\vn\otimes \vu\|_{L_{t,x}^{2} (Q_\rho)}\right)$ (see inequality (\ref{EstimationUtileUL103})) and since by (\ref{Definition_Averaged_Quantities}) we have the identities $\|\vu\|_{L_t^{\infty}L_x^{2} (Q_\rho)}=\rho^{\frac12}\mathcal{A}_\rho^{\frac12}$, $\|\vn\otimes \vu\|_{L_{t,x}^{2} (Q_\rho)}=\rho^{\frac12}\mathcal{B}_\rho^{\frac12}$, $\rho^{\frac{1}{2}}\mathcal{D}_\rho^{\frac{7}{10}}=\|\vf\|_{L_{t,x}^{\frac{10}{7}}(Q_\rho)}$ and $\rho^{\frac{1}{2}}\delta_\rho^{\frac{7}{10}}=\|\vg\|_{L_{t,x}^{\frac{10}{7}}(Q_\rho)}$, we obtain:
\begin{equation}
(V)\leq C\frac{\rho}{r} \left(\mathcal{D}_\rho^{\frac{7}{10}} (\mathcal{A}_\rho + \mathcal{B}_\rho)^{\frac12} + \delta_\rho^{\frac{7}{10}}(\alpha_\rho + \beta_\rho)^{\frac12}  \right) \label{I4}.
\end{equation} 
Gathering the estimates \eqref{I1},\eqref{I2},\eqref{I3} and \eqref{I4}, we have
\begin{multline*}
\frac{1}{r} \int_{B_r} |\vu(\tau,y)|^2 + |\vb(\tau,y)|^2dy 
+\frac{1}{r} \iint_{Q_r} |\vn \otimes \vu(s,y)|^2+ |\vn \otimes \vb(s,y)|^2 dyds\\
\leq  C\frac{r^2}{\rho^2}(\mathcal{A}_\rho + \alpha_\rho)+C \frac{\rho^2}{r^2} \left( (\mathcal{A}_\rho + \alpha_\rho + \beta_\rho)  \mathcal{B}_\rho^{\frac12} + (\alpha_\rho + \mathcal{A}_\rho + \mathcal{B}_\rho)  \beta_\rho^{\frac12}\right)\\
+C \frac{\rho^2}{r^2} \mathcal{P}_\rho^{\frac{1}{q_0}} \left( (\mathcal{A}_\rho + \mathcal{B}_\rho)^{\frac12} + (\alpha_\rho + \beta_\rho)^{\frac12}\right)\hspace{3cm}\\
+C\frac{\rho}{r} \left(\mathcal{D}_\rho^{\frac{7}{10}} (\mathcal{A}_\rho + \mathcal{B}_\rho)^{\frac12} + \delta_\rho^{\frac{7}{10}}(\alpha_\rho + \beta_\rho)^{\frac12}\right).\hspace{6.1cm}
\end{multline*}
Since this estimate is uniform with respect of the time variable of the left-hand side, we finally can write:
\begin{equation*}
\begin{split}
\mathcal{A}_{r} +  \mathcal{B}_{r}+ \alpha_{r} + \beta_{r} & \leq  C \frac{r^2}{\rho^2}  (\mathcal{A}_\rho + \alpha_\rho)  +C \frac{\rho^2}{r^2} \left( (\mathcal{A}_\rho + \alpha_\rho + \beta_\rho)  \mathcal{B}_\rho^{\frac 12} + (\alpha_\rho + \mathcal{A}_\rho + \mathcal{B}_\rho)  \beta_\rho^{\frac12} \right) \\
& +  C \frac{\rho^2}{r^2} \mathcal{P}_\rho^{\frac{1}{q_0}} \left( (\mathcal{A}_\rho+ \mathcal{B}_\rho)^{\frac 12} + (\alpha_\rho + \beta_\rho)^{\frac12}\right) \\
& +  C\frac{\rho}{r} \left(\mathcal{D}_\rho^{\frac{7}{10}}  (\mathcal{A}_\rho + \mathcal{B}_\rho)^{\frac 12} + \delta_\rho^{\frac{7}{10}}  (\alpha_\rho + \beta_\rho)^{\frac12}  \right),
\end{split}
\end{equation*}
and Proposition \ref{Propo_FirstEstimate} is proven. \hfill$\blacksquare$\\

The second estimate that we need relies on a detailed study of the properties of the pressure and following Kukavica \cite{Kukavica} we have:
\begin{Proposition}[Second Estimate]\label{Propo_SecondEstimate}
With the quantities defined in (\ref{Definition_Averaged_Quantities}), under the hypotheses of Theorem \ref{Teorem1} and for $0<r<\frac{\rho}{2} \leq  \frac{R_3}{2}$ we have the estimate:
\begin{equation}\label{pestimscalq0}
\mathcal{P}_{r} \leq C \left( \left(\frac{\rho}{r} \right)^{3-\q} (\mathcal{A}_\rho \beta_\rho)^{\frac{q_0}{2}} + \left(\frac{r}{\rho} \right)^{2\q-2} \mathcal{P}_{\rho} \right).
\end{equation}
\end{Proposition}
In order to obtain the previous inequality we will first study a general estimate stated in the lemma below and then (\ref{pestimscalq0}) will follow by a scaling argument. 
\begin{Lemme}\label{Lemma_EstimationPression1}
For $0 < \sigma \leq \tfrac12$ and for a parabolic ball $Q_\sigma$, there is a constant $C$ such that whenever $P \in L_{t,x}^{q_0} (Q_1)$ for $1<q_0 \leq \frac{3}{2}$, $\Delta P= - \displaystyle{\sum^3_{i,j= 1} \partial_i \partial_j (u_i b_j)}$ in $Q_\sigma$, $\vu\in L_t^{\infty}L_x^{2} (Q_1)$ and $\vn \otimes \vb\in L_{t,x}^{2} (Q_1)$, then we have the following control
\begin{equation}\label{pressureesti32}
\norm{P}_{L_{t,x}^{q_0}  (Q_\sigma)} \leq C \left( \sigma^{\frac{2}{q_0}-1}\|\vu\|_{L_t^{\infty}L_x^{2} (Q_1)} \|\vn \otimes \vb\|_{L_{t,x}^{2} (Q_1)} +  
\sigma^{\frac{3}{q_0}}\|P\|_{L_{t,x}^{q_0} (Q_1)}\right),
\end{equation}
where $Q_\sigma$ and $Q_1$ are parabolic balls of radius $\sigma$ and $1$ respectively.
\end{Lemme}
{\bf Proof.} To obtain this inequality we introduce $\eta : \R \longrightarrow [0, 1]$ a smooth function supported in the ball $B_1$ such that $\eta \equiv 1$ on the ball $B_{\frac35}$ and  $\eta \equiv 0$ outside the ball $B_{\frac45}$. Note in particular that on $Q_\sigma$ we have the identity $P=\eta P$. Now a straightforward calculation shows that we have the identity
$$ - \Delta (\eta P) = -\eta \Delta P + (\Delta \eta)P - 2  \sum^3_{i= 1}\partial_i ((\partial_i \eta) P),$$
from which we deduce the inequality
\begin{equation}\label{FormuleEtaPression}
\|P\|_{L_{t,x}^{q_0}(Q_\sigma)} =\|\eta P\|_{L_{t,x}^{q_0}(Q_\sigma)} \leq \underbrace{\left\|\frac{\big(-\eta \Delta P\big)}{(-\Delta)}\right\|_{L_{t,x}^{q_0}(Q_\sigma)}}_{(I)} + \underbrace{\left\|\frac{(\Delta \eta) P}{(-\Delta)}\right\|_{L_{t,x}^{q_0}(Q_\sigma)} }_{(II)} +2\sum^3_{i= 1}\underbrace{\left\|\frac{\partial_i ( (\partial_i \eta) P)}{(-\Delta)}\right\|_{L_{t,x}^{q_0}(Q_\sigma)}}_{(III)}.
\end{equation}
For the first term of (\ref{FormuleEtaPression}), since $\Delta P=-\displaystyle{\sum^3_{i,j= 1}}\partial_i \partial_j (u_i b_j)$ on $Q_\sigma$, if  we denote by $N_{i,j} = u_i (b_j - (b_j)_1)$ where $ (b_j)_1$ is the average of $b_j$ over the ball of radius $1$ (recall the definition (\ref{Def_Averages}))  since $\vu$ is divergence free we have $\displaystyle{\sum^3_{i,j= 1}}\partial_i \partial_j (u_i b_j)=\displaystyle{\sum^3_{i,j= 1}}\partial_i \partial_j N_{i,j}$ and thus we can write
\begin{eqnarray}
(I)&=&\left\|\frac{\big(-\eta \Delta P\big)}{(-\Delta)}\right\|_{L_{t,x}^{q_0}(Q_\sigma)}\leq C\sigma^{5(\frac{1}{q_0}-\frac{2}{3})}\left\|\frac{1}{(-\Delta)}\Big(\eta \,  \sum^3_{i,j= 1}\partial_i \partial_j  N_{i,j}  \Big)\right\|_{L_{t,x}^{\frac 32}(Q_\sigma)}\notag\\
&\leq &C\sigma^{5(\frac{1}{q_0}-\frac{2}{3})}\sum^3_{i,j= 1}\left\|\frac{1}{(-\Delta)} \left(\partial_i \partial_j(\eta N_{i,j} ) - \partial_i \big((\partial_j \eta) N_{i,j} \big) - \partial_j \big((\partial_i \eta) N_{i,j} \big) + 2(\partial_i \partial_j \eta) N_{i,j}\right)\right\|_{L_{t,x}^{\frac 32}(Q_\sigma)}\label{FormulePression_I}
\end{eqnarray}
Denoting by $\mathcal{R}_i=\frac{\partial_i}{\sqrt{-\Delta}}$ the usual Riesz transforms on $\mathbb{R}^3$, by the boundedness of these operators in Lebesgue spaces and using the support properties of the auxiliary function $\eta$, we have for the first term above:
\begin{eqnarray*}
	\left\|\frac{\partial_i \partial_j}{(-\Delta)} \eta N_{i,j}(t,\cdot) \right\|_{L^{\frac32} (B_\sigma)} &\leq &\|\mathcal{R}_i \mathcal{R}_j (\eta N_{i,j} )(t,\cdot) \|_{L^{\frac32} (\R)}  \leq  C\|\eta N_{i,j}(t,\cdot)\|_{L^{\frac32} (B_1)}\\ 
	&\leq &C\|u_i(t,\cdot) \|_{L^{2}(B_1)} \|b_j(t,\cdot) - (b_j)_1\|_{L^{6} (B_1)}\\ &\leq &C\|\vu(t,\cdot)\|_{L^{2} (B_1)} \|\vn \otimes \vb(t,\cdot)\|_{L^{2} (B_1)},
\end{eqnarray*}
where we used H\"older and Poincaré inequalities in the last line. Now taking the $L^{\frac32}$-norm in the time variable of the previous inequality we obtain\\
\begin{equation}\label{p11}
\left\|\frac{\partial_i \partial_j}{(-\Delta)} \eta N_{i,j} \right\|_{L_{t,x}^{\frac32} (Q_\sigma)}\leq C\sigma^{\frac{1}{3}}\|\vu\|_{L^{\infty}_tL^{2}_x (Q_1)} \|\vn \otimes \vb\|_{L_{t,x}^{2} (Q_1)}. 
\end{equation}
The remaining terms of (\ref{FormulePression_I}) can all be studied in a similar manner. Indeed, noting that $\partial_i \eta$ vanishes on $B_{\frac35} \cup B^c_{\frac45}$ and since $B_\sigma \subset B_{\frac12}\subset B_{\frac{3}{5}}$, using the integral representation for the operator $\frac{\partial_i}{(- \Delta )}$ we have for the second term of (\ref{FormulePression_I}) the estimate
\begin{eqnarray}
\left\|\frac{\partial_i}{(- \Delta )}\big((\partial_j\eta)N_{i,j}\big)(t,\cdot)\right\|_{L^{\frac32}(B_\sigma)} &\leq &C\sigma^2\left\|\frac{\partial_i}{(- \Delta )}\big((\partial_j\eta)N_{i,j}\big)(t,\cdot)\right\|_{L^{\infty}(B_\sigma)}\notag\\
&\leq &C \, \sigma^2 \left\|\int_{\{\frac35<|y|<\frac45\}} \frac{x_i - y_i }{|x-y|^3} \big((\pj \eta) N_{i,j} \big)(t,y)  \, dy\right\|_{L^{\infty}(B_\sigma)}\notag\\
& \leq &C \, \sigma^2 \|N_{i,j}(t,\cdot)\|_{L^{1} (B_1)}\label{KernelEstimate1}\\
& \leq &C \, \sigma^2 \|u_i(t,\cdot) \|_{L^{2} (B_1)} \|b_j(t,\cdot) - (b_j)_1\|_{L^{2}(B_1)} \notag\\
&\leq& C \,  \|\vu(t,\cdot)\|_{L^{2} (B_1)}\|\vn \otimes \vb(t,\cdot)\|_{L^{2} (B_1)},\notag
\end{eqnarray}
where we used the same ideas as previously and the fact that $0<\sigma<1$, and with the same arguments as in (\ref{p11}) before, taking the $L^{\frac32}$-norm in the time variable, we obtain
\begin{equation}\label{p12}
\left\|\frac{\partial_i}{(- \Delta )}\big((\partial_j\eta)N_{i,j}\big)\right\|_{L^{\frac32}_{t,x}(Q_\sigma)} \leq C\sigma^{\frac{1}{3}}\|\vu\|_{L^{\infty}_tL^{2}_x (Q_1)} \|\vn \otimes \vb\|_{L_{t,x}^{2} (Q_1)}. 
\end{equation}
A symmetric argument gives 
\begin{equation}\label{p13}
\left\|\frac{\partial_{j}}{(- \Delta )}\big((\partial_i\eta)N_{i,j}\big)\right\|_{L^{\frac32}_{t,x}(Q_\sigma)} \leq C\sigma^{\frac{1}{3}}\|\vu\|_{L^{\infty}_tL^{2}_x (Q_1)} \|\vn \otimes \vb\|_{L_{t,x}^{2} (Q_1)}, \end{equation}
and observing that the convolution kernel associated to the operator $\frac{1}{(-\Delta)}$ is $\frac{C}{|x|}$, following the same ideas we have for the last term of (\ref{FormulePression_I}) the inequality
\begin{equation}\label{p14}
\left\|\frac{(\partial_i \partial_j \eta) N_{i,j}}{(-\Delta)}\right\|_{L_{t,x}^{\frac 32}(Q_\sigma)}\leq C\sigma^{\frac{1}{3}}\|\vu\|_{L^{\infty}_tL^{2}_x (Q_1)} \|\vn \otimes \vb\|_{L_{t,x}^{2} (Q_1)}.
\end{equation}
Therefore, combining the estimates \eqref{p11}, \eqref{p12},  \eqref{p13} and  \eqref{p14} and getting back to (\ref{FormulePression_I}) we finally have:
\begin{eqnarray}
(I)&=&\left\|\frac{\big(-\eta \Delta P\big)}{(-\Delta)}\right\|_{L_{t,x}^{\frac 32}(Q_\sigma)}\leq C\sigma^{5(\frac{1}{q_0}-\frac{2}{3})}\left(\sigma^{\frac{1}{3}}\|\vu\|_{L^{\infty}_tL^{2}_x (Q_1)} \|\vn \otimes \vb\|_{L_{t,x}^{2} (Q_1)}\right)\notag\\
&\leq &\sigma^{\frac{5}{q_0}-3}\|\vu\|_{L^{\infty}_tL^{2}_x (Q_1)} \|\vn \otimes \vb\|_{L_{t,x}^{2} (Q_1)}\label{FormulePression_I1}
\end{eqnarray}
We continue our study of expression (\ref{FormuleEtaPression}) and for the  term $(II)$ we first treat the space variable. Recalling the support properties of the auxiliary function $\eta$ and properties of the convolution kernel associated to the operator $\frac{1}{(-\Delta)}$, we can write as before (see (\ref{KernelEstimate1})):
$$\left\|\frac{(\Delta \eta) P(t,\cdot)}{(-\Delta)}\right\|_{L^{q_0}(B_\sigma)}\leq C\sigma^{\frac{3}{q_0}} \| P(t,\cdot)\|_{L^1(B_1)}\leq C\sigma^{\frac{3}{q_0}} \| P(t,\cdot)\|_{L^{q_0}(B_1)},$$
and thus, taking the $L^{q_0}$-norm in the time variable we obtain:
\begin{equation}\label{FormulePression_II}
(II)=\left\|\frac{(\Delta \eta) P}{(-\Delta)}\right\|_{L^{q_0}_{t,x}(Q_\sigma)}\leq C\sigma^{\frac{3}{q_0}} \|P\|_{L^{q_0}_{t,x}(Q_1)}.
\end{equation}
For the last term of expression (\ref{FormuleEtaPression}), following the same ideas developed in (\ref{KernelEstimate1}) we can write
$$\left\|\frac{\partial_i}{(-\Delta)}  (\partial_i \eta) P(t,\cdot) \right\|_{L^{q_0}(B_\sigma)}\leq C\sigma^{\frac{3}{q_0}} \| P(t,\cdot)\|_{L^1(B_1)}\leq C\sigma^{\frac{3}{q_0}} \| P(t,\cdot)\|_{L^{q_0}(B_1)},$$
and we obtain 
\begin{equation}\label{FormulePression_III}
(III)=\left\|\frac{\partial_i ( (\partial_i \eta) P )}{(-\Delta)}\right\|_{L_{t,x}^{q_0}(Q_\sigma)}\leq C\sigma^{\frac{3}{q_0}} \|P\|_{L^{q_0}_{t,x}(Q_1)}.
\end{equation}
Now, gathering the estimates (\ref{FormulePression_I1}), (\ref{FormulePression_II}) and (\ref{FormulePression_III}) we obtain the inequality 
$$\norm{P}_{L_{t,x}^{q_0}  (Q_\sigma)} \leq C \left( \sigma^{\frac{5}{q_0}-3}\|\vu\|_{L_t^{\infty}L_x^{2} (Q_1)} \|\vn \otimes \vb\|_{L_{t,x}^{2} (Q_1)} +  
\sigma^{\frac{3}{q_0}}\|P\|_{L_{t,x}^{q_0} (Q_1)}\right),$$
recalling at this point that since $1<q_0<\tfrac{3}{2}$, we have $\frac{2}{q_0}-1< \frac{5}{q_0}-3$ and since $0<\sigma\leq \tfrac{1}{2}$ we have $\sigma^{\frac{5}{q_0}-3}\leq \sigma^{\frac{2}{q_0}-1}$ and we finally obtain the estimate 
$$\norm{P}_{L_{t,x}^{q_0}  (Q_\sigma)} \leq C \left( \sigma^{\frac{2}{q_0}-1}\|\vu\|_{L_t^{\infty}L_x^{2} (Q_1)} \|\vn \otimes \vb\|_{L_{t,x}^{2} (Q_1)} +  
\sigma^{\frac{3}{q_0}}\|P\|_{L_{t,x}^{q_0} (Q_1)}\right),$$ 
and the proof of Lemma \ref{Lemma_EstimationPression1} is finished.\hfill $\blacksquare$\\

\noindent{\bf Proof of Proposition \ref{Propo_SecondEstimate}}. Once we have established the estimate (\ref{pressureesti32}) it is quite simple to deduce inequality \eqref{pestimscalq0}. Indeed, if we fix $\sigma = \frac{r}{\rho} \leq \frac12$ and if we introduce the functions $P_\rho(t,x)=P(\rho^2t, \rho x)$, $\vu_\rho(t,x)=\vu(\rho^2t, \rho x)$ and $\vb_\rho(t,x)=\vb(\rho^2t, \rho x)$ then from (\ref{pressureesti32}) we have
$$\|P_\rho\|_{L_{t,x}^{q_0}  (Q_{\frac{r}{\rho}})} \leq C \left( \left(\frac{r}{\rho}\right)^{\frac{2}{q_0}-1}\|\vu_\rho\|_{L_t^{\infty}L_x^{2} (Q_1)} \|\vn \otimes \vb_\rho\|_{L_{t,x}^{2} (Q_1)} +\left(\frac{r}{\rho}\right)^{\frac{3}{q_0}}\|P_\rho\|_{L_{t,x}^{q_0} (Q_1)}\right),$$
and by a convenient change of variable we obtain
$$\|P\|_{L_{t,x}^{q_0}  (Q_{r})}\rho^{-\frac{5}{q_0}} \leq C \left( \left(\frac{r}{\rho}\right)^{\frac{2}{q_0}-1}\rho^{-\frac32}\|\vu\|_{L_t^{\infty}L_x^{2} (Q_\rho)} \rho^{-\frac32}\|\vn \otimes \vb\|_{L_{t,x}^{2} (Q_\rho)} +\left(\frac{r}{\rho}\right)^{\frac{3}{q_0}}\rho^{-\frac{5}{q_0}}\|P\|_{L_{t,x}^{q_0} (Q_\rho)}\right).$$
Now, recalling that by (\ref{Definition_Averaged_Quantities}) we have the identities 
$$r^{\frac{5}{q_0}-2}\mathcal{P}_r^{\frac{1}{q_0}}=\|P\|_{L_{t,x}^{q_0}  (Q_r)}, \quad \rho^{\frac{1}{2}}\mathcal{A}_\rho^{\frac{1}{2}}=\|\vu\|_{L_t^{\infty}L_x^{2} (Q_\rho)} \quad \mbox{and}\quad \rho^{\frac{1}{2}}\beta_\rho^{\frac{1}{2}}=\|\vn \otimes \vb\|_{L_{t,x}^{2} (Q_\rho)},$$
we obtain
$$\mathcal{P}_r^{\frac{1}{q_0}}\leq C\left(\left(\frac{\rho}{r}\right)^{\frac{3}{q_0}-1}(\mathcal{A}_\rho\beta_\rho)^{\frac{1}{2}}+\left(\frac{r}{\rho}\right)^{2-\frac{2}{q_0}}\mathcal{P}_\rho^{\frac{1}{q_0}}\right),$$
and we finish the proof of Proposition \ref{Propo_SecondEstimate} by taking all this inequality to the $q_0$-power.\hfill $\blacksquare$
\section{Inductive argument}\label{Secc_InductiveArgument}
In Section \ref{Secc_AveragedBounds}, we have proven the following relationships between the averaged quantities defined in the expression (\ref{Definition_Averaged_Quantities}):
\begin{equation}\label{EstimationAveragedQuantities1}
\begin{split}
\mathcal{A}_{r} +  \mathcal{B}_{r}+ \alpha_{r} + \beta_{r} & \leq  C \frac{r^2}{\rho^2}  (\mathcal{A}_\rho + \alpha_\rho)  +C \frac{\rho^2}{r^2} \left( (\mathcal{A}_\rho + \alpha_\rho + \beta_\rho)  \mathcal{B}_\rho^{\frac 12} + (\alpha_\rho + \mathcal{A}_\rho + \mathcal{B}_\rho)  \beta_\rho^{\frac12} \right) \\
& +  C \frac{\rho^2}{r^2} \mathcal{P}_{\rho}^{\frac{1}{q_0}} \left( (\mathcal{A}_\rho+ \mathcal{B}_\rho)^{\frac 12} + (\alpha_\rho + \beta_\rho)^{\frac12}\right) \\
& +  C\frac{\rho}{r} \left(\mathcal{D}_\rho^{\frac{7}{10}}  (\mathcal{A}_\rho + \mathcal{B}_\rho)^{\frac 12} + \delta_\rho^{\frac{7}{10}}  (\alpha_\rho + \beta_\rho)^{\frac12}  \right)\\[4mm]
\mbox{and}\\
\mathcal{P}_{r} &\leq C \left( \left(\frac{\rho}{r} \right)^{3-\q} (\mathcal{A}_\rho \beta_\rho)^{\frac{q_0}{2}} + \left(\frac{r}{\rho} \right)^{2\q-2} \mathcal{P}_{\rho} \right).
\end{split}
\end{equation}
In this section we will see how to use  these relationships to obtain some of the local Morrey information assumed in Proposition \ref{Proposition_Principale1}. Indeed, we have:
\begin{Proposition}\label{Propo_indu_argu}
Let $(\vu, P, \vb)$ be a suitable solution of MHD equations (\ref{EquationMHD}) over $\Omega$ in the sense of Definition \ref{Def_SuitableSolutions}. Recall that in the framework of the general assumptions of Theorem \ref{Teorem1}, we have the following local information on the pressure $P \in L^{q_{0}}_{t,x}(\Omega) \mbox{ with } 1<q_0\leq  \tfrac{3}{2}$ and on the external forces $\vf$ and $\vg$: $\mathds{1}_{\Omega}\vf \in \mathcal{M}_{t,x}^{\frac{10}{7}, \tau_{a}}$ and $ \mathds{1}_{\Omega}\vg \in \mathcal{M}_{t,x}^{\frac{10}{7}, \tau_{b}}$ for some $\tau_c=\min\{\tau_{a}, \tau_b\}>\frac{5}{2-\alpha}>\frac{5}{3}$ with $0<\alpha<\frac13$.\\
	
Define now a real parameter $\tau_{0}$ such that $\frac{5}{1-\alpha}<\tau_{0}< 5\q$ and $2-\frac{5}{\tau_c}+\frac{5}{\tau_{0}}>0$. There exists a positive constant $\epsilon^{*}$ which depends only on $\tau_{a}, \tau_b$ and $\tau_{0}$ such that, if $\left(t_{0}, x_{0}\right) \in \Omega$ and
\begin{equation}\label{smallgrad}
\limsup _{r \rightarrow 0} \frac{1}{r} \iint_{\left]t_{0}-r^{2}, t_{0}+r^{2}\right[ \times B(x_0, r)}|\vn \otimes \vu \,(s, y)|^{2} + |\vn \otimes \vb \,(s, y)|^{2} dyds<\epsilon^{*}, 
\end{equation}
then there exists a (parabolic) neighborhood $Q_{R_3}$ of $(t_{0}, x_{0})$ with $0<R_0<R_1<R_2 <R_3< 2R_0$ such that 
$$\mathds{1}_{Q_{R_3}} \vec{u} \in \mathcal{M}_{t,x}^{3, \tau_{0}}, \quad \mathds{1}_{Q_{R_3}} \vb \in \mathcal{M}_{t,x}^{3, \tau_{0}}\quad \text{and} \quad \mathds{1}_{Q_{R_3}} P \in \mathcal{M}_{t,x}^{\q, \frac{\tau_{0}}{2}}.$$
\end{Proposition}
Note that the conclusion of this proposition gives exactly the information on $\vu$ and $\vb$ that was assumed in the first point of Proposition \ref{Proposition_Principale1}. However, although we have some information on the pressure $P$, this is not enough to obtain the fourth hypothesis of  Proposition \ref{Proposition_Principale1}. This term will be studied in detail in Section \ref{Secc_LastHypotheses} below. \\[5mm]
{\bf Proof of Proposition \ref{Propo_indu_argu}.}
Recall that from the global hypothesis of Theorem \ref{Teorem1}, we have a local control over the set $\Omega$ (see the set of hypotheses (\ref{HypothesesLocal1})), thus as we want to obtain a local information and since we assumed $Q_{4R_0}(t_0, x_0)\subset \Omega$, by the definition of Morrey spaces given in (\ref{Morreyparabolic}), we only need to prove that there exists a radius $R_3$ small enough such that for all $0<r< R_3$ and for all $(t,x) \in Q_{R_3}(t_0, x_0)$ we have the following controls
\begin{equation}\label{EstimationMorreyPoint1}
\iint_{Q_r (t,x) } |\vu(s,y)|^{3} + |\vb(s,y)|^3\, dy\,ds \leq C \, r^{5(1-\frac{3}{\tau_{0}})} \quad \text{and} \quad
\iint_{Q_r (t,x) } |P(s,y)|^{q_0}\,dy\, ds \leq C \, r^{5(1-\frac{2\q}{\tau_{0}})},
\end{equation}
indeed, for larger values of $r$ theses quantities will be controlled by the information over the set $\Omega$.\\

In order to obtain these estimates, we will implement an inductive argument  using the averaged quantities defined in (\ref{Definition_Averaged_Quantities}) and the inequalities (\ref{EstimationAveragedQuantities1}) obtained in the previous section. Indeed, in a first step, we remark that by Lemma \ref{Lemme_L3norm}, we can write
$$\iint_{Q_r (t,x) } |\vu(s,y)|^{3} + |\vb(s,y)|^3\,dy\,ds\leq C\, r^2(\mathcal{A}_r + \mathcal{B}_r + \alpha_r + \beta_r)^{\frac32}(t,x),$$
moreover, since we have the identity $r^{5-2q_0}\mathcal{P}_r(t,x)=\|P\|_{L_{t,x}^{q_0}  (Q_r)}^{q_0}$, we see that in order to obtain (\ref{EstimationMorreyPoint1}) for all small $0<r<R_3$ and all point $(t, x)\in Q_{R_3}$, it is enough to show the estimates
$$(\mathcal{A}_r + B_r + \alpha_r + \beta_r)(t,x)\leq Cr^{2(1-\frac{5}{\tau_{0}})} \quad \text{and} \quad \mathcal{P}_{r}(t,x)\leq C\, r^{2 \q (1-\frac{5}{\tau_{0}})}.$$
Let us now introduce the following quantities:
\begin{equation}\label{Def_MainGoal1}
\Ar (t,x)= \frac{1}{r^{2(1-\ftau)}} \left( \mathcal{A}_r + \mathcal{B}_r + \alpha_r + \beta_r \right) (t,x)
\quad \text{and} \quad \Qrq (t,x) = \frac{1}{r^{2\q(1-\ftau)}} \mathcal{P}_{r} (t,x),
\end{equation}
again, to prove (\ref{EstimationMorreyPoint1}) we only need to show that there exists $0<\kappa<1$ and $0<R_3 < 2R_0$ such that for all $n\in \mathbb{N}$ and $(t,x) \in Q_{R_3}$, we have
\begin{equation}\label{maingoal1}
\mathbf{A}_{\kappa^n R_3} (t,x) \leq C \quad \text{and} \quad
\mathbf{Q}_{\kappa^n R_3}  (t,x) \leq C,
\end{equation}
and the whole idea here is to use an inductive argument that ensures that we have these two previous estimates for all radii of the type $\kappa^n R_3>0$. This idea will be implemented in two steps by studying separately each one of the quantities of (\ref{maingoal1}).\\

In order to simplify the arguments, we shall also need the quantities:
\begin{equation}\label{Def_QuantitesenGRAS}
\BBr (t,x) = (\mathcal{B}_r + \beta_r )(t,x), \quad\Prq (t,x) = \frac{1}{r^{\q(1-\ftau)}} \mathcal{P}_{r} (t,x), \quad
\Dr (t,x) = \frac{1}{r^{3 -\frac{5}{\tau_c} }} \left(\mathcal{D}_r^{\frac{7}{10}} + \delta_r^{\frac{7}{10}} \right) (t,x).
\end{equation}
With these new quantities, we can rewrite the two inequalities of expression \eqref{EstimationAveragedQuantities1} as follows
\begin{equation}\label{Arexpression}
\Ar\leq C\left(\left(\frac{r}{\rho}\right)^{\frac{10}{\tau_{0}}} \Ao  
+\left(\frac{\rho}{r}\right)^{4- \frac{10}{\tau_{0}}} \Ao   \BBo^{\frac{1}{2}}  
+\left(\frac{\rho}{r}\right)^{4- \frac{10}{\tau_{0}}}  \Poq^{\frac{1}{\q}} \Ao^{\frac{1}{2}}+\left(\frac{\rho}{r} \right)^{3-\frac{10}{\tau_{0}}} \rho^{2+\frac{5}{\tau_0}-\frac{5}{\tau_c}} \Do \Ao^{\frac{1}{2}}\right),
\end{equation}
and 
\begin{equation}\label{Prexpression}
\Prq \leq C \left( \left(\frac{\rho}{r} \right)^{3-\frac{5\q}{\tau_{0}}} (\Ao \BBo)^{\frac{\q}{2}} + \left(\frac{r}{\rho} \right)^{q_0(1+\frac{5}{\tau_0})-2} \Poq \right).
\end{equation}
Observe that these two estimates essentially give us the estimate for $\vu, \vb$ and $P$ within the (small) parabolic ball $Q_r$ in terms of $\vu, \vb$ and $P$ within the (larger) parabolic ball $Q_\rho$.\\ 

We define now a new expression that will help us to set up the inductive argument:
\begin{equation}\label{thetaexpression}
\Tr (t,x) = \Ar (t,x) + \left(\kappa^{5\q (\frac{2}{\tau_{0}} - 1)}\Prq (t,x) \right)^{\frac{2}{q_0}} \quad \text{with} \quad \kappa = \frac{r}{\rho}<1,
\end{equation}
and we will see how to obtain from (\ref{Arexpression}) and (\ref{Prexpression}) a recursive equation in terms of $\Tr $ from which we will deduce (\ref{maingoal1}). Indeed, we have the following lemma:
\begin{Lemme}\label{Lemme_InegaliteIteration}
For all $(t,x)\in Q_{2R_0}(t_0,x_0)$, for all $0<r<\frac{\rho}{2}$ and for $\rho$ small enough we have the inequality
\begin{equation}\label{iterationform}
\Tr(t,x) \leq \frac{1}{2} \TTo (t,x) + \epsilon, 
\end{equation}
where $\epsilon$ is a small constant that depends on the information available on the external forces $\vf$ and $\vg$ through the quantity $\Dr$ given in (\ref{Def_QuantitesenGRAS}).
\end{Lemme}
It is worth noting here that since $0<r<\rho$ this inequality expresses a control of the quantities $\mathcal{A}_r, B_r, \alpha_r, \beta_r$ and $\mathcal{P}_r$ on small domains through the information on larger domains.\\[5mm]
{\bf Proof of Lemma \ref{Lemme_InegaliteIteration}.} As announced, this inequality relies on the controls (\ref{Arexpression}) and (\ref{Prexpression}) obtained previously. 
In order to construct $\Tr$ we first multiply expression (\ref{Prexpression}) by 
$\kappa^{5\q (\frac{2}{\tau_{0}}-1)}$, we take the $\frac{2}{q_0}$-power of it and then we sum the resulting inequality to (\ref{Arexpression}) and we obtain (recall that $\kappa = \frac{r}{\rho}$):
\begin{eqnarray*}
\Tr&=&\Ar+\left(\kappa^{5\q (\frac{2}{\tau_{0}} - 1)}\Prq (t,x) \right)^{\frac{2}{q_0}}\\
&\leq & C\left(\kappa^{\frac{10}{\tau_{0}}} \Ao +\kappa^{\frac{10}{\tau_{0}}-4} \Ao   \BBo^{\frac{1}{2}} +\kappa^{\frac{10}{\tau_{0}}-4}  \Poq^{\frac{1}{\q}} \Ao^{\frac{1}{2}}+\kappa^{\frac{10}{\tau_{0}}-3} \rho^{2+\frac{5}{\tau_0}-\frac{5}{\tau_c}} \Do \Ao^{\frac{1}{2}}\right)\\
&&+ C \left(\kappa^{\frac{15q_0}{\tau_{0}}-5q_0 - 3} (\Ao \BBo)^{\frac{\q}{2}} + \kappa^{\frac{15q_0}{\tau_0}-4q_0-2} \Poq \right)^{\frac{2}{q_0}}.
\end{eqnarray*}
As it is clear from the definition of $\Tr$ given in (\ref{thetaexpression}) that we have $\Ar\leq \Tr$, we can write
\begin{eqnarray}
\Tr&\leq & C\left(\kappa^{\frac{10}{\tau_{0}}}  \TTo    
+\kappa^{\frac{10}{\tau_{0}}-4} \BBo^{\frac{1}{2}}   \TTo  
+\underbrace{\kappa^{\frac{10}{\tau_{0}}-4}  \Poq^{\frac{1}{\q}} \Ao^{\frac{1}{2}}}_{(I)}+\underbrace{\kappa^{\frac{10}{\tau_{0}}-3} \rho^{2+\frac{5}{\tau_0}-\frac{5}{\tau_c}} \Do \Ao^{\frac{1}{2}}}_{(II)}\right)\notag\\
&&+ C \underbrace{\left(\kappa^{\frac{15q_0}{\tau_{0}}-5q_0 - 3} (\Ao \BBo)^{\frac{\q}{2}} + \kappa^{\frac{15q_0}{\tau_0}-4q_0-2} \Poq \right)^{\frac{2}{q_0}}}_{(III)}.\label{iterationform0}
\end{eqnarray}
We now study the terms $(I)$, $(II)$ and $(III)$. The first one is easy to handle since we have
\begin{eqnarray}
\kappa^{\frac{10}{\tau_{0}}-4} \Poq^{\frac{1}{\q}} \Ao^{\frac{1}{2}}&=&\kappa^{\frac{10}{\tau_{0}}-4} \left(\kappa^{5(\frac{1}{\tau_{0}} -\frac{1}{2})} \Poq^{\frac{1}{\q}} \times\kappa^{5(\frac{1}{2} - \frac{1}{\tau_{0}} )} \Ao^{\frac{1}{2}}\right)\leq\kappa^{\frac{10}{\tau_{0}}-4} \left( \kappa^{10(\frac{1}{2} - \frac{1}{\tau_{0}} )} \Ao + \kappa^{10(\frac{1}{\tau_{0}} - \frac{1}{2} )} \Poq^{\frac{2}{\q}} \right)\notag \\
&\leq & \kappa \left(  \Ao + \left(\kappa^{5\q (\frac{2}{\tau_{0}} - 1)}\Poq \right)^{\frac{2}{\q}} \right) \leq \kappa  \TTo.\label{Estimation(I)Iteration}
\end{eqnarray}
For the term $(II)$ of (\ref{iterationform0}) we simply write
\begin{equation}\label{Estimation(II)Iteration}
\kappa^{\frac{10}{\tau_{0}}-3} \rho^{2+\frac{5}{\tau_0}-\frac{5}{\tau_c}} \Do \Ao^{\frac{1}{2}}\leq \kappa^{\frac{10}{\tau_{0}}-3} \rho^{2+\frac{5}{\tau_0}-\frac{5}{\tau_c}} ( \Do^2+\Ao)\leq \kappa^{\frac{10}{\tau_{0}}-3} \rho^{2+\frac{5}{\tau_0}-\frac{5}{\tau_c}} ( \Do^2+\TTo).
\end{equation}
The last term $(III)$ of (\ref{iterationform0}) is treated in the following way.
\begin{eqnarray}
\left(\kappa^{\frac{15q_0}{\tau_{0}}-5q_0 - 3} (\Ao \BBo)^{\frac{\q}{2}} + \kappa^{\frac{15q_0}{\tau_0}-4q_0-2} \Poq \right)^{\frac{2}{q_0}}&\leq &C\left(\kappa^{\frac{30}{\tau_{0}}-10-\frac{6}{\q}}\Ao \BBo +\kappa^{\frac{10}{\tau_0}+2-\frac{4}{q_0}}\left(\kappa^{5\q (\frac{2}{\tau_{0}} - 1)}\Poq (t,x) \right)^{\frac{2}{q_0}} \right)\notag\\
&\leq & C\kappa^{\frac{30}{\tau_{0}}-10-\frac{6}{\q}} \BBo\TTo + C\kappa^{\frac{10}{\tau_0}+2-\frac{4}{q_0}}\TTo.\label{Estimation(III)Iteration}
\end{eqnarray}
Plugging estimates (\ref{Estimation(I)Iteration}), (\ref{Estimation(II)Iteration}) and (\ref{Estimation(III)Iteration}) in inequality (\ref{iterationform0}) we obtain
\begin{eqnarray}
\Tr&\leq & C\left(\kappa^{\frac{10}{\tau_{0}}}  \TTo    
+\kappa^{\frac{10}{\tau_{0}}-4} \BBo^{\frac{1}{2}} \TTo  
+\kappa \TTo  + \kappa^{\frac{10}{\tau_{0}}-3} \rho^{2+\frac{5}{\tau_0}-\frac{5}{\tau_c}}\Do^2+ \kappa^{\frac{10}{\tau_{0}}-3} \rho^{2+\frac{5}{\tau_0}-\frac{5}{\tau_c}}\TTo\right.\notag\\
& &+\left. \kappa^{\frac{30}{\tau_{0}}-10-\frac{6}{\q}} \BBo\TTo + \kappa^{\frac{10}{\tau_0}+2-\frac{4}{q_0}}\TTo\right)\notag\\
&\leq &C \left(\kappa^{\frac{10}{\tau_{0}}}+\kappa^{\frac{10}{\tau_{0}}-4} \BBo^{\frac{1}{2}}+\kappa+ \kappa^{\frac{10}{\tau_{0}}-3}\rho^{2+\frac{5}{\tau_0}-\frac{5}{\tau_c}}+\kappa^{\frac{30}{\tau_{0}}-10-\frac{6}{\q}} \BBo+\kappa^{\frac{10}{\tau_0}+2-\frac{4}{q_0}}\right)\TTo\notag\\
&&+ C \kappa^{\frac{10}{\tau_{0}}-3} \rho^{2+\frac{5}{\tau_0}-\frac{5}{\tau_c}}\Do^2.\label{iterationform2}
\end{eqnarray}
Now we want to fix a small $\kappa$ such that (\ref{iterationform}) can be deduced from this inequality: we thus want to prove
\begin{eqnarray}
C \left(\kappa^{\frac{10}{\tau_{0}}}+\kappa^{\frac{10}{\tau_{0}}-4} \BBo^{\frac{1}{2}}+\kappa+ \kappa^{\frac{10}{\tau_{0}}-3}\rho^{2+\frac{5}{\tau_0}-\frac{5}{\tau_c}}+\kappa^{\frac{30}{\tau_{0}}-10-\frac{6}{\q}} \BBo+\kappa^{\frac{10}{\tau_0}+2-\frac{4}{q_0}}\right)&\leq &\frac{1}{2}\quad \mbox{and}\notag\\
C \kappa^{\frac{10}{\tau_{0}}-3} \rho^{2+\frac{5}{\tau_0}-\frac{5}{\tau_c}}\Do^2&\leq & \epsilon.\label{EstimationPourInteration1}
\end{eqnarray}
At this point we remark that due to the hypothesis (\ref{smallgrad}) and to the definition of quantities $\mathcal{B}_\rho$ and $\beta_\rho$ given in (\ref{Definition_Averaged_Quantities}) we have 
$$\underset{\rho\to 0}{\lim \sup}\;\BBo=\underset{\rho\to 0}{\lim \sup}(\mathcal{B}_\rho+\beta_\rho)<2\epsilon^*,$$
and thus, although we have $0<\kappa<1$ and $\frac{10}{\tau_{0}}-4<0$ and $\frac{30}{\tau_{0}}-10-\frac{6}{\q}<0$, then the following terms of (\ref{iterationform2})
\begin{equation}\label{ConditionVoisinage1}
\kappa^{\frac{10}{\tau_{0}}-4} \BBo^{\frac{1}{2}}\qquad \mbox{and}\qquad \kappa^{\frac{30}{\tau_{0}}-10-\frac{6}{\q}} \BBo,
\end{equation}
can be made very small if $\rho$ is small enough. Moreover, since by hypothesis we have $2+\frac{5}{\tau_0}-\frac{5}{\tau_c}>0$, then the fourth term of (\ref{iterationform2}) 
$\kappa^{\frac{10}{\tau_{0}}-3}\rho^{2+\frac{5}{\tau_0}-\frac{5}{\tau_c}}$ can also be made small. Finally we observe that $\frac{10}{\tau_0}+2-\frac{4}{q_0}>0$ which is equivalent to $\frac{4\tau_0}{10+2\tau_0}<q_0$, but since $1<q_0<\frac{3}{2}$ and $\frac{5}{1-\alpha}<\tau_0<5q_0$, then this condition is fulfilled and the term $\kappa^{\frac{10}{\tau_0}+2-\frac{4}{q_0}}$ can be small if $\kappa$ is small. We also observe that by (\ref{Def_QuantitesenGRAS}) and by the definition of the averaged quantities $\mathcal{D}_\rho$ and $\delta_\rho$ given in (\ref{Definition_Averaged_Quantities}) we have that $\Do$ is bounded and controlled by the Morrey norms of the external forces $\vf$ and $\vg$. Thus since $2+\frac{5}{\tau_0}-\frac{5}{\tau_c}>0$ the term
\begin{equation}\label{ConditionVoisinage2} 
\kappa^{\frac{10}{\tau_{0}}-3} \rho^{2+\frac{5}{\tau_0}-\frac{5}{\tau_c}}\Do^2,
\end{equation}
can also be made very small if $\rho$ is small enough. With all these observations we have the controls (\ref{EstimationPourInteration1}) and the inequality (\ref{iterationform2}) can thus be rewritten in the following form:
$$\Tr (t,x) \leq \frac{1}{2} \TTo (t,x) + \epsilon,$$
and Lemma \ref{Lemme_InegaliteIteration} is proven. \hfill $\blacksquare$
\begin{Remarque}
The terms (\ref{ConditionVoisinage1}) and (\ref{ConditionVoisinage2}) can be made small if $R_0$ is small enough and we can see here that all the techniques displayed here will only be valid on a small neighborhood of the point $(t_0, x_0)$. 
\end{Remarque}
Now, we turn to the proof of Proposition \ref{Propo_indu_argu}. 
With the inequality \eqref{iterationform} at hand we can obtain the first estimate of \eqref{maingoal1}. We first consider the estimates centered at the fixed point $(t_0, x_0)$. Indeed, notice that for any radius $\rho$ such that $0<\rho < 2R_0$ and since we have $Q_{4R_0}(t_0,x_0)\subset \Omega$ (recall formula (\ref{Defi_BouleParaboliquedeBase})), by the hypotheses given in (\ref{HypothesesLocal1}) we have the bounds: 
\begin{equation*}
\begin{split}
\|\vu\|_{L^\infty_tL^2_x(Q_{\rho}(t_0,x_0))}\leq \|\vu\|_{L^\infty_tL^2_x(\Omega)}<+\infty,
\qquad \|\vn\otimes\vu\|_{L_{t,x}^{2}(Q_{\rho}(t_0,x_0))}\leq \|\vn\otimes\vu\|_{L_{t,x}^{2}(\Omega)}<+\infty,\\
\text{and}\qquad
\|P\|_{L_{t,x}^{q_0}(Q_{\rho}(t_0,x_0))}\leq \|P\|_{L_{t,x}^{q_0}(\Omega)}<+\infty.\end{split}
\end{equation*}
Now, by the definition of the quantities $\mathcal{A}_\rho(t_0, x_0)$, $\mathcal{B}_\rho(t_0, x_0)$, $\alpha_\rho(t_0, x_0)$, $\beta_\rho(t_0, x_0)$ and $\mathcal{P}_\rho (t_0, x_0)$ given in (\ref{Definition_Averaged_Quantities}) we have 
\begin{eqnarray*}
\rho\mathcal{A}_\rho (t_0, x_0)=\|\vu\|_{L^\infty_tL^2_x(Q_{\rho}(t_0,x_0))}^2,\qquad \rho\mathcal{B}_\rho(t_0, x_0)=\|\vn\otimes \vu\|_{L_{t,x}^{2} (Q_{\rho}(t_0,x_0))}^2,\\
\rho\alpha_\rho (t_0, x_0)=\|\vb\|_{L^\infty_tL^2_x(Q_{\rho}(t_0,x_0))}^2, \qquad \rho\beta_\rho (t_0, x_0)=\|\vn\otimes \vb\|_{L_{t,x}^{2}(Q_{\rho}(t_0,x_0))}^2,\\
\text{and}\qquad \rho^{5-2q_0}\mathcal{P}_\rho(t_0, x_0)=\|P\|_{L_{t,x}^{q_0}(Q_{\rho}(t_0,x_0))}, \hspace{6cm}
\end{eqnarray*}
and thus we have the following uniform bounds
$$\underset{0<\rho<2R_0}{\sup}\Big\{\rho\mathcal{A}_\rho(t_0,x_0),\;\rho\mathcal{B}_\rho(t_0,x_0),\; \rho\alpha_\rho(t_0,x_0), \;\rho\beta_\rho(t_0,x_0),\; \rho^{5-2q_0}\mathcal{P}_\rho(t_0, x_0) \Big\}<+\infty,$$
from which we can deduce, by the definition of the quantities $\mathbf{A}_\rho(t_0, x_0)$ given in (\ref{Def_MainGoal1}) and $\mathbf{P}_\rho(t_0, x_0)$ given in (\ref{Def_QuantitesenGRAS}), the uniform bounds
\begin{equation}\label{UniformBoundsInitialisation}
\underset{0<\rho<2R_0}{\sup}\; \rho^{3-\frac{10}{\tau_0}}\mathbf{A}_\rho(t_0, x_0)<+\infty \qquad\text{and}\qquad \underset{0<\rho<2R_0}{\sup}\; \rho^{5-q_0(1+\frac{5}{\tau_0})}\mathbf{P}_\rho(t_0, x_0)<+\infty.
\end{equation}
Note now, that there exists a $0<\kappa<\tfrac12$ and a fixed $0<\rho_0 <2R_0$ small such that, on one hand, by (\ref{UniformBoundsInitialisation}) the quantities $\mathbf{A}_{\rho_0}(t_0, x_0)$ and $\mathbf{P}_{\rho_0}(t_0, x_0)$ are bounded (then the quantity $\mathbf{\Theta}_{\rho_0}$ defined by expression  (\ref{thetaexpression}) is itself bounded) and, on the other hand, if  $\rho_0$ is small enough, then the inequality \eqref{iterationform} in Lemma \ref{Lemme_InegaliteIteration} holds true and we can write
$$\mathbf{\Theta}_{\kappa\rho_0}(t_0, x_0)\leq \frac{1}{2}\mathbf{\Theta}_{\rho_0}(t_0, x_0)+\epsilon.$$
We can iterate this process and we obtain for all $n\geq 1$:
$$\mathbf{\Theta}_{\kappa^n\rho_0}(t_0, x_0)\leq \frac{1}{2^n}\mathbf{\Theta}_{\rho_0}(t_0, x_0)+\epsilon\sum_{j=0}^{n-1}2^{-j},$$
therefore, there exists $N\geq 1$ such that for all $n\geq N$ we have 
\begin{equation}\label{EstimationThetaBorne1}
\mathbf{\Theta}_{\kappa^n\rho_0}(t_0, x_0)\leq 4\epsilon,
\end{equation}
from which we obtain (from formula (\ref{thetaexpression})) that 
\begin{equation}\label{EstimationARBorne1}
\mathbf{A}_{\kappa^N\rho_0}(t_0, x_0)\leq \frac{1}{8} C \quad \text{and} \quad \mathbf{P}_{\kappa^N\rho_0}(t_0, x_0)\leq \frac{1}{32} C.
\end{equation}
This information is centered at the point $(t_0,x_0)$, in order to treat the uncentered bound, we can let $\frac{1}{2}\kappa^N\rho_0$ to be the radius $R_3$ we want to find, thus for all points $(t,x) \in Q_{R_3}$ we have $Q_{R_3} (t,x) \subset Q_{2R_3}$, which implies   
$$\mathbf{A}_{R_3} (t,x)  \leq 2^{3 - \frac{10}{\tau_{0}}}\mathbf{A}_{2R_3} (t_0,x_0) \leq 8 \,\mathbf{A}_{2R_3} (t_0,x_0) 
\leq 8 \,\mathbf{A}_{\kappa^N\rho_0} (t_0,x_0)
<C,$$ 
and
$$\mathbf{P}_{R_3} (t,x)  \leq 2^{5-q_0(1+\frac{5}{\tau_{0}})}\mathbf{P}_{2R_3} (t_0,x_0) \leq 32 \,\mathbf{P}_{2R_3} (t_0,x_0) 
\leq 32 \,\mathbf{P}_{\kappa^N\rho_0} (t_0,x_0)
<C,$$
by definition of $\mathbf{\Theta}_{R_3}$, we thus get $\mathbf{\Theta}_{R_3} (t,x)\leq C$. Applying Lemma \ref{Lemme_InegaliteIteration} and iterating once more, we find that the same will be true for $\kappa R_3$ and then for all $\kappa^n R_3$, $n\in \mathbb{N}$,i.e. 
$$\mathbf{A}_{\kappa^n R_3} (t,x)\leq C, \quad \text{for all}  \quad n\in \mathbb{N}  \quad \text{and} \quad (t,x) \in Q_{R_3}$$
and the first inequality of (\ref{maingoal1}) is proven.\\

The second inequality of (\ref{maingoal1}) requires a different treatment since from (\ref{EstimationThetaBorne1}) and by the definition of the quantity $\mathbf{Q}_{\kappa^n\rho_0}$  given in (\ref{Def_MainGoal1}), we can only deduce that for all $n\geq N$ we have the bound
$$\kappa^{nq_0(1-\frac{5}{\tau_0})+5q_0(\frac{2}{\tau_0}-1)}\mathbf{Q}_{\kappa^n\rho_0}(t_0, x_0)\leq C,$$
which is not enough to ensure that the quantity $\mathbf{Q}_{\kappa^n\rho_0}(t_0, x_0)$ is bounded (since $q_0(1-\frac{5}{\tau_0})<0$).\\ 
To overcome this issue, using the definition of the quantities $\mathbf{A}_r$ and $\mathbf{Q}_{r}$ given in (\ref{Def_MainGoal1}) and using the estimate (\ref{pestimscalq0}) we can write (recalling that $\kappa=\frac{r}{\rho}$):
\begin{eqnarray}
\Qrq (t_0,x_0)&\leq &C\left( \left(\frac{\rho}{r} \right)^{3+ \q (1-\frac{10}{\tau_{0}})} \Ao^{\q} (t_0,x_0) +  \left(\frac{r}{\rho} \right)^{2 (\frac{5 \q}{\tau_{0}} -1)} \mathbf{Q}_{\rho} (t_0,x_0) \right)\notag\\
&\leq& C\left( \kappa^{-3-\q (1-\frac{10}{\tau_{0}})} \Ao^{\q} (t_0,x_0) +\kappa^{2 (\frac{5 \q}{\tau_{0}} -1)} \mathbf{Q}_{\rho} (t_0,x_0) \right).\label{Qrexpression}
\end{eqnarray}
We need to impose a smallness condition on $0<\tilde{\kappa}<1$ and we will assume that we have
$$C\, \tilde{\kappa}^{2 (\frac{5\q}{\tau_{0}} -1)} < \frac{1}{2},$$
which is possible since $\frac{5}{1-\alpha}<\tau_{0}<5q_0$. Now, from (\ref{EstimationARBorne1}) we know that the quantity $\Ao(t_0, x_0)$ can be made small enough if $\rho$ is small enough, and thus the estimate \eqref{Qrexpression} becomes 
\begin{equation}\label{iterationQ}
\Qrq (t_0,x_0) \leq \frac{1}{2}  \mathbf{Q}_{\rho} (t_0,x_0) + \tilde{\epsilon},
\end{equation}
where $\tilde{\epsilon}$ is a small constant.\\

As before, for a small $0<\tilde{\kappa}<\tfrac12$ and for a fixed small $0<\tilde{\rho_0} <2R_0$, we have by (\ref{UniformBoundsInitialisation}) that the quantity $\mathbf{Q}_{\tilde{\rho_0}} (t_0,x_0)$ is bounded  and if $\tilde{\kappa}$ and $\tilde{\rho_0}$ are small enough then inequality  \eqref{iterationQ} holds true, \emph{i.e.}, 
\begin{equation*}
\mathbf{Q}_{\tilde{\kappa}\tilde{\rho_0}}(t_0, x_0)\leq \frac{1}{2}\mathbf{Q}_{\tilde{\rho_0}}(t_0, x_0)+\tilde{\epsilon},
\end{equation*}
Iterating the inequality above, we obtain that for all $n\geq 1$
$$\mathbf{Q}_{\tilde{\kappa}^n\tilde{\rho_0}}(t_0, x_0)\leq \frac{1}{2^n}\mathbf{Q}_{\tilde{\rho_0}}(t_0, x_0)+\tilde{\epsilon}\sum_{j=0}^{n-1}2^{-j},$$
and there exists $\tilde{N}\geq 1$ such that for all $n\geq \tilde{N}$ we have 
$$\mathbf{Q}_{\tilde{\kappa}^n\tilde{\rho_0}}(t_0, x_0)\leq 4\tilde{\epsilon} \leq \frac{1}{32} C.$$
In order to treat the uncentered bound, we proceed as before: let $\frac{1}{2}\tilde{\kappa}^{\tilde{N}}\tilde{\rho_0}$ be the radius $R_3$, thus for all points $(t,x) \in Q_{R_3}$ we have $Q_{R_3} (t,x) \subset Q_{2R_3}$, which implies   
$$\mathbf{Q}_{R_3} (t,x)  \leq 2^{(5 - \frac{10q_0}{\tau_{0}})}\mathbf{Q}_{2R_3} (t_0,x_0) \leq 32 \,\mathbf{Q}_{2 R_3} (t_0,x_0) \leq \mathbf{Q}_{\tilde{\kappa}^n\tilde{\rho_0}}(t_0, x_0) \leq C,$$  
and the second inequality of (\ref{maingoal1}) is now proved since the inequality above holds true for all $\kappa^n R_3$, $n\in \mathbb{N}$. We conclude the proof of Proposition \ref{Propo_indu_argu} by choosing $R_3 = \max \{\frac{1}{2}\kappa^{N}\rho_0,  \frac{1}{2}\tilde{\kappa}^{\tilde{N}}\tilde{\rho_0}\}$.\hfill$\blacksquare$\\

From the proof of Proposition \ref{Propo_indu_argu}, we can deduce a more specific result on $\vn \otimes \vu$ and $\vn \otimes \vb$. Indeed, we can obtain the following result that gives the assumption \emph{2)} of Proposition \ref{Proposition_Principale1}:
\begin{Corollaire}\label{cor-indu-argu}
Under all the assumptions of Proposition \ref{Propo_indu_argu}, we have 
$$\mathds{1}_{Q_{R_3}} \vn \otimes \vu \in \mathcal{M}_{t,x}^{2, \tau_{1}} \quad \text{and} \quad \mathds{1}_{Q_{R_3}}  \vn \otimes \vb \in \mathcal{M}_{t,x}^{2, \tau_{1}},$$ 
with $\frac{1}{\tau_{1}} = \frac{1}{\tau_{0}} + \frac{1}{5}$.
\end{Corollaire}
{\bf Proof. } From the definition of $\mathbf{A}_r$ in (\ref{Def_MainGoal1}) and from the first estimate of (\ref{maingoal1}), in the proof of Proposition \ref{Propo_indu_argu}, for all $0<r$ small and for all $(t,x) \in Q_{R_3}(t_0, x_0)$, we have shown that we have
$$( B_r + \beta_r)(t,x)=\frac{1}{r}  \iint_{Q_r(t,x)} |\vn \otimes \vu(s,y)|^2dyds + \frac{1}{r}  \iint_{Q_r(t,x)} |\vn \otimes \vb(s,y)|^2dyds \leq C \,r^{2(1-\frac{5}{\tau_{0}})}=C\,r^{4 - \frac{10}{\tau_{1}}}, $$
where we used the relationship $\frac{1}{\tau_{1}} = \frac{1}{\tau_{0}} + \frac{1}{5}$. We obtain then, for all $0<r$ small, the estimate 
$$ \iint_{Q_r(t,x)} |\vn \otimes \vu(s,y)|^2dyds + \iint_{Q_r(t,x)} |\vn \otimes \vb(s,y)|^2dyds \leq  C\,r^{5(1 - \frac{2}{\tau_{1}})},$$
and to conclude it is enough to recall the definition of the Morrey space $\mathcal{M}_{t,x}^{2, \tau_{1}}$ given in (\ref{Morreyparabolic}). \hfill$\blacksquare$
\section{Further estimates}\label{Secc_LastHypotheses}
In the previous sections we have proven so far the points \emph{1)}, \emph{2)} and \emph{5)} that were assumed in Proposition \ref{Proposition_Principale1} (for a small $R_3<2R_0$):
\begin{equation}\label{Estimations135}
\begin{split}
&\mathds{1}_{Q_{R_3}}\vu,\;\mathds{1}_{Q_{R_3}} \vb\in \mathcal{M}_{t,x}^{3, \tau_{0}}\quad \mbox{for some} \quad\tau_{0} > \frac{5}{1-\alpha},\\[2mm]
&\mathds{1}_{Q_{R_3}} \vn \otimes \vu, \;\mathds{1}_{Q_{R_3}}\vn \otimes \vb \in \mathcal{M}_{t,x}^{2, \tau_{1}}\quad \mbox{with} \quad\frac{1}{\tau_{1}} = \frac{1}{\tau_{0}} + \frac{1}{5},\\[2mm]
&\mathds{1}_{Q_{R_3}} \vf\in \mathcal{M}_{t,x}^{\frac{10}{7},\tau_{a}},\quad\mathds{1}_{Q_{R_3}}\vg \in \mathcal{M}_{t,x}^{\frac{10}{7},\tau_{b}}\quad \mbox{for some}\quad \tau_a,\tau_b > \frac{5}{2-\alpha}.
\end{split}
\end{equation}
Our current task consists in proving the remaining points \emph{3)} and \emph{4)} using all the information available up to now, \emph{i.e.} we need to study the following assertions (with $R_0<R_1<R_2<R_3$): \\
\begin{eqnarray}
&&\mathds{1}_{Q_{R_2}}\vu, \;\mathds{1}_{Q_{R_2}} \vb  \in \mathcal{M}_{t,x}^{3, \delta}\quad\mbox{with} \quad\frac{1}{\delta} + \frac{1}{\tau_{0}} \leq \frac{1-\alpha}{5},\label{Estimations24}\\[2mm]
&&\mbox{for}\; 1\leq i,j\leq 3 \; \mbox{we have}\; \mathds{1}_{Q_{R_1}} \frac{\vn\partial_i \partial_j}{(-\Delta)}(u_ib_j)\in\mathcal{M}_{t,x}^{\mathfrak{p}, \mathfrak{q}}\quad\mbox{with}\; \mathfrak{p}_0\leq \mathfrak{p}<+\infty,\; \mathfrak{q}_0\leq \mathfrak{q}<+\infty,\label{Estimations241}
\end{eqnarray}
where $1 \leq \mathfrak{p}_0\leq \tfrac{6}{5}$ and $\tfrac52 <\mathfrak{q}_0< 3$ where $\frac{1}{\mathfrak{q}_0}=\frac{2-\alpha}{5}$ with $0< \alpha < \frac13$.\\

These two points are actually related and in order to study them we need to recall some tools of harmonic analysis in the setting of parabolic spaces. 
Let us now introduce, for $0<\mathfrak{a} < 5$, the parabolic Riesz potential $\mathcal{I}_{\mathfrak{a}}$ of a locally integrable function $\vf:\mathbb{R} \times \R\longrightarrow \mathbb{R}^3$  which is given by the expression
\begin{equation}\label{Def_Riesz_potential}
\mathcal{I}_{\mathfrak{a}} (\vf) (t,x) =\int_{\mathbb{R}} \int_{\R} \frac{1}{(\vert t-s\vert^{\frac{1}{2}} + \vert x-y \vert)^{5-\mathfrak{a}}} \vf (s,y) dy\, ds.
\end{equation}
As for the standard Riesz Potential in $\mathbb{R}^3$, we have a corresponding  boundedness property:
\begin{Lemme}[Adams-Hedberg's inequality]\label{Lemme_Hed}
If $0 <\mathfrak{a} < \frac{5}{q}$, $1 < p \leq q < +\infty$ and $\vf \in \mathcal{M}_{t,x}^{p,q} (\mathbb{R} \times \R)$ then for $\lambda=1- \frac{\mathfrak{a} q}{5}$, we have the inequality
$$\|\mathcal{I}_\mathfrak{a}(\vf)\|_{\mathcal{M}_{t,x}^{\frac{p}{\lambda} , \frac{q}{\lambda}}} \leq C\|\vf\|_{\mathcal{M}_{t,x}^{p,q}}.$$
\end{Lemme}
See \cite{Adams} for a proof of this fact.\\

We can state now the main proposition of this section which will focus on the information (\ref{Estimations24}):
\begin{Proposition}\label{Propo_Points24}
Let $(\vu, p, \vb)$ be a suitable solution of MHD equations (\ref{EquationMHD}) over $\Omega$ in the sense of Definition \ref{Def_SuitableSolutions}. Assume the general hypotheses (\ref{HypothesesLocal1}) and assume moreover the local informations (\ref{Estimations135}) for a parabolic ball $Q_{R_3}$. Then for some $R_2$ such that $R_0<R_2<R_3$ we have 
$$\mathds{1}_{Q_{R_2}} \vu \in \mathcal{M}_{t,x}^{3, \delta}\quad \mbox{and} \quad\mathds{1}_{Q_{R_2}} \vb \in \mathcal{M}_{t,x}^{3, \delta},$$ 
with $\frac{1}{\delta}+\frac{1}{\tau_{0}} < \frac{1-\alpha}{5}$.
\end{Proposition}
{\bf Proof.} For a point $(t_0, x_0)$ that satisfies the hypothesis (\ref{HypothesePetitesseGrad}), consider the following radii 
$$0<R_0<R_2<\bar{R}<\widetilde{R}<R_3<2R_0<t_0,$$ 
and the corresponding parabolic balls (recall formula (\ref{Def_ParabolicBalldeBase}))
$$Q_{R_2}(t_0,x_0)\subset Q_{\bar{R}}(t_0,x_0)\subset Q_{\widetilde{R}}(t_0,x_0)\subset Q_{R_3}(t_0,x_0) \subset Q_{2R_0}(t_0,x_0).$$ 
We introduce now two test functions $\bar\phi,\bar\varphi:\mathbb{R}\times \R\longrightarrow \mathbb{R}$ that belong to the space $\mathcal{C}^{\infty}_{0}(\mathbb{R}\times \R)$ and such that 
\begin{eqnarray}
\bar\phi\equiv 1 \; \text{on} \; Q_{R_2} \quad \text{and} \quad \text{supp}(\bar\phi)\subset Q_{\bar R},\label{DefSoporteFuncTest1}\\[2mm]
\bar\varphi\equiv 1 \; \text{on}\; Q_{\widetilde{R}} \quad  \text{and} \quad \text{supp}(\bar\varphi) \subset  Q_{R_3}.\label{DefSoporteFuncTest2}
\end{eqnarray}
Note that since $R_3<2R_0<t_0$ we have $\bar\phi(0,\cdot)=\bar\varphi (0, \cdot)=0$ and remark that we have by construction the identity $\bar\phi \bar\varphi\equiv\bar\phi$. We define the variable $\vVV$ by the expression
\begin{equation*}
\vVV= \bar\phi (\vu+ \vb),
\end{equation*}
and if we study the equation satisfied by $\vVV$ we obtain
\begin{equation}\label{Equation_ChaleurLocale1}
\begin{cases}
\partial_t\vVV(t,x)=\Delta \vVV(t,x)+\NN(t,x),\\[3mm]
\vVV(0, x)=0,
\end{cases}
\end{equation}
where 
\begin{equation}\label{Definition_TermeNN}
\NN=(\partial_t \bar\phi-\Delta \bar\phi) (\vu+\vb)-2\sum_{i=1}^{3}(\partial_{i}\bar\phi) (\partial_i(\vu+\vb))-\bar\phi \left( (\vb\cdot\vn)\vu+(\vu\cdot\vn)\vb\right)-2 \bar\phi(\vn P)+ \bar\phi (\vf+\vg).
\end{equation}
Although this problem is very similar to the one studied with the variable $\vUU$ defined in (\ref{Def_VariableU}) which satisfies equation (\ref{Equation_ChaleurLocale}), we will perform different computations in order to obtain the conclusion of Proposition \ref{Propo_Points24}. The main point is to express the pressure $P$ in a very specific manner, indeed, since  $P = \bar\varphi P$ on the cylinder $Q_{\widetilde{R}}$ (see (\ref{DefSoporteFuncTest2})), then over the parabolic ball $Q_{R_2}$ we have the identity 
$$\displaystyle{- \Delta (\bar\varphi P)=-\bar\varphi \Delta P + (\Delta \bar\varphi)P-2\sum^3_{i= 1}\partial_i ( (\partial_i \bar\varphi) P)},$$ 
from which we deduce the formula
\begin{equation}\label{Definition_LocalPression}
\bar\phi (\vn P) =\bar\phi \frac{\vn \big(-\bar\varphi \Delta P\big)}{(-\Delta)} + \bar\phi\frac{\vn\big((\Delta \bar\varphi)P\big)}{(-\Delta)} -2\sum^3_{i= 1}\bar\phi\frac{\vn\big(\partial_i ( (\partial_i \bar\varphi) P )\big)}{(-\Delta)}.
\end{equation}
Recalling that we have the identity $\displaystyle{\Delta P=-\sum^3_{i,j= 1}\partial_i \partial_j (u_i b_j)}$, then the first term of (\ref{Definition_LocalPression}) can be rewritten in the following manner:
\begin{eqnarray}
\bar\phi \frac{\vn \big(-\bar\varphi \Delta P\big)}{(-\Delta)}
&= &\bar\phi \frac{\vn}{(- \Delta )}\Big(\bar\varphi \,  \sum^3_{i,j= 1}\partial_i \partial_j  (u_i b_j)  \Big) \notag\\
&= &\sum^3_{i,j= 1} \bar\phi \frac{\vn}{(- \Delta )} \Big(\partial_i  \partial_j (\bar\varphi u_i b_j ) - \partial_i  \big((\partial_j \bar\varphi) u_i b_j \big) - \partial_j  \big((\partial_i \bar\varphi) u_i b_j \big) + (\partial_i \partial_j  \bar\varphi) (u_i b_j)  \Big),\label{Identite_Pression1}
\end{eqnarray}
note that the first term of the right-hand side above satisfies the identity 
\begin{equation}\label{Identite_Pression2}
\bar\phi \frac{\vn}{(- \Delta )} \partial_i  \partial_j (\bar\varphi u_i b_j )= \left[\bar\phi, \, \frac{\vn \partial_i  \partial_j}{(- \Delta )}  \right] (\bar\varphi u_i b_j)  + \frac{\vn \partial_i  \partial_j}{(- \Delta )} ( \bar\phi u_i b_j ),
\end{equation}
where in the last term above we used the identity $\bar\phi=\bar\phi\bar\varphi$. Now, plugging the identity (\ref{Identite_Pression2}) in (\ref{Identite_Pression1}) and modifying 	accordingly expression (\ref{Definition_LocalPression}), we obtain the following formula for the term $\NN$ defined in (\ref{Definition_TermeNN}):
\begin{equation*}
\begin{split}
\NN &= \sum_{k=1}^{11} \NN_k=(\partial_t \bar\phi - \Delta \bar\phi)  (\vu + \vb)-2\sum_{i=1}^{3}(\partial_{i}\bar\phi) (\partial_i(\vu+\vb))-\bar\phi \left( (\vb\cdot\vn)\vu+(\vu\cdot\vn)\vb\right)\\
&-2\sum^3_{i,j= 1} \left[\bar\phi, \, \frac{ \vn \partial_j \partial_k}{(-\Delta)}\right] (\bar\varphi u_i b_j)-2\sum^3_{i,j= 1}  \frac{\vn  \pai \pj }{(- \Delta )}( \bar\phi u_i b_j )+2 \sum^3_{i,j= 1}  \frac{\bar\phi \vn\pai}{(- \Delta )} (\pj \bar\varphi) u_i b_j\\
&+2\sum^3_{i,j= 1}  \frac{\bar\phi \vn\partial_j}{(- \Delta )} (\partial_i \bar\varphi) u_i b_j- 2 \sum^3_{i,j= 1} \bar\phi \frac{\vn}{(-\Delta )}(\pai \pj\bar \varphi) (u_i b_j)\\
&-2 \bar\phi \frac{\vn \big(  (\Delta \bar\varphi) P \big) }{(-\Delta)}
+4 \sum^3_{i= 1}\bar\phi\frac{\vn\big(\partial_i ( (\partial_i \bar\varphi) P )\big)}{(-\Delta)}+\bar\phi(\vf+\vg).
\end{split}
\end{equation*}
Once we have obtained this expression for the term $\NN$, we study the solutions of the equation (\ref{Equation_ChaleurLocale1}) and we obtain
\begin{equation*}
\vVV = \displaystyle\int_{0}^{t} e^{ (t-s)\Delta}  \NN (s,\cdot)\,ds = \sum_{k=1}^{11}  \int_{0}^{t} e^{ (t-s)\Delta}  \NN_{k} (s,\cdot)\,ds :=\sum_{k=1}^{11} \vVV_k,
\end{equation*}
where we have\\
\begin{equation}\label{Def_VariablesVV1}
\begin{split}
\sum_{k=1}^{11} \vVV_k&=\underbrace{\int_{0}^{t} e^{ (t-s)\Delta} (\partial_t \bar\phi - \Delta \bar\phi)  (\vu + \vb)ds}_{\vVV_1}-2\underbrace{\sum_{i=1}^{3}\int_{0}^{t} e^{ (t-s)\Delta} (\partial_{i}\bar\phi) (\partial_i(\vu+\vb))ds}_{\vVV_2}\\
&-\underbrace{\int_{0}^{t} e^{ (t-s)\Delta} \bar\phi \left( (\vb\cdot\vn)\vu+(\vu\cdot\vn)\vb\right)ds}_{\vVV_3}-2\underbrace{\sum^3_{i,j= 1} \int_{0}^{t} e^{ (t-s)\Delta} \left[\bar\phi, \, \frac{ \vn \partial_i \partial_j}{(-\Delta)}\right] (\bar\varphi u_i b_j)ds}_{\vVV_4}\\
&-2\underbrace{\sum^3_{i,j= 1} \int_{0}^{t} e^{ (t-s)\Delta} \frac{\vn  \pai \pj }{(- \Delta )}( \bar\phi u_i b_j )ds}_{\vVV_5}+2\underbrace{\sum^3_{i,j= 1}  \int_{0}^{t} e^{ (t-s)\Delta} \frac{\bar\phi \vn\pai}{(- \Delta )} (\pj \bar\varphi) u_i b_jds}_{\vVV_6}\\
&+2\underbrace{\sum^3_{i,j= 1}  \int_{0}^{t} e^{ (t-s)\Delta} \frac{\bar\phi \vn\partial_j}{(- \Delta )} (\partial_i \bar\varphi) u_i b_jds}_{\vVV_7}- 2 \underbrace{\sum^3_{i,j= 1} \int_{0}^{t} e^{ (t-s)\Delta} \bar\phi \frac{\vn}{(-\Delta )}(\pai \pj\bar \varphi) (u_i b_j)ds}_{\vVV_8}\\
&-2 \underbrace{\int_{0}^{t} e^{ (t-s)\Delta} \bar\phi \frac{\vn }{(-\Delta)}\big(  (\Delta \varphi) P \big) ds}_{\vVV_9}
+4 \underbrace{\sum^3_{i= 1}\int_{0}^{t} e^{ (t-s)\Delta} \bar\phi \frac{\vn\partial_i  }{(-\Delta)}((\partial_i \bar\varphi) P )ds}_{\vVV_{10}}\\
&+\underbrace{\int_{0}^{t} e^{ (t-s)\Delta} \bar\phi(\vf+\vg)ds}_{\vVV_{11}}.
\end{split}
\end{equation}
We will study each one of these terms with the following lemma. 
\begin{Lemme}\label{Lemme_EstimationMorrey3Sigma}
In addition to the general hypotheses of Theorem \ref{Teorem1} let us further assume that we have all the informations stated in \eqref{Estimations135}. Then for all $k=1,…,11$ we have 
$$\mathds{1}_{Q_{R_2}}\vVV_k\in \mathcal{M}_{t,x}^{3, \sigma},$$
where $\tau_0<\sigma<10$.
\end{Lemme}
\begin{Remarque}
The upper bound $\sigma<10$ is given here to fix the possible values for this parameter. However we will see later on that in fact $\sigma$ must be quite close to $\tau_0$. See Remark \ref{Remarque_ConditionSigma} below for more details.
\end{Remarque}
{\bf Proof.}
\begin{itemize}
\item For the term $\vVV_1$, recalling that $e^{ (t-s)\Delta}f=\mathfrak{g}_{t-s}\ast f$ where $\mathfrak{g}_t$ is the usual $3D$-heat kernel, we can write
$$|\mathds{1}_{Q_{R_2}}\vVV_1(t,x)|=\left|\mathds{1}_{Q_{R_2}}\int_{0}^{t} \int_{\mathbb{R}^3}  \mathfrak{g}_{t-s}(x-y)[(\partial_t \bar\phi - \Delta \bar\phi)  (\vu + \vb)](s,y)dyds\right|,$$
and using the decay properties of the heat kernel as well as the properties of the test function $\bar\phi$ (see (\ref{DefSoporteFuncTest1})), we have
$$|\mathds{1}_{Q_{R_2}}\vVV_1(t,x)|\leq C\mathds{1}_{Q_{R_2}}\int_{\mathbb{R}} \int_{\mathbb{R}^3} \frac{1}{(|t-s|^{\frac{1}{2}}+|x-y|)^3}
\left| \mathds{1}_{Q_{\bar R}}  (\vu + \vb) (s,y)
\right| \,dy \,ds.$$
Now, recalling the definition of the Riesz potential given in (\ref{Def_Riesz_potential}) and since $Q_{R_2}\subset Q_{\bar R}$ we obtain the pointwise estimate
$$|\mathds{1}_{Q_{R_2}}\vVV_1(t,x)|\leq C\mathds{1}_{Q_{\bar R}}\mathcal{I}_{2}(  |\mathds{1}_{Q_{\bar R}}  (\vu + \vb)|)(t,x),$$
thus, taking Morrey $\mathcal{M}_{t,x}^{3, \sigma}$ norm in this inequality, we have 
$$\|\mathds{1}_{Q_{R_2}}\vVV_1(t,x)\|_{\mathcal{M}_{t,x}^{3, \sigma}}\leq C\|\mathds{1}_{Q_{\bar R}}\mathcal{I}_{2}(  |\mathds{1}_{Q_{\bar R}}  (\vu + \vb)|)\|_{\mathcal{M}_{t,x}^{3, \sigma}}.$$
Now, for some $2<q<\frac{5}{2}$ we set $\lambda=1-\frac{2q}{5}$ and we define $3=\frac{a}{\lambda}$ and $\sigma< 10 <\frac{q}{\lambda}$ (remark that $a\leq q$). Thus, by Lemma \ref{Lemme_locindi} and by Lemma \ref{Lemme_Hed} we can write:
\begin{eqnarray}
\|\mathds{1}_{Q_{\bar R}}\mathcal{I}_{2}(  |\mathds{1}_{Q_{\bar R}}  (\vu + \vb)|)\|_{\mathcal{M}_{t,x}^{3, \sigma}}&\leq &C\|\mathcal{I}_{2}(  |\mathds{1}_{Q_{\bar R}}  (\vu + \vb)|)\|_{\mathcal{M}_{t,x}^{\frac{a}{\lambda},\frac{q}{\lambda}}},\notag\\
&\leq &C\|\mathds{1}_{Q_{\bar R}}  (\vu + \vb)\|_{\mathcal{M}_{t,x}^{a, q}}\leq C\|\mathds{1}_{Q_{R_3}}  (\vu + \vb)\|_{\mathcal{M}_{t,x}^{3, \tau_0}}<+\infty,\label{EstimationRieszTermVV1}
\end{eqnarray}
where in the last estimate we applied Lemma \ref{Lemme_locindi} again noting that $a\leq 3$ and $q<\tau_0$.
\item For the second term of (\ref{Def_VariablesVV1}) we start writing $(\partial_{i}\bar\phi) (\partial_i(\vu+\vb))=\partial_i((\partial_i\bar\phi)(\vu+\vb))-(\partial_i^2\bar\phi)(\vu+\vb)$, and we have
\begin{equation}\label{EstimationPonctuelleVV2}
|\mathds{1}_{Q_{R_2}}\vVV_2(t,x)|\leq \sum_{i=1}^{3}\left|\mathds{1}_{Q_{R_2}}\int_{0}^{t} e^{ (t-s)\Delta} \partial_i\big((\partial_{i}\bar\phi) (\vu+\vb)\big)ds\right|+\left|\mathds{1}_{Q_{R_2}}\int_{0}^{t} e^{ (t-s)\Delta} (\partial_{i}^2\bar\phi)(\vu+\vb)ds\right|.
\end{equation}
For the first term above, by the properties of the heat kernel and by the definition of the Riesz potential $\mathcal{I}_{1}$ (see (\ref{Def_Riesz_potential})), we obtain
\begin{eqnarray*}
\left|\mathds{1}_{Q_{R_2}}\int_{0}^{t} e^{ (t-s)\Delta} \partial_i\big((\partial_{i}\bar\phi) (\vu+\vb)\big)ds\right|&=&\left|\mathds{1}_{Q_{R_2}}\int_{0}^{t} \int_{\mathbb{R}^3}\partial_i\mathfrak{g}_{t-s}(x-y)(\partial_{i}\bar\phi) (\vu+\vb)(s,y)dyds\right|\\
&\leq &C\mathds{1}_{Q_{R_2}}\int_{\mathbb{R}}\int_{\mathbb{R}^3} \frac{|\mathds{1}_{Q_{\bar R}} (\vu+\vb)(s,y)|}{(|t-s|^{\frac{1}{2}}+|x-y|)^4}dyds\\
&\leq & C\mathds{1}_{Q_{R_2}}(\mathcal{I}_1(|\mathds{1}_{Q_{\bar R}} (\vu+\vb)|))(t,x).
\end{eqnarray*}
The second term of (\ref{EstimationPonctuelleVV2}) can be treated as the term $\vVV_1$ and we have the pointwise estimate
$$\left|\mathds{1}_{Q_{R_2}}\int_{0}^{t} e^{ (t-s)\Delta} (\partial_{i}^2\bar\phi)(\vu+\vb)ds\right|\leq C\mathds{1}_{Q_{\bar R}}\mathcal{I}_{2}(  |\mathds{1}_{Q_{\bar R}}  (\vu + \vb)|)(t,x),$$
and gathering these two estimates we have 
$$|\mathds{1}_{Q_{R_2}}\vVV_2(t,x)|\leq C\mathds{1}_{Q_{R_2}}(\mathcal{I}_1(|\mathds{1}_{Q_{\bar R}} (\vu+\vb)|))(t,x)+C\mathds{1}_{Q_{\bar R}}\mathcal{I}_{2}(  |\mathds{1}_{Q_{\bar R}}  (\vu + \vb)|)(t,x),$$
and taking the Morrey $\mathcal{M}_{t,x}^{3, \sigma}$ we obtain
$$\|\mathds{1}_{Q_{R_2}}\vVV_2\|_{\mathcal{M}_{t,x}^{3, \sigma}}\leq C\|\mathds{1}_{Q_{R_2}}(\mathcal{I}_1(|\mathds{1}_{Q_{\bar R}} (\vu+\vb)|))\|_{\mathcal{M}_{t,x}^{3, \sigma}}+C\|\mathds{1}_{Q_{\bar R}}\mathcal{I}_{2}(|\mathds{1}_{Q_{\bar R}}  (\vu + \vb)|)\|_{\mathcal{M}_{t,x}^{3, \sigma}}.$$
The second term of the right-hand above can be treated in the same manner as (\ref{EstimationRieszTermVV1}), thus we only study now the quantity $\|\mathds{1}_{Q_{R_2}}(\mathcal{I}_1(|\mathds{1}_{Q_{\bar R}} (\vu+\vb)|))\|_{\mathcal{M}_{t,x}^{3, \sigma}}$. For some $4\leq q <5$ we define $\lambda=1-\frac{q}{5}$, noting that $3\leq \frac{3}{\lambda}$ and $\sigma <10<\frac{q}{\lambda}$, by Lemma \ref{Lemme_Hed}, we can write
\begin{eqnarray*}
\|\mathds{1}_{Q_{R_2}}(\mathcal{I}_1(|\mathds{1}_{Q_{\bar R}} (\vu+\vb)|))\|_{\mathcal{M}_{t,x}^{3, \sigma}}&\leq &C\|\mathcal{I}_1(|\mathds{1}_{Q_{\bar R}} (\vu+\vb)|)\|_{\mathcal{M}_{t,x}^{\frac{3}{\lambda}, \frac{q}{\lambda}}}\leq C\|\mathds{1}_{Q_{\bar R}} (\vu+\vb)\|_{\mathcal{M}_{t,x}^{3, q}}\\
&\leq & C\|\mathds{1}_{Q_{R_3}} (\vu+\vb)\|_{\mathcal{M}_{t,x}^{3, \tau_0}}<+\infty,
\end{eqnarray*}
from which we deduce that $\|\mathds{1}_{Q_{R_2}}\vVV_2\|_{\mathcal{M}_{t,x}^{3, \sigma}}<+\infty$.
\item For the term $\vVV_3$ in (\ref{Def_VariablesVV1}), in a similar manner we obtain the inequality
\begin{eqnarray*}
|\mathds{1}_{Q_{R_2}}\vVV_3(t,x)|&=&\left|\mathds{1}_{Q_{R_2}}\int_{0}^{t} \int_{\mathbb{R}^3}  \mathfrak{g}_{t-s}(x-y)\left[\bar\phi \left( (\vb\cdot\vn)\vu+(\vu\cdot\vn)\vb\right)\right](s,y)dyds\right|\\
&\leq & C\mathds{1}_{Q_{R_2}}\int_{\mathbb{R}} \int_{\mathbb{R}^3} \frac{\left|\bar\phi \left( (\vb\cdot\vn)\vu+(\vu\cdot\vn)\vb\right)\right|(s,y)}{(|t-s|^{\frac{1}{2}}+|x-y|)^3} \,dy \,ds\\
&\leq & C\mathds{1}_{Q_{R_2}}\mathcal{I}_2\left(\left|\mathds{1}_{Q_{\bar R}}\left( (\vb\cdot\vn)\vu+(\vu\cdot\vn)\vb\right)\right|\right)(t,x),
\end{eqnarray*}
from which we deduce 
\begin{equation}\label{Decomposition2termesVV3}
\|\mathds{1}_{Q_{R_2}}\vVV_3\|_{\mathcal{M}_{t,x}^{3, \sigma}}\leq C\left\|\mathds{1}_{Q_{R_2}}\mathcal{I}_2\left(|\mathds{1}_{Q_{\bar R}} (\vb\cdot\vn)\vu|\right)\right\|_{\mathcal{M}_{t,x}^{3, \sigma}}+C\left\|\mathds{1}_{Q_{R_2}}\mathcal{I}_2\left(|\mathds{1}_{Q_{\bar R}} (\vu\cdot\vn)\vb|\right)\right\|_{\mathcal{M}_{t,x}^{3, \sigma}}.
\end{equation}
As we have completely symmetric information on $\vu$ and $\vb$ it is enough the study one of these terms and we will treat the first one. We set now $\frac{5}{3-\alpha}<q<\frac{5}{2}$ and $\lambda=1-\frac{2q}{5}$. Since $3\leq \frac{6}{5\lambda}$ and $\tau_0<\sigma<\frac{q}{\lambda}$, applying Lemma \ref{Lemme_locindi} and Lemma \ref{Lemme_Hed} we have
$$\left\|\mathds{1}_{Q_{R_2}}\mathcal{I}_2\left(|\mathds{1}_{Q_{\bar R}} (\vb\cdot\vn)\vu|\right)\right\|_{\mathcal{M}_{t,x}^{3, \sigma}}\leq C\left\|\mathds{1}_{Q_{R_2}}\mathcal{I}_2\left(|\mathds{1}_{Q_{\bar R}} (\vb\cdot\vn)\vu|\right)\right\|_{\mathcal{M}_{t,x}^{\frac{6}{5\lambda},\frac{q}{\lambda}}}\leq C\left\|\mathds{1}_{Q_{\bar R}} (\vb\cdot\vn)\vu\right\|_{\mathcal{M}_{t,x}^{\frac{6}{5}, q}}.$$
Recall that we have $\frac{5}{1-\alpha}<\tau_0<\sigma< 10$ and by the H\"older inequality in Morrey spaces (see Lemma \ref{Lemme_Product}) we obtain
$$\left\|\mathds{1}_{Q_{\bar R}} (\vb\cdot\vn)\vu\right\|_{\mathcal{M}_{t,x}^{\frac{6}{5}, q}}\leq \left\|\mathds{1}_{Q_{R_3}}\vb\right\|_{\mathcal{M}_{t,x}^{3, \tau_0}}\left\|\mathds{1}_{Q_{R_3}}\vn \otimes \vu\right\|_{\mathcal{M}_{t,x}^{2, \tau_1}}<+\infty,$$
where $\frac{1}{q}=\frac{1}{\tau_0}+\frac{1}{\tau_1}=\frac{2}{\tau_0}+\frac{1}{5}$. Note that the condition $\frac{5}{1-\alpha}<\tau_0<\sigma< 10$ and the relationship $\frac{1}{q}=\frac{2}{\tau_0}+\frac{1}{5}$ are compatible with the fact that $\frac{5}{3-\alpha}<q < \frac{5}{2}$. Applying exactly the same ideas in the second term of (\ref{Decomposition2termesVV3}) we finally obtain 
$$\|\mathds{1}_{Q_{R_2}}\vVV_3\|_{\mathcal{M}_{t,x}^{3,\sigma}}<+\infty.$$
\begin{Remarque}\label{Remarque_ConditionSigma}
The condition $\frac{5}{3-\alpha}<q<\frac{5}{2}$ required to apply the H\"older inequality jointly with the constraint $\tau_0<\sigma<\frac{q}{\lambda}$ implies that $\sigma$ must be very close from $\tau_0$.
\end{Remarque}
\item The quantity $\vVV_4$ in (\ref{Def_VariablesVV1}) is the most technical one and it will be treated as follows
\begin{eqnarray*}
|\mathds{1}_{Q_{R_2}}\vVV_4(t,x)|&\leq &\sum^3_{i,j= 1} \mathds{1}_{Q_{R_2}}\int_{\mathbb{R}} \int_{\mathbb{R}^3}\frac{\left|\left[\bar\phi, \, \frac{ \vn \partial_i\partial_j}{(-\Delta)}\right] (\bar\varphi u_i b_j)(s,y)\right|}{(|t-s|^{\frac{1}{2}}+|x-y|)^3}dyds\\
&\leq &\sum^3_{i,j= 1}\mathds{1}_{Q_{R_2}}\mathcal{I}_{2}\left(\left|\left[\bar\phi, \, \frac{ \vn \partial_i\partial_j}{(-\Delta)}\right] (\bar\varphi u_i b_j)\right|\right)(t,x),
\end{eqnarray*}
and taking the Morrey $\mathcal{M}_{t,x}^{3, \sigma}$ norm we have
$$\|\mathds{1}_{Q_{R_2}}\vVV_4\|_{\mathcal{M}_{t,x}^{3, \sigma}}\leq \sum^3_{i,j= 1}\left\|\mathds{1}_{Q_{R_2}}\mathcal{I}_{2}\left(\left|\left[\bar\phi, \, \frac{ \vn \partial_i\partial_j}{(-\Delta)}\right] (\bar\varphi u_i b_j)\right|\right)\right\|_{\mathcal{M}_{t,x}^{3, \sigma}}.$$
If we set $\frac{1}{q}=  \frac{2}{\tau_0} + \frac{1}{5}$ and $\lambda=1-\frac{2q}{5}$ then we have $3\leq \frac{3}{2\lambda}$ and $\sigma \leq  \frac{q}{\lambda} = \frac{5 \tau_0}{10- \tau_0}$ and by Lemma \ref{Lemme_locindi} and Lemma \ref{Lemme_Hed} we obtain:
\begin{eqnarray*}
\left\|\mathds{1}_{Q_{R_2}}\mathcal{I}_{2}\left(\left|\left[\bar\phi, \, \frac{ \vn \partial_i\partial_j}{(-\Delta)}\right] (\bar\varphi u_i b_j)\right|\right)\right\|_{\mathcal{M}_{t,x}^{3, \sigma}}&\leq &C\left\|\mathds{1}_{Q_{R_2}}\mathcal{I}_{2}\left(\left|\left[\bar\phi, \, \frac{ \vn \partial_i\partial_j}{(-\Delta)}\right] (\bar\varphi u_i b_j)\right|\right)\right\|_{\mathcal{M}_{t,x}^{\frac{3}{2\lambda}, \frac{q}{\lambda}}}\\
&\leq &C\left\|\left[\bar\phi, \, \frac{ \vn \partial_i\partial_j}{(-\Delta)}\right] (\bar\varphi u_i b_j)\right\|_{\mathcal{M}_{t,x}^{\frac{3}{2}, q}},
\end{eqnarray*}
We will study this norm and by the definition of Morrey spaces (\ref{Morreyparabolic}), if we introduce a threshold $\mathfrak{r}=\frac{\bar R-R_2}{2}$, we have
\end{itemize}
\begin{equation}\label{CommutatorEstimatevVV4}
\begin{split}
\left\|\left[\bar\phi, \, \frac{ \vn \partial_i\partial_j}{(-\Delta)}\right] (\bar\varphi u_i b_j)\right\|_{\mathcal{M}_{t,x}^{\frac{3}{2}, q}}^{\frac{3}{2}}&\leq\underset{\underset{0<r<\mathfrak{r}}{(\mathfrak{t},\bar{x})}}{\sup}\;\frac{1}{r^{5(1-\frac{3}{2q})}}\int_{Q_r(\mathfrak{t},\bar x)}\left|\left[\bar\phi, \, \frac{ \vn \partial_i\partial_j}{(-\Delta)}\right] (\bar\varphi u_i b_j)\right|^{\frac{3}{2}}dxdt\qquad\\
&+\underset{\underset{\mathfrak{r}< r}{(\mathfrak{t},\bar{x})}}{\sup}\;\frac{1}{r^{5(1-\frac{3}{2q})}}\int_{Q_r(\mathfrak{t},\bar{x})}\left|\left[\bar\phi, \, \frac{ \vn \partial_i\partial_j}{(-\Delta)}\right] (\bar\varphi u_i b_j)\right|^{\frac{3}{2}}dxdt.\qquad
\end{split}
\end{equation}
\begin{itemize}
\item[]
Now, we study the second term of the right-hand side above, which is easy to handle as we have $\mathfrak{r}<r$ and we can write
$$\underset{\underset{\mathfrak{r}< r}{(\mathfrak{t},\bar{x}) \in \mathbb{R}\times \R}}{\sup}\;\frac{1}{r^{5(1-\frac{3}{2q})}}\int_{Q_r(\mathfrak{t},\bar{x})}\left|\left[\bar\phi, \, \frac{ \vn \partial_i\partial_j}{(-\Delta)}\right] (\bar\varphi u_i b_j)\right|^{\frac{3}{2}}dxdt\leq C_{\mathfrak{r}}\left\|\left[\bar\phi, \, \frac{ \vn \partial_i\partial_j}{(-\Delta)}\right] (\bar\varphi u_i b_j)\right\|_{L^{\frac{3}{2}}_{t,x}}^{\frac{3}{2}},$$
and since $\bar\phi$ is a regular function and $\frac{ \vn \partial_i\partial_j}{(-\Delta)}$ is a Calder\'on-Zydmund operator, by the Calder\'on commutator theorem (see \cite{PGLR0}), we have that the operator $\left[\bar\phi, \, \frac{ \vn \partial_i\partial_j}{(-\Delta)}\right]$ is bounded in the space $L^{\frac{3}{2}}_{t,x}$ and we can write
\begin{eqnarray*}
\left\|\left[\bar\phi, \, \frac{ \vn \partial_i\partial_j}{(-\Delta)}\right] (\bar\varphi u_i b_j)\right\|_{L^{\frac{3}{2}}_{t,x}}&\leq &C\left\|\bar\varphi u_i b_j\right\|_{L^{\frac{3}{2}}_{t,x}}\leq C\|\mathds{1}_{Q_{R_3}} u_i b_j\|_{\mathcal{M}^{\frac{3}{2}, \frac{3}{2}}_{t,x}}\\
&\leq & C\|\mathds{1}_{Q_{R_3}} \vu\|_{\mathcal{M}^{3,3}_{t,x}}\|\mathds{1}_{Q_{R_3}} \vb\|_{\mathcal{M}^{3, 3}_{t,x}}\leq C\|\mathds{1}_{Q_{R_3}} \vu\|_{\mathcal{M}^{3,\tau_0}_{t,x}}\|\mathds{1}_{Q_{R_3}} \vb\|_{\mathcal{M}^{3, \tau_0}_{t,x}}<+\infty,
\end{eqnarray*}
where in the last line we used H\"older inequalities in Morrey spaces and we applied Lemma \ref{Lemme_locindi}.\\

The first term of the right-hand side of  (\ref{CommutatorEstimatevVV4}) requires some extra computations: indeed, as we are interested to obtain information over the parabolic ball $Q_{r}(\mathfrak{t}, \bar{x})$ we can write
for some  $0<r<\mathfrak{r}$:
\begin{equation}\label{CommutatorEstimatevVV401}
\mathds{1}_{Q_{r}}\left[\bar\phi, \, \frac{ \vn \partial_i\partial_j}{(-\Delta)}\right] (\bar\varphi u_i b_j))=\mathds{1}_{Q_{r}}\left[\bar\phi, \, \frac{ \vn \partial_i\partial_j}{(-\Delta)}\right] (\mathds{1}_{Q_{2r}}\bar\varphi u_i b_j)+\mathds{1}_{Q_{r}}\left[\bar\phi, \, \frac{ \vn \partial_i\partial_j}{(-\Delta)}\right] ((\mathbb{I}-\mathds{1}_{Q_{2r}})\bar\varphi u_i b_j),
\end{equation}
and as before we will study the $L^{\frac{3}{2}}_{t,x}$ norm of these two terms. For the first quantity in the right-hand side of (\ref{CommutatorEstimatevVV401}), by the Calder\'on commutator theorem, by the definition of Morrey spaces and by the H\"older inequalities we have 
\begin{eqnarray*}
\left\|\mathds{1}_{Q_{r}}\left[\bar\phi, \, \frac{ \vn \partial_i\partial_j}{(-\Delta)}\right] (\mathds{1}_{Q_{2r}}\bar\varphi u_i b_j)\right\|_{L^{\frac{3}{2}}_{t,x}}^{\frac{3}{2}}&\leq& C\|\mathds{1}_{Q_{2r}}\bar\varphi u_i b_j\|_{L^{\frac{3}{2}}_{t,x}}^{\frac{3}{2}}\leq Cr^{5 (1-\frac{3}{\tau_0})} \|\mathds{1}_{Q_{R_3}}u_i b_j\|_{\mathcal{M}^{\frac{3}{2}, \frac{\tau_0}{2}}_{t,x}}^{\frac{3}{2}}\\
&\leq & Cr^{5 (1-\frac{3}{\tau_0})} \|\mathds{1}_{Q_{R_3}}\vu\|_{\mathcal{M}^{3, \tau_0}_{t,x}}^{\frac{3}{2}}\|\mathds{1}_{Q_{R_3}}\vb\|_{\mathcal{M}^{3, \tau_0}_{t,x}}^{\frac{3}{2}},
\end{eqnarray*}
for all $0<r<\mathfrak{r}$, from which we deduce that 
$$\underset{\underset{0<r<\mathfrak{r}}{(\mathfrak{t},\bar{x}) }}{\sup}\;\frac{1}{r^{5(1-\frac{3}{2q})}}\int_{Q_r(\mathfrak{t},\bar{x})}\left|\mathds{1}_{Q_{r}}\left[\bar\phi, \, \frac{ \vn \partial_i\partial_j}{(-\Delta)}\right] (\mathds{1}_{Q_{2r}}\bar\varphi u_i b_j)\right|^{\frac{3}{2}}dxdt\leq C \|\mathds{1}_{Q_{R_3}}\vu\|_{\mathcal{M}^{3, \tau_0}_{t,x}}^{\frac{3}{2}}\|\mathds{1}_{Q_{R_3}}\vb\|_{\mathcal{M}^{3, \tau_0}_{t,x}}^{\frac{3}{2}}<+\infty.$$
We study now the second term of the right-hand side of (\ref{CommutatorEstimatevVV401}) and for this we consider the following operator:
$$T: f \mapsto  \left(\mathds{1}_{Q_{r}}  \left[\bar\phi, \, \frac{ \vn \pai \pj}{- \Delta }\right] (\mathbb{I} - \mathds{1}_{Q_{2r}}) \bar \varphi
\right) f,$$
and by the properties of the convolution kernel of the operator $\frac{1}{(-\Delta)}$ we obtain
$$|T(f)(x)|\leq C\mathds{1}_{Q_{r}}(x)\int_{\mathbb{R}^3}\frac{(\mathbb{I} - \mathds{1}_{Q_{2r}})(y) \mathds{1}_{Q_{R_3}}(y) |f(y)| |\bar \phi(x)-\bar \phi(y)|}{|x-y|^4} dy.$$
Recalling that $0<r<\mathfrak{r}=\frac{\bar R-R_2}{2}$, by the support properties of the test function $\bar\phi$ (see (\ref{DefSoporteFuncTest1})), the integral above is meaningful if $|x-y|>r$ and thus we can write
\begin{eqnarray*}
\left\|\mathds{1}_{Q_{r}}\left[\bar\phi, \, \frac{ \vn \partial_i\partial_j}{(-\Delta)}\right] ((\mathbb{I}-\mathds{1}_{Q_{2r}})\bar\varphi u_i b_j)\right\|_{L^\frac{3}{2}_{t,x}}^\frac{3}{2}&\leq&C\left\|\mathds{1}_{Q_{r}} \int_{\mathbb{R}^3} \frac{\mathds{1}_{|x-y| > r}}{|x-y|^4}(\mathbb{I} - \mathds{1}_{Q_{2r}})(y) \mathds{1}_{Q_{R_3}}(y)|u_i b_j |dy\right\|_{L^\frac{3}{2}_{t,x}}^\frac{3}{2}\\
&\leq &C\left(\int_{|y|>r}\frac{1}{|y|^4}\| \mathds{1}_{Q_{R_3}}
|u_i b_j |(\cdot-y)\|_{L^\frac{3}{2}_{t,x}(Q_r)}dy\right)^\frac{3}{2}\\
&\leq & Cr^{-\frac{3}{2}}\| \mathds{1}_{Q_{R_3}}u_i b_j \|_{L^\frac{3}{2}_{t,x}(Q_{r})}^\frac{3}{2},
\end{eqnarray*}
with this estimate at hand and using the definition of Morrey spaces, we can write
\begin{eqnarray*}
\int_{Q_r(\mathfrak{t},\bar{x})}\left|\mathds{1}_{Q_{r}}\left[\bar\phi, \, \frac{ \vn \partial_i\partial_j}{(-\Delta)}\right] ((\mathbb{I}-\mathds{1}_{Q_{2r}})\bar\varphi u_i b_j)\right|^{\frac{3}{2}}dxdt&\leq &Cr^{-\frac{3}{2}}r^{5 (1-\frac{3}{\tau_0})}\| \mathds{1}_{Q_{R_3}}u_i b_j \|_{\mathcal{M}^{\frac{3}{2}, \frac{\tau_0}{2}}_{t,x}}^\frac{3}{2}\\
&\leq &Cr^{5(1-\frac{3}{2q})}\| \mathds{1}_{Q_{R_3}}u_i b_j \|_{\mathcal{M}^{\frac{3}{2}, \frac{\tau_0}{2}}_{t,x}}^\frac{3}{2},
\end{eqnarray*}
where in the last inequality we used the fact that $\frac{1}{q}=  \frac{2}{\tau_0} + \frac{1}{5}$, which implies  $r^{-\frac{3}{2}}r^{5 (1-\frac{3}{\tau_0})}= r^{5(1-\frac{3}{2q})}$. Thus we finally obtain
$$\underset{\underset{0<r<\mathfrak{r}}{(\mathfrak{t},\bar{x}) }}{\sup}\;\frac{1}{r^{5(1-\frac{3}{2q})}}\int_{Q_r(\mathfrak{t},\bar{x})}\left|\mathds{1}_{Q_{r}}\left[\bar\phi, \, \frac{ \vn \partial_i\partial_j}{(-\Delta)}\right] ((\mathbb{I}-\mathds{1}_{Q_{2r}})\bar\varphi u_i b_j)\right|^{\frac{3}{2}}dxdt\leq C \|\mathds{1}_{Q_{R_3}}\vu\|_{\mathcal{M}^{3, \tau_0}_{t,x}}^{\frac{3}{2}}\|\mathds{1}_{Q_{R_3}}\vb\|_{\mathcal{M}^{3, \tau_0}_{t,x}}^{\frac{3}{2}}<+\infty.$$
We have proven so far that all the term in (\ref{CommutatorEstimatevVV4}) are bounded and we can conclude that 
$$\|\mathds{1}_{Q_{R_2}}\vVV_4\|_{\mathcal{M}_{t,x}^{3, \sigma}}<+\infty.$$

\item For the quantity $\vVV_5$ in (\ref{Def_VariablesVV1}) we write
\begin{eqnarray*}
|\mathds{1}_{Q_{R_2}}\vVV_5(t,x)|&\leq&\sum^3_{i,j= 1} \mathds{1}_{Q_{R_2}}\left|\int_{0}^{t}\int_{\mathbb{R}^3}\vn\mathfrak{g}_{t-s}(x-y)\frac{\pai}{\sqrt{- \Delta}}\frac{ \pj }{\sqrt{- \Delta}}( \bar\phi u_i b_j )(s,y)dyds\right|\\
&\leq &C\sum^3_{i,j= 1} \mathds{1}_{Q_{R_2}}\int_{\mathbb{R}}\int_{\mathbb{R}^3} \frac{|\mathcal{R}_i\mathcal{R}_j( \bar\phi u_i b_j )(s,y)|}{(|t-s|^{\frac{1}{2}}+|x-y|)^4}dyds\\
&\leq &C\sum^3_{i,j= 1} \mathds{1}_{Q_{R_2}} \mathcal{I}_{1}\left(|\mathcal{R}_i\mathcal{R}_j( \bar\phi u_i b_j )|\right)(t,x),
\end{eqnarray*}
where we used the decaying properties of the heat kernel (recall that $\mathcal{R}_i=\frac{\pai}{\sqrt{- \Delta}}$ are the Riesz transforms). 
Now taking the Morrey $\mathcal{M}^{3, \sigma}_{t,x}$ norm and by Lemma \ref{Lemme_locindi} (with $\nu=\frac{4\tau_0+5}{5\tau_0}$, $p=3$, $q=\tau_0$ such that $\frac{p}{\nu}>3$ and $\frac{q}{\nu}>\sigma$ which is compatible with the condition $\tau_0<\sigma$) we have
\begin{eqnarray*}
\|\mathds{1}_{Q_{R_2}}\vVV_5\|_{\mathcal{M}^{3, \sigma}_{t,x}}&\leq &C \sum^3_{i,j= 1}\|  \mathds{1}_{Q_{R_2}}\mathcal{I}_{1}\left(|\mathcal{R}_i\mathcal{R}_j( \bar\phi u_i b_j )|\right)\|_{\mathcal{M}^{\frac{p}{\nu}, \frac{q}{\nu}}_{t,x}}
\end{eqnarray*}
Then by Lemma \ref{Lemme_Hed} with $\lambda= 1 - \tfrac{\tau_0 /2}{5}$ (recall $\frac{5}{1-\alpha}<\tau_0 <10$ so that $\nu > 2 \lambda$) and by the boundedness of Riesz transforms in Morrey spaces we obtain:
\begin{eqnarray*}
\|\mathds{1}_{Q_{R_2}} \mathcal{I}_{1}\left(|\mathcal{R}_i\mathcal{R}_j( \bar\phi u_i b_j )|\right)\|_{\mathcal{M}^{\frac{p}{\nu}, \frac{q}{\nu}}_{t,x}}
&\leq&C\|\mathcal{I}_{1}\left(|\mathcal{R}_i\mathcal{R}_j( \bar\phi u_i b_j )|\right)\|_{\mathcal{M}^{\frac{p}{2\lambda}, \frac{q}{2\lambda}}_{t,x}}
\leq C \|\mathcal{R}_i\mathcal{R}_j( \bar\phi u_i b_j )\|_{\mathcal{M}^{\frac{3}{2},\frac{\tau_0}{2}}_{t,x}}\\
&\leq &\| \mathds{1}_{Q_{R_3}} u_i b_j \|_{\mathcal{M}^{\frac{3}{2},\frac{\tau_0}{2}}_{t,x}}
\leq  C\|\mathds{1}_{Q_{R_3}} \vu\|_{\mathcal{M}^{3,\tau_0}_{t,x}}\|  \mathds{1}_{Q_{R_3}} \vb\|_{\mathcal{M}^{3,\tau_0}_{t,x}}<+\infty.
\end{eqnarray*}

\item The quantities $\vVV_6$ and $\vVV_7$ in (\ref{Def_VariablesVV1}) can be treated in a very similar fashion since their inner structure is essentially the same. We thus only treat here the term $\vVV_6$ and following the same ideas we have
$$|\mathds{1}_{Q_{R_2}}\vVV_6|\leq C\sum^3_{i,j= 1}  \mathds{1}_{Q_{R_2}}\int_{\mathbb{R}}  \int_{\mathbb{R}^3}\frac{\left|\frac{\bar\phi \vn\pai}{(- \Delta )} (\pj \bar\varphi) u_i b_j(s,y)\right|}{(|t-s|^{\frac{1}{2}}+|x-y|)^3} dyds=C\sum^3_{i,j= 1}  \mathds{1}_{Q_{R_2}}\mathcal{I}_{2}\left(\left|\frac{\bar\phi \vn\pai}{(- \Delta )} (\pj \bar\varphi) u_i b_j\right|\right).$$
For $2<q<\frac{5}{2}$, define $\lambda=1-\frac{2q}{5}$, we thus have $3\leq \frac{3}{2\lambda}$ and $\sigma<10\leq \frac{q}{\lambda}$. Then, by Lemma \ref{Lemme_locindi} and Lemma \ref{Lemme_Hed} we can write
\begin{eqnarray*}
\left\|\mathds{1}_{Q_{R_2}}\mathcal{I}_{2}\left(\left|\frac{\bar\phi \vn\pai}{(- \Delta )} (\pj \bar\varphi) u_i b_j\right|\right)\right\|_{\mathcal{M}^{3, \sigma}_{t,x}}&\leq &C\left\|\mathds{1}_{Q_{R_2}}\mathcal{I}_{2}\left(\left|\frac{\bar\phi \vn\pai}{(- \Delta )} (\pj \bar\varphi) u_i b_j\right|\right)\right\|_{\mathcal{M}^{\frac{3}{2\lambda}, \frac{q}{\lambda}}_{t,x}}\\
&\leq & C\left\|\frac{\bar\phi \vn\pai}{(- \Delta )} (\pj \bar\varphi) u_i b_j\right\|_{\mathcal{M}^{\frac{3}{2}, q}_{t,x}},
\end{eqnarray*}
but since the operator $\frac{\bar\phi \vn\pai}{(- \Delta )}$ is bounded in Morrey spaces and since $2<q<\frac{5}{2}< \tfrac{\tau_0}{2}$, one has by Lemma \ref{Lemme_locindi} and by the H\"older inequalities
\begin{eqnarray*}
\left\|\frac{\bar\phi \vn\pai}{(- \Delta )} (\pj \bar\varphi) u_i b_j\right\|_{\mathcal{M}^{\frac{3}{2}, q}_{t,x}}&\leq &C\left\| \mathds{1}_{Q_{R_3}}u_i b_j\right\|_{\mathcal{M}^{\frac{3}{2}, q}_{t,x}}\leq \left\| \mathds{1}_{Q_{R_3}}u_i b_j\right\|_{\mathcal{M}^{\frac{3}{2}, \frac{\tau_0}{2}}_{t,x}}\\
&\leq &C\|\mathds{1}_{Q_{R_3}} \vu\|_{\mathcal{M}^{3,\tau_0}_{t,x}}\|  \mathds{1}_{Q_{R_3}} \vb\|_{\mathcal{M}^{3,\tau_0}_{t,x}}<+\infty,
\end{eqnarray*}
from which we deduce $\|\mathds{1}_{Q_{R_2}}\vVV_6\|_{\mathcal{M}^{3, \sigma}_{t,x}}<+\infty$. The same computations can be performed to obtain that $\|\mathds{1}_{Q_{R_2}}\vVV_7\|_{\mathcal{M}^{3, \sigma}_{t,x}}<+\infty$.

\item The quantity $\vVV_8$ in (\ref{Def_VariablesVV1}) is treated in the following manner: we first write
\begin{eqnarray*}
|\mathds{1}_{Q_{R_2}}\vVV_8(t,x)|&\leq &\sum^3_{i,j= 1}\left| \mathds{1}_{Q_{R_2}}\int_{0}^{t} e^{ (t-s)\Delta} \bar\phi \frac{\vn}{(-\Delta )}(\pai \pj\bar \varphi) (u_i b_j)ds\right|\\
&\leq &C\sum^3_{i,j= 1} \mathds{1}_{Q_{R_2}}\left(\mathcal{I}_2\left|\bar\phi \frac{\vn}{(-\Delta )}(\pai \pj\bar \varphi) (u_i b_j)\right|\right)(t,x),
\end{eqnarray*}
from which we deduce the estimate
$$\|\mathds{1}_{Q_{R_2}}\vVV_8\|_{\mathcal{M}^{3, \sigma}_{t,x}}
\leq C \sum^3_{i,j= 1}\left\| \mathds{1}_{Q_{R_2}}\left(\mathcal{I}_2\left|\bar\phi \frac{\vn}{(-\Delta )}(\pai \pj\bar \varphi) (u_i b_j)\right|\right)\right\|_{\mathcal{M}^{3,\sigma}_{t,x}}.$$
We set $1<\nu<\frac{3}{2}$, $2\nu <q<\frac{5\nu}{2}$ and $\lambda=1-\frac{2q}{5\nu}$, thus we have $3\leq \frac{\nu}{\lambda}$ and $\sigma <10< \frac{q}{\lambda}$, then, by Lemma \ref{Lemme_locindi} and by Lemma \ref{Lemme_Hed} we can write
\begin{eqnarray}
\left\|\mathds{1}_{Q_{R_2}}\left(\mathcal{I}_2\left|\bar\phi \frac{\vn}{(-\Delta )}(\pai \pj\bar \varphi) (u_i b_j)\right|\right)\right\|_{\mathcal{M}^{3,\sigma}_{t,x}}\leq C\left\| \mathds{1}_{Q_{R_2}}\left(\mathcal{I}_2\left|\bar\phi \frac{\vn}{(-\Delta )}(\pai \pj\bar \varphi) (u_i b_j)\right|\right)\right\|_{\mathcal{M}^{\frac{\nu}{\lambda},\frac{q}{\lambda}}_{t,x}}&&\notag\\
\leq C\left\| \bar\phi \frac{\vn}{(-\Delta )}(\pai \pj\bar \varphi) (u_i b_j)\right\|_{\mathcal{M}^{\nu,q}_{t,x}}\leq C\left\| \bar\phi \frac{\vn}{(-\Delta )}(\pai \pj\bar \varphi) (u_i b_j)\right\|_{\mathcal{M}^{\nu,\frac{5\nu}{2}}_{t,x}}&&\notag\\
\leq C\left\| \bar\phi \frac{\vn}{(-\Delta )}(\pai \pj\bar \varphi) (u_i b_j)\right\|_{L^{\nu}_tL^{\infty}_x},&&\label{Formula_intermediairevVV80}
\end{eqnarray}
where in the last estimate we used the space inclusion $L^{\nu}_tL^{\infty}_x\subset \mathcal{M}^{\nu,\frac{5\nu}{2}}_{t,x}$. Let us focus now in the $L^\infty$ norm above (\emph{i.e.} without considering the time variable). Remark that due to the support properties of the auxiliary function $\bar\varphi$ given in  (\ref{DefSoporteFuncTest2}) we have $supp(\pai \pj\bar \varphi) =Q_{R_3}\setminus Q_{\widetilde{R}}$ and recall by (\ref{DefSoporteFuncTest1}) we have $supp\; \bar \phi = Q_{\bar R}$ where $\bar{R}<\widetilde{R}<R_3$, thus by the properties of the kernel of the operator $\frac{\vn}{(-\Delta)}$ we can write
\begin{eqnarray}
\left| \bar\phi \frac{\vn}{(-\Delta )}(\pai \pj\bar \varphi) (u_i b_j)\right|&\leq& C\left|\int_{\mathbb{R}^3} \frac{1}{|x-y|^2}\mathds{1}_{Q_{\bar R}}(x)\mathds{1}_{Q_{R_3}\setminus Q_{\widetilde{R}}}(y)(\pai \pj\bar \varphi) (u_i b_j)(\cdot,y)dy\right|\notag\\
&\leq & C\left|\int_{\mathbb{R}^3} \frac{\mathds{1}_{|x-y|>(\widetilde{R}-\bar R)}}{|x-y|^2}\mathds{1}_{Q_{\bar R}}(x)\mathds{1}_{Q_{R_3}\setminus Q_{\widetilde{R}}}(y)(\pai \pj\bar \varphi) (u_i b_j)(\cdot,y)dy\right|,\label{Formula_intermediairevVV801}
\end{eqnarray}
and the previous expression is nothing but the convolution between the function $(\pai \pj\bar \varphi) (u_i b_j)$ and a $L^\infty$-function, thus we have 
\begin{equation}\label{Formula_intermediairevVV81}
\left\| \bar\phi \frac{\vn}{(-\Delta )}(\pai \pj\bar \varphi) (u_i b_j) (t,\cdot)\right\|_{L^\infty}\leq C\|(\pai \pj\bar \varphi) (u_i b_j)(t,\cdot)\|_{L^1}\leq C\|\mathds{1}_{Q_{R_3}}(u_i b_j)(t,\cdot)\|_{L^{\nu}},
\end{equation}
and taking the $L^\nu$-norm in the time variable we obtain
\begin{eqnarray*}
\left\| \bar\phi \frac{\vn}{(-\Delta )}(\pai \pj\bar \varphi) (u_i b_j)\right\|_{L^{\nu}_t L^{\infty}_{x}}
&\leq &C\|\mathds{1}_{Q_{R_3}}u_i b_j\|_{L^{\nu}_{t,x}}= C\|\mathds{1}_{Q_{R_3}}u_i b_j\|_{\mathcal{M}^{\nu,\nu}_{t,x}}\\
&\leq&C\|\mathds{1}_{Q_{R_3}}\vu\|_{\mathcal{M}^{3,\tau_0}_{t,x}}\|\mathds{1}_{Q_{R_3}}\vb\|_{\mathcal{M}^{3,\tau_0}_{t,x}}<+\infty,
\end{eqnarray*}
where we used the fact that $1<\nu<\frac{3}{2}<\frac{\tau_0}{2}$ and we applied Hölder's inequality. Gathering together all these estimates we obtain 
$\|\mathds{1}_{Q_{R_2}}\vVV_8\|_{\mathcal{M}^{3, \sigma}_{t,x}}<+\infty$.

\item The quantity $\vVV_9$ in (\ref{Def_VariablesVV1}) can be treated in the same way as the term $\vVV_8$. Indeed, by the same arguments displayed to deduce (\ref{Formula_intermediairevVV80}), we can write (recall that $1<\nu<\frac{3}{2}$ and thus we can replace $\nu$ by $q_0$ without loss of generality, see (\ref{HypothesesLocal1})): $\displaystyle{
\|\mathds{1}_{Q_{R_2}}\vVV_9\|_{\mathcal{M}^{3, \sigma}_{t,x}}
\leq C\left\| \bar\phi \frac{\vn}{(-\Delta )}( (\Delta \bar\varphi) P)\right\|_{L^{q_0}_t L^{\infty}_{x}}}$ and if we study the $L^\infty$-norm in the space variable of this term, by the same ideas used in (\ref{Formula_intermediairevVV801})-(\ref{Formula_intermediairevVV81}) we obtain
$$\left\| \bar\phi \frac{\vn}{(-\Delta )}( (\Delta \bar\varphi) P) (t,\cdot)\right\|_{L^{\infty}}\leq C\|(\Delta \bar\varphi) P (t,\cdot)\|_{L^1}\leq C\|\mathds{1}_{Q_{R_2}}P (t,\cdot)\|_{L^{q_0}}.$$
Thus, taking the $L^{q_0}$-norm in the time variable we have
$$\|\mathds{1}_{Q_{R_2}}\vVV_9\|_{\mathcal{M}^{3, \sigma}_{t,x}}
\leq C\left\| \bar\phi \frac{\vn}{(-\Delta )}( (\Delta \bar\varphi) P)\right\|_{L^{q_0}_t L^{\infty}_{x}}\leq C \|\mathds{1}_{Q_{R_3}}P\|_{L^{q_0}_{t,x}}<+\infty.$$
\item The study of the quantity $\vVV_{10}$ in (\ref{Def_VariablesVV1}) follows the same lines as the terms $\vVV_{8}$ and $\vVV_{9}$. However instead of (\ref{Formula_intermediairevVV801}) we have
\begin{equation}\label{Formula_intermediairevVV802}
\left|\bar\phi \frac{\vn\partial_i }{(-\Delta)} ( (\partial_i \bar\varphi) P )\right|
\leq  C\left|\int_{\mathbb{R}^3} \frac{\mathds{1}_{|x-y|>(\widetilde{R}-\bar R)}}{|x-y|^3}\mathds{1}_{Q_{\bar R}}(x)\mathds{1}_{Q_{R_3}\setminus Q_{\widetilde{R}}}(y)(\pai \bar \varphi) P(y)dy\right|,
\end{equation}
and thus we can write (again replacing $\nu$ by $q_0$)
$$\|\mathds{1}_{Q_{R_2}}\vVV_{10}\|_{\mathcal{M}^{3, \sigma}_{t,x}}\leq  \left\|\bar\phi \frac{\vn\partial_i }{(-\Delta)} ( (\partial_i \bar\varphi) P)\right\|_{L^{q_0}_t L^{\infty}_{x}}\leq C \|\mathds{1}_{Q_{R_3}}P\|_{L^{q_0}_{t,x}}<+\infty.$$
\item The last term of  (\ref{Def_VariablesVV1}) is easy to handle, indeed, we have 
\begin{eqnarray*}
|\mathds{1}_{Q_{R_2}}\vVV_{11}(t,x)|&\leq &\left|\mathds{1}_{Q_{R_2}}\int_{0}^{t} e^{ (t-s)\Delta} \bar\phi(\vf+\vg)ds\right|\leq C\mathds{1}_{Q_{R_2}}\int_{\mathbb{R}}\int_{\mathbb{R}^3}\frac{|\bar\phi(\vf+\vg)(s,y)|}{(|t-s|^{\frac{1}{2}}+|x-y|)^3}dyds\\
&\leq& C\mathds{1}_{Q_{R_2}}\mathcal{I}_{2}(\mathds{1}_{Q_{R_3}}|\vf+\vg|)(t,x),
\end{eqnarray*} 
and taking the Morrey $\mathcal{M}_{t,x}^{3, \sigma}$ norm we obtain
$\|\mathds{1}_{Q_{R_2}}\vVV_{11}\|_{\mathcal{M}_{t,x}^{3, \sigma}}\leq C\|\mathds{1}_{Q_{R_2}}\mathcal{I}_{2}(\mathds{1}_{Q_{R_3}}|\vf+\vg|)\|_{\mathcal{M}_{t,x}^{3, \sigma}}$, then if we set $\frac{11}{5}<q<\frac{5}{2}$ and $\lambda=1-\frac{2q}{5}$ we thus have $3\leq \frac{10}{7\lambda}$ and $\sigma<10<\frac{q}{\lambda}$. Now by Lemma \ref{Lemme_locindi} and Lemma \ref{Lemme_Hed} we have
$$\|\mathds{1}_{Q_{R_2}}\mathcal{I}_{2}(\mathds{1}_{Q_{R_3}}|\vf+\vg|)\|_{\mathcal{M}_{t,x}^{3, \sigma}}\leq C\|\mathcal{I}_{2}(\mathds{1}_{Q_{R_3}}|\vf+\vg|)\|_{\mathcal{M}_{t,x}^{\frac{10}{7\lambda}, \frac{q}{\lambda}}}\leq C\|\mathds{1}_{Q_{R_3}}|\vf+\vg|\|_{\mathcal{M}_{t,x}^{\frac{10}{7}, q}},$$
but since $q<\frac{5}{2}<\frac{5}{2-\alpha}<\tau_a,\tau_b$, by  Lemma \ref{Lemme_locindi} we obtain 
$$\|\mathds{1}_{Q_{R_3}}|\vf+\vg|\|_{\mathcal{M}_{t,x}^{\frac{10}{7}, q}}\leq C\left(\|\mathds{1}_{Q_{R_3}}\vf\|_{\mathcal{M}_{t,x}^{\frac{10}{7}, \tau_a}}+\|\mathds{1}_{Q_{R_3}}\vg\|_{\mathcal{M}_{t,x}^{\frac{10}{7}, \tau_b}}\right)<+\infty,$$
thus, gathering all the estimates above we have $\|\mathds{1}_{Q_{R_2}}\vVV_{11}\|_{\mathcal{M}_{t,x}^{3, \sigma}}<+\infty$.
\end{itemize}
Lemma \ref{Lemme_EstimationMorrey3Sigma} is now completely proven. \hfill $\blacksquare$\\

\noindent{\bf End of the proof of Proposition \ref{Propo_Points24}.} With this Lemma \ref{Lemme_EstimationMorrey3Sigma} we have proven so far that $\mathds{1}_{Q_{R_2}} \vu \in \mathcal{M}_{t,x}^{3, \sigma}$ and $\mathds{1}_{Q_{R_2}} \vb \in \mathcal{M}_{t,x}^{3, \sigma}$ for $\tau_0<\sigma$ with $\sigma$ very close to $\tau_0$, say $\sigma=\tau_0+\epsilon$ (see Remark \ref{Remarque_ConditionSigma})). But this is not enough to ensure the hypothesis \eqref{Estimations24}, \emph{i.e.} the condition $\frac{1}{\delta}+\frac{1}{\tau_{0}} < \frac{1-\alpha}{5}$ stated in Proposition \ref{Propo_Points24}.\\ 

In order to obtain this relationship, we will iterate the arguments above. Indeed, plugging the information $\mathds{1}_{Q_{R_2}} \vu \in \mathcal{M}_{t,x}^{3, \tau_0+\epsilon}$ and $\mathds{1}_{Q_{R_2}} \vb \in \mathcal{M}_{t,x}^{3, \tau_0+\epsilon}$ in the set of hypotheses (\ref{Estimations135}) and reapplying Lemma \ref{Lemme_EstimationMorrey3Sigma}, we will obtain $\mathds{1}_{Q_{\bar R_2}} \vu \in \mathcal{M}_{t,x}^{3, \sigma_1}$ and $\mathds{1}_{Q_{\bar R_2}} \vb \in \mathcal{M}_{t,x}^{3, \sigma_1}$ where $\bar R_2<R_2$ and $\sigma_1=\sigma+\epsilon=\tau_0+2\epsilon$. Repeating these arguments until obtaining $\mathds{1}_{Q_{\bar {\bar{R_2}}}} \vu \in \mathcal{M}_{t,x}^{3, \sigma_n}$ and $\mathds{1}_{Q_{\bar{\bar{R_2}}}} \vb \in \mathcal{M}_{t,x}^{3, \sigma_n}$ where $\sigma_n=\tau_0+(n+1) \epsilon$ such that $\frac{1}{\sigma_n}+\frac{1}{\tau_0}<\frac{1-\alpha}{5}$ with $\bar {\bar{R_2}}<\bar{R_2}$. As we can see, at each iteration we have to consider smaller parabolic balls and without fear of confusion we can set $\delta=\sigma_n$ with the corresponding radius to be $R_2$. We thus have 
$\mathds{1}_{Q_{ R_2}} \vu \in \mathcal{M}_{t,x}^{3, \delta}$ and $\mathds{1}_{Q_{ R_2}} \vb \in \mathcal{M}_{t,x}^{3, \delta}$ with $\frac{1}{\delta}+\frac{1}{\tau_0}<\frac{1-\alpha}{5}$ and the proof of  Proposition \ref{Propo_Points24} is finished. \hfill $\blacksquare$\\
\begin{Remarque}\label{RemarkDeltaGrand}
Note that the parameter $\delta$ can be made big enough in order to satisfy $\delta\geq 10$.
\end{Remarque}
We study now the information stated in (\ref{Estimations241}) and we have the following corollary. 
\begin{Corollaire}\label{CorollaireFinale}
Let $(\vu, p, \vb)$ be a suitable solution of MHD equations (\ref{EquationMHD}) over $\Omega$ in the sense of Definition \ref{Def_SuitableSolutions}. Assume the general hypotheses (\ref{HypothesesLocal1}) and assume moreover the local informations (\ref{Estimations135}) for a parabolic ball $Q_{R_3}$ and (\ref{Estimations24}) for a parabolic ball $Q_{R_2}$. Then for some $R_1$ such that $R_0<R_1<R_2<R_3$ and for all $1\leq i,j\leq 3$ we have 
$$\mathds{1}_{Q_{R_1}} \frac{\vn\partial_i \partial_j}{(-\Delta)}(u_ib_j)\in\mathcal{M}_{t,x}^{\mathfrak{p}, \mathfrak{q}},$$
with $\mathfrak{p}_0\leq \mathfrak{p}<+\infty$, $\mathfrak{q}_0\leq \mathfrak{q}<+\infty$ and where $1 \leq \mathfrak{p}_0\leq \tfrac{6}{5}$ and $\tfrac52 <\mathfrak{q}_0< 3$ where $\frac{1}{\mathfrak{q}_0}=\frac{2-\alpha}{5}$ with $0< \alpha <  \frac{1}{3}$.
\end{Corollaire}
{\bf Proof.} We consider here two auxiliary functions $\widetilde{\phi}$ and $\widetilde{\varphi}$ satisfying the same properties stated in (\ref{DefSoporteFuncTest1}) and $(\ref{DefSoporteFuncTest2})$ where we replace $R_3$ by $R_2$ and $R_2$ by $R_1$, respectively. Thus, by definition of the auxiliary function $\widetilde{\phi}$ we have the identity $\mathds{1}_{Q_{R_1}}=\widetilde{\phi} \mathds{1}_{Q_{R_1}}$ and then we can study the term $\widetilde{\phi} \vn P=\widetilde{\phi}  \displaystyle{\sum^3_{i,j= 1}}\frac{\vn}{(- \Delta )} \partial_i  \partial_j  (u_i b_j )$,  but due to the computations performed in (\ref{Definition_LocalPression})-(\ref{Identite_Pression1}) we have the identity
\begin{equation}\label{DerniereFormule}
\begin{split}
\widetilde{\phi} \vn P=&\underbrace{\sum^3_{i,j= 1} \widetilde{\phi} \frac{\vn \partial_i  \partial_j }{(- \Delta )} (\widetilde{\varphi} u_i b_j )}_{(a)}-\underbrace{\sum^3_{i,j= 1}  \frac{\widetilde{\phi} \vn\pai}{(- \Delta )} (\pj \widetilde{\varphi}) u_i b_j}_{(b)}-\underbrace{\sum^3_{i,j= 1}  \frac{\widetilde{\phi} \vn\partial_j}{(- \Delta )} (\partial_i \widetilde{\varphi}) u_i b_j}_{(c)}\\
&+ 2 \underbrace{\sum^3_{i,j= 1} \widetilde{\phi} \frac{\vn}{(-\Delta )}(\pai \pj\widetilde{\varphi}) (u_i b_j)}_{(d)}+ \underbrace{\widetilde{\phi} \frac{\vn \big(  (\Delta \widetilde{\varphi}) P \big) }{(-\Delta)}}_{(e)}
-2 \underbrace{\sum^3_{i= 1}\widetilde{\phi}\frac{\vn\big(\partial_i ( (\partial_i \widetilde{\varphi}) P )\big)}{(-\Delta)}}_{(f)}
\end{split}
\end{equation}
and we only need to prove that each one of these terms belong to the Morrey space $\mathcal{M}^{\frac{6}{5}, 3}_{t,x}$.\\

\begin{itemize}
\item[$\bullet$] The term $(a)$ is treated as follows: since the Riesz transforms are bounded in Morrey spaces we obtain
$$\left\|\widetilde{\phi}\frac{\vn  \pai \pj }{(- \Delta )}( \widetilde{\varphi} u_i b_j )\right\|_{\mathcal{M}^{\frac{6}{5}, 3}_{t,x}}\leq C\left\|\vn( \widetilde{\varphi} u_i b_j )\right\|_{\mathcal{M}^{\frac{6}{5}, 3}_{t,x}},$$
now, for $1\leq k\leq 3$, by Remark \ref{RemarkDeltaGrand}, using all the information available and by H\"older's inequality in Morrey spaces (recall that $0<\alpha<\frac{1}{3}$ and $\frac{1}{\delta}+\frac{1}{\tau_0}<\frac{1-\alpha}{5}$), we have
\begin{eqnarray*}
\left\|(\partial_k\widetilde{\varphi}) u_i b_j \right\|_{\mathcal{M}^{\frac{6}{5}, 5}_{t,x}}&\leq& C\left\|\mathds{1}_{Q_{R_2}}u_i b_j \right\|_{\mathcal{M}^{\frac{3}{2}, 5}_{t,x}}\leq C\|\mathds{1}_{Q_{R_2}}u_i\|_{\mathcal{M}^{3, \delta}_{t,x}}\|\mathds{1}_{Q_{R_2}}b_j\|_{\mathcal{M}^{3, \delta}_{t,x}}<+\infty\\
\|\widetilde{\varphi}(\partial_k u_i) b_j \|_{\mathcal{M}^{\frac{6}{5}, 3}_{t,x}}&\leq &C\|\mathds{1}_{Q_{R_3}}\vn \otimes \vu\|_{\mathcal{M}^{2, \tau_1}_{t,x}}\|\mathds{1}_{Q_{R_3}} u_j\|_{\mathcal{M}^{3, \delta}_{t,x}}<+\infty\\
 \|\widetilde{\varphi} u_i(\partial_k b_j) \|_{\mathcal{M}^{\frac{6}{5}, 3}_{t,x}}&\leq &C\|\mathds{1}_{Q_{R_3}}u_i\|_{\mathcal{M}^{3, \delta}_{t,x}}\|\mathds{1}_{Q_{R_3}} \vn \otimes \vb\|_{\mathcal{M}^{2, \tau_1}_{t,x}}<+\infty,
\end{eqnarray*}
thus we can deduce that we have the estimate 
$$\left\|\widetilde{\phi}\frac{\vn  \pai \pj }{(- \Delta )}( \widetilde{\varphi} u_i b_j )\right\|_{\mathcal{M}^{\frac{6}{5}, 3}_{t,x}}<+\infty.$$

\item[$\bullet$] The terms $(b)$ and $(c)$ of (\ref{DerniereFormule}) can be treated in a similar fashion and we have:
\begin{eqnarray*}
\left\|\frac{\widetilde{\phi} \vn\pai}{(- \Delta )} (\pj \bar\varphi) u_i b_j\right\|_{\mathcal{M}^{\frac{6}{5}, 3}_{t,x}}&\leq &C\|\mathds{1}_{Q_{R_2}}u_i b_j\|_{\mathcal{M}^{\frac{6}{5}, 5}_{t,x}}\leq C \|\mathds{1}_{Q_{R_2}}u_i b_j\|_{\mathcal{M}^{\frac{3}{2}, 5}_{t,x}}\\
&\leq &C\|\mathds{1}_{Q_{R_2}}u_i\|_{\mathcal{M}^{3, 10}_{t,x}}\|\mathds{1}_{Q_{R_2}}b_j\|_{\mathcal{M}^{3, 10}_{t,x}}\leq C\|\mathds{1}_{Q_{R_2}}u_i\|_{\mathcal{M}^{3, \delta}_{t,x}}\|\mathds{1}_{Q_{R_2}}b_j\|_{\mathcal{M}^{3, \delta}_{t,x}}<+\infty.
\end{eqnarray*}
\item[$\bullet$] The term $(d)$ is treated as follows.
$$\left\|\widetilde{\phi} \frac{\vn}{(-\Delta )}(\pai \pj \widetilde{\varphi}) (u_i b_j)\right\|_{\mathcal{M}^{\frac{6}{5}, 3}_{t,x}}\leq C\left\|\widetilde{\phi} \frac{\vn}{(-\Delta )}(\pai \pj \widetilde{\varphi}) (u_i b_j)\right\|_{\mathcal{M}^{\frac{3}{2}, \frac{15}{4}}_{t,x}}\leq C\left\|\widetilde{\phi} \frac{\vn}{(-\Delta )}(\pai \pj  \widetilde{\varphi}) (u_i b_j)\right\|_{L^{\frac{3}{2}}_t L^{\infty}_{x}},$$
where we used the space inclusion $L^{\frac{3}{2}}_t L^{\infty}_{x}\subset \mathcal{M}^{\frac{3}{2}, \frac{15}{4}}_{t,x}$. Following the same ideas displayed in formulas (\ref{Formula_intermediairevVV80})-(\ref{Formula_intermediairevVV81}), due to the support properties of the auxiliary functions we obtain
$$\left\|\widetilde{\phi} \frac{\vn}{(-\Delta )}(\pai \pj \widetilde{\varphi}) (u_i b_j)\right\|_{L^{\frac{3}{2}}_t L^{\infty}_{x}}\leq \|\mathds{1}_{Q_{R_3}}u_i b_j\|_{L^{\frac{3}{2}}_{t,x}}\leq C\|\mathds{1}_{Q_{R_3}}\vu\|_{\mathcal{M}^{3,\tau_0}_{t,x}}\|\mathds{1}_{Q_{R_3}}\vb\|_{\mathcal{M}^{3,\tau_0}_{t,x}}<+\infty.$$
\item[$\bullet$] The term $(e)$ of (\ref{DerniereFormule}) follows the same ideas as previous one, and we have
$$\left\|\widetilde{\phi} \frac{\vn \big(  (\Delta \widetilde{\varphi}) P \big) }{(-\Delta)}\right\|_{\mathcal{M}^{\frac{6}{5}, 3}_{t,x}}\leq C\left\|\widetilde{\phi} \frac{\vn \big(  (\Delta \widetilde{\varphi}) P \big) }{(-\Delta)}\right\|_{L^{\frac{3}{2}}_t L^{\infty}_{x}}\leq C \|\mathds{1}_{Q_{R_3}}P\|_{L^{\frac{3}{2}}_{t,x}}<+\infty.$$
\item[$\bullet$] The last term of (\ref{DerniereFormule}) is estimated in a very similar manner (see also (\ref{Formula_intermediairevVV802})):
$$\left\|\widetilde{\phi}\frac{\vn\big(\partial_i ( (\partial_i \widetilde{\varphi}) P )\big)}{(-\Delta)}\right\|_{\mathcal{M}^{\frac{6}{5}, 3}_{t,x}}\leq C\left\|\widetilde{\phi}\frac{\vn\big(\partial_i ( (\partial_i \widetilde{\varphi}) P )\big)}{(-\Delta)}\right\|_{L^{\frac{3}{2}}_t L^{\infty}_{x}}\leq C \|\mathds{1}_{Q_{R_3}}P\|_{L^{\frac{3}{2}}_{t,x}}<+\infty.$$
\end{itemize}
The proof of Corollary \ref{CorollaireFinale} is finished. 
\hfill $\blacksquare$\\

From the set of hypotheses of Theorem \ref{Teorem1} we have now deduced the framework used to prove Proposition \ref{Proposition_Principale1} and thus the proof of Theorem \ref{Teorem1} is finished.
\appendix 
\section{Useful Properties of Morrey spaces}\label{Secc_AppendixA}

\begin{Lemme}\label{Lemme_Product}
\begin{itemize}
\item[]
\item[1)]If $\vf, \vg:\mathbb{R} \times \R\longrightarrow \R$ are two functions such that $\vf\in \mathcal{M}_{t,x}^{p,q} (\mathbb{R} \times \R)$ and $\vg\in L^{\infty}_{t,x} (\mathbb{R} \times \R)$, then for all $1\leq p\leq q<+\infty$ we have
$$\|\vf\cdot\vg\|_{\mathcal{M}_{t,x}^{p,q}}\leq  C\|\vf\|_{\mathcal{M}_{t,x}^{p, q}} \|\vg\|_{L^{\infty}_{t,x}}.$$
\item[2)]If $\vf, \vg:\mathbb{R} \times \R\longrightarrow \R$ are two functions that belong to the space $\mathcal{M}_{t,x}^{p,q} (\mathbb{R} \times \R)$ then we have the inequality
$$\|\vf\cdot\vg\|_{\mathcal{M}_{t,x}^{\frac{p}{2}, \frac{q}{2}}}\leq  C\|\vf\|_{\mathcal{M}_{t,x}^{p, q}} \|\vg\|_{\mathcal{M}_{t,x}^{p, q}}.$$
\item[3)] More generally, let $1\leq p_0 \leq q_0 <+\infty$, $1\leq p_1\leq q_1<+\infty$ and $1\leq p_2\leq q_2<+\infty$. If $\tfrac{1}{p_1}+\tfrac{1}{p_2}\leq \frac{1}{p_0}$ and $\tfrac{1}{q_1}+\tfrac{1}{q_2}=\tfrac{1}{q_0}$, then for two measurable functions $\vf, \vg:\mathbb{R} \times \R\longrightarrow \R$ such that $\vf\in \mathcal{M}^{p_1, q_1}_{t,x}$ and $\vg\in \mathcal{M}^{p_2, q_2}_{t,x}$, we have the following version of the H\"older inequality in Morrey spaces:
$$\|\vf\cdot \vg\|_{\mathcal{M}^{p_0, q_0}_{t,x}}\leq \|\vf\|_{\mathcal{M}^{p_1, q_1}_{t,x}}\|\vg\|_{\mathcal{M}^{p_2, q_2}_{t,x}}.$$
\end{itemize}
\end{Lemme} 
Our next lemma (which is a particular case of the previous lemma) explains the behaviour of parabolic Morrey spaces with respect to localization in time and space. 
\begin{Lemme}\label{Lemme_locindi}
Let $\Omega$ be a bounded set of $\mathbb{R} \times \R$. If we have $1\leq p_0 \leq p_1$, $1\leq p_0\leq q_0 \leq q_1<+\infty$ and if the function $\vf:\mathbb{R} \times \R\longrightarrow \R$  belongs to the space $\mathcal{M}_{t,x}^{p_1,q_1} (\mathbb{R} \times \R)$ then we have the following localization property 
$$\|\mathds{1}_{\Omega}\vf\|_{\mathcal{M}_{t,x}^{p_0, q_0}} \leq C\|\mathds{1}_{\Omega}\vf\|_{\mathcal{M}_{t,x}^{p_1,q_1}}\leq C\|\vf\|_{\mathcal{M}_{t,x}^{p_1,q_1}}.$$
\end{Lemme} 

\section{A technical Lemma}\label{Secc_AppendixB}
{\bf Proof of Lemma \ref{Lemme_PropertiesTestFunctionEnergyInequality}.}
\begin{itemize}
\item[1)] This point holds true thanks to the properties of the test function $\phi$ (see (\ref{Definition_FuncTestInegaliteEnergie})) and the properties of the heat kernel $\mathfrak{g}_t(x)$. Indeed, for all $(s,y)\in Q_{r}(t,x) $ we have $3r^2 < 4r^2+t-s < 5r^2$ and $|x-y| < r$ and thus we obtain:
$$\mathfrak{g}_{(4r^2+t-s)}(x-y) = \frac{1}{(4 \pi (4r^2+t-s))^{\frac{3}{2}}} e^{- \frac{|x-y|^2}{4(4r^2+t-s)}}\geq \frac{C}{r^3}. $$ 
Thus, estimate the estimate $\omega(s,y)\geq \frac{C}{r}$ holds due to the definition of the auxiliary functions $\phi$ and $\theta$.
\item[2)] For the second point, for $s< t+r^2$, by the usual heat kernel estimates, we have
\begin{equation}\label{heatkernel}
\mathfrak{g}_{(4r^2+t-s)}(x-y) \leq   \frac{C}{(4r^2+t-s)^{\frac{3}{2}} + |x-y|^{3} } \leq \frac{C}{r^3},
\end{equation}
hence, the estimate $\omega(s,y)\leq \frac{C}{r}$ is valid for all $(s, y) \in Q_{\rho}(t,x)$. 
\item[3)] Note that we have 
\begin{equation}\label{heatkernelderi}
|\vn \mathfrak{g}_{(4r^2+t-s)}(x-y)| \leq   \frac{C}{(4r^2+t-s)^{2} + |x-y|^{4} }\leq \frac{C}{r^4} \quad \text{for} \quad s< t+r^2.
\end{equation}
Since $\nabla \phi$ is supported outside the cylinder $Q_{\frac{1}{2}}$, we shall only consider the case $(s,y)\in  Q_{\rho}(t,x)\setminus Q_{\frac{\rho}{2}}(t,x)$. Using \eqref{heatkernel} again, we find the following estimate
\begin{equation}\label{heatkernel1}
\mathfrak{g}_{(4r^2+t-s)}(x-y) \leq \frac{C}{\rho^3} \leq \frac{C}{r^3} \quad \text{for} \quad (s,y)\in  Q_{\rho}(t,x)\setminus Q_{\frac{\rho}{2}}(t,x).
\end{equation}
This estimate and \eqref{heatkernelderi} imply by contruction (recall (\ref{Definition_FuncTestInegaliteEnergie})) that $|\vn \omega(s,y)|\leq \frac{C}{r^2}.$
\item[4)] Regarding the last estimate, we first note that $(\partial_s + \Delta_y) \mathfrak{g}_{(4r^2+t-s)}(x-y) = 0$, so it remains to treat the term involving $\vn \mathfrak{g}$ and the  case when time derivative and space derivative fall on the two test function $\phi$ and $\theta$. For the time derivative, we see that $\partial_s \left(\theta\left(\frac{s-t}{r^2}\right)\right)$ is neglected for all $s < t+ r^2$. For space derivative, we have 
$$|\vn \mathfrak{g}_{(4r^2+t-s)}(x-y)| \leq \frac{C}{\rho^4}, \quad \text{for} \quad (s,y)\in  Q_{\rho}(t,x)\setminus Q_{\frac{\rho}{2}}(t,x).$$
For the same reason as before, since $\nabla \phi $ vanishes on $Q_{\frac{1}{2}}$, the estimate $|(\partial_s + \Delta) \omega(s,y)|\leq C \frac{r^2}{\rho^5}$ follows from the estimate above and \eqref{heatkernel1}.\hfill$\blacksquare$\\
\end{itemize}

\noindent{\bf Acknowledgments:} We would like to thank Professor Pierre-Gilles Lemari\'e-Rieusset for fruitful discussions.  J. \textsc{He} is supported by the program \emph{Sophie Germain} of the \emph{Fondation Math\'ematique Jacques Hadamard}.

\end{document}